\documentclass[reqno]{amsart}

\usepackage{amsmath,amsfonts,amscd,amsthm,graphicx,amssymb,epic,eepic}

\theoremstyle{plain}
\newtheorem{theorem}{Theorem}
\newtheorem{proposition}{Proposition}

\newtheorem{lemma}{Lemma}
\numberwithin{equation}{section}

\newcommand{\eps}{\varepsilon}
\newcommand{\Z}{{\mathbb{Z}}}

\renewcommand{\P}{{\mathbb{P}}}
\newcommand{\R}{{\mathbb{R}}}

\newcommand{\T}{{\mathbb{T}}}
\newcommand{\N}{{\mathbb{N}}}

\newcommand{\TT}{{\mathcal{T}}}
\newcommand{\FQ}{{\mathcal{F}_Q}}
\newcommand{\FI}{{\mathcal{F}_Q (I)}}
\newcommand{\II}{{\mathcal{I}}}
\newcommand{\JJ}{{\mathcal{J}}}
\renewcommand{\SS}{{\mathfrak{S}}}
\newcommand{\ar}{\operatorname{area}}

\newcommand{\dist}{\operatorname{dist}}
\renewcommand{\t}{\operatorname{t}}
\newcommand{\G}{\mathbb{G}}

\title[The Periodic Lorentz Gas in the
Small-Scatterer Limit]{The Distribution of the Free Path Lengths
in the Periodic Two-Dimensional Lorentz Gas in the Small-Scatterer
Limit}

\author{Florin P. Boca and Alexandru Zaharescu}

\address{Department of Mathematics, University of Illinois at
Urbana-Champaign, 1409 W. Green Street, Urbana, IL 61801, USA}

\address{Institute of Mathematics ``Simion Stoilow" of the Romanian
Academy, P.O. Box 1-764, RO-014700 Bucharest, Romania}

\email{fboca@math.uiuc.edu}

\email{zaharesc@math.uiuc.edu}


\begin{document}

\maketitle

\begin{abstract}
We study the free path length and the geometric free path length
in the model of the periodic two-dimensional Lorentz gas (Sinai
billiard). We give a complete and rigorous proof for the existence
of their distributions in the small-scatterer limit and explicitly
compute them. As a corollary one gets a complete proof for the
existence of the constant term $c=2-3\ln
2+\frac{27\zeta(3)}{2\pi^2}$ in the asymptotic formula $h(T)=-2\ln
\eps+c+o(1)$ of the KS entropy of the billiard map in this model,
as conjectured by P. Dahlqvist.
\end{abstract}

\section{Introduction and main results}
A periodic two-dimensional Lorentz gas (Sinai billiard) is a
billiard system on the two-dimensional torus with one or more
circular regions (scatterers) removed. This model in classical
mechanics was introduced by Lorentz \cite{Lor} in 1905 to describe
the dynamics of electrons in metals. The associated dynamical
system is simple enough to allow a comprehensive study, yet
complex enough to exhibit chaos. According to Gutzwiller
\cite{Gu}: ``The original billiard of Sinai was designed to
imitate, in the most simple-minded manner, a gas of hard spherical
balls which bounce around inside a finite enclosure. The
formidable technical difficulties of this fundamental problem were
boiled down to the shape of a square for the enclosure, and the
collisions between the balls were reduced to a single point
particle hitting a circular hard wall at the center of the
enclosure.''
\par
The model was intensively studied from the point of view of
dynamical systems \cite{Bu,Ch2,Dah,FOK,Gal,Gol,Sin}. Our primary
goal here is to estimate the \emph{free-path length} (\emph{first
return time}) in this periodic two-dimensional model in the
small-scatterer limit. We solve the following three open problems:
\begin{itemize}
\item[(1)] the existence and computation of the distribution of
the free path length, previously considered in \cite{BGW,CG,DDG}.
\item[(2)] the existence and computation of the distribution of
the geometric free path length, previously shown, but not fully
proved, in \cite{Dah}. \item[(3)] the existence and computation of
the second (constant) term in the asymptotic formula of the KS
entropy $h(T_\eps)$ of the billiard map in this model, previously
studied in \cite{Ch1,Ch2,Dah,FOK}.
\end{itemize}
\par
For each $\eps \in (0,\frac{1}{2})$ let
\begin{equation*}
Z_\eps =\{ x\in \R^2 \, ;\, \dist (x,\Z^2)\geq \eps \},
\end{equation*}
denote by $\partial Z_\eps$ the boundary $\Z^2+\eps \T$ of
$Z_\eps$, and define the \emph{free path length} (also called
\emph{first exit time}) as the Borel map given by
\begin{equation*}
\tau_\eps (x,\omega)=\inf \{ \tau >0 \, ;\, x+\tau \omega \in
\partial Z_\eps \} ,
\qquad x\in Z_\eps, \ \omega \in \T .
\end{equation*}
If $\tan \omega$ is irrational, then $\tau_\eps (x,\omega)<\infty$
for every $x\in Z_\eps$. We consider the probability space
$(Y_\eps,\mu_\eps)$, with $Y_\eps=Z_\eps /\Z^2 \subseteq [0,1)^2$
and $\mu_\eps$ the normalized Lebesgue measure on $Y_\eps$. Let
$e_t=e_{(t,\infty)}$ denote the characteristic function of
$(t,\infty)$. For every $t>0$ the probability that $\tau_\eps
(x,\omega)>\frac{t}{2\eps}$ is given by
\begin{equation*}
\P_\eps (t) =\mu_\eps ( \{ (x,\omega) \in Y_\eps \times [0,2\pi)
\, ;\, 2\eps \tau_\eps (x,\omega)> t \}) =\int_{Y_\eps \times \T}
e_t (2\eps \tau_\eps) \, d\mu_\eps .
\end{equation*}
Lower and upper bounds for $\P_\eps$ of correct order of magnitude
were established by Bourgain, Golse and Wennberg \cite{BGW}, using
the rational channels introduced by Bleher \cite{Ble}. More
recently, Caglioti and Golse \cite{CG} have proved the existence
of the Cesaro $\limsup$ and $\liminf$ means, proving for large $t$
that
\begin{equation}\label{1.1}
\begin{split}
\limsup\limits_{\delta \rightarrow 0^+} \frac{1}{\vert \ln \delta
\vert} \int_\delta^{1/4} \P_\eps (t)\, \frac{d\eps}{\eps} &
=\frac{2}{\pi^2 t} +O\left( \frac{1}{t^2} \right)=
\liminf\limits_{\delta \rightarrow 0^+} \frac{1}{\vert \ln \delta
\vert} \int_\delta^{1/4} \P_\eps (t)\, \frac{d\eps}{\eps} .
\end{split}
\end{equation}
\par
In Sections 2-7 below we prove the existence of the limit $\P(t)$
of $\P_\eps (t)$ as $\eps \rightarrow 0^+$ and explicitly compute
it.

\begin{theorem}\label{T1.1}
For every $t>0$ and $\delta >0$
\begin{equation*}
\P_\eps(t)=\P(t)+O_\delta (\eps^{1/8-\delta})\qquad
(\varepsilon\rightarrow 0^+),
\end{equation*}
with
\begin{equation*}
\begin{split}
\P (t) & =\frac{6}{\pi^2}\begin{cases}\vspace{0.1cm}
\mbox{$\displaystyle \frac{\pi^2}{6} (1-t)+\frac{t^2}{2}$}
 & \mbox{if $0<t\leq 1;$} \\ \vspace{0.1cm}
\mbox{$\displaystyle \int_0^{t-1} \psi(x,t)\, dx + \int_{t-1}^1
\phi (x,t)\, dx $} &
\mbox{if $1<t\leq 2;$} \\
\mbox{$\displaystyle \int_0^1 \psi(x,t)\, dx$} & \mbox{if $ t>2,$}
\end{cases} \\
\psi (x,t) & =\frac{(1-x)^2}{x} \left( 2\ln \frac{t-x}{t-2x}
-\frac{t}{x} \, \ln \frac{(t-x)^2}{t(t-2x)} \right)
, \\
\phi(x,t) & =\frac{1-t}{x}\,\ln \frac{1}{t-x} +
\frac{(t-x)(x-t+1)}{x} +\frac{(1-x)^2}{x} \left( 2\ln
\frac{t-x}{1-x}-\frac{t}{x}\, \ln \frac{t-x}{t(1-x)} \right) .
\end{split}
\end{equation*}
\end{theorem}
After a direct computation the above formula for $\P (t)$ yields
\begin{equation*}
\P (t)=\frac{24}{\pi^2} \sum\limits_{n=1}^\infty
\frac{2^n-1}{n^2(n+1)^2(n+2) t^n} ,\qquad t\geq 2,
\end{equation*}
and thus for large $t$ we find
\begin{equation*}
\P (t)=\frac{2}{\pi^2 t}+O\left( \frac{1}{t^2} \right),
\end{equation*}
which agrees with \eqref{1.1}.
\par
The related ``homogeneous" problem when the trajectory starts at
the origin $O$ and the phase space is a subinterval of the
velocity range $[0,2\pi)$ was studied by Gologan and the authors.
The limit distribution
\begin{equation*}
H(t)=\lim\limits_{\eps \rightarrow 0^+} \frac{1}{2\pi} \big\vert
\{ \omega \in [0,2\pi) \, ;\, \eps \tau_\eps (O,\omega) >t \}
\big\vert =\lim\limits_{\eps \rightarrow 0^+} \int_{\T} e_t \big(
\eps \tau_\eps (O,\omega)\big) \, d\omega,
\end{equation*}
where $\vert \ \vert$ denotes the Lebesgue measure, was shown to
exist and explicitly computed in \cite{BGZ0,BGZ}. Unlike $\P$, the
function $H$ is compactly supported on the interval $[0,1]$.
Interestingly, in the particular situation where the scatterers
are vertical segments, this case is related to some old problems
in diophantine approximation investigated by Erd\" os, Sz\" usz
and Tur\' an \cite{ERD,EST}, Friedman and Niven \cite{FN}, and by
Kesten \cite{KES}.
\par
The main tools used to prove Theorem \ref{T1.1} are a certain
three-strip partition of $[0,1)^2$ and the Weil-Sali\' e estimate
for Kloosterman sums \cite{Est,Hoo,Wei}. The latter is used in
infinitesimal form with respect to the parameter $\omega$ to count
the number of solutions of equations of form $xy=1\pmod{q}$ in
various regions in $\R^2$. This approach, somehow reminiscent of
the circle method, produces good estimates, allowing us to keep
under control the error terms. It was developed and used recently
in many situations to study problems related to the spacing
statistics of Farey fractions and lattice points in $\R^2$
\cite{ABCZ,BCZ,BCZ1,BGZ0,BGZ}. A possible source for getting
better estimates for the error terms might come from further
cancellations in certain sums of Kloosterman sums, of the form
\cite{Des,GS,Kuz}
\begin{equation*}
S=\sum\limits_{a,b} \sum\limits_c h_{a,b}(c) S(a,b;c).
\end{equation*}
The three-strip partition of $\T^2$ is related to the continued
fraction decomposition of the slope of the trajectory. Following
work of Blank and Krikorian \cite{BK} on the longest orbit of the
billiard, Caglioti and Golse explicitly introduced this partition
and used it in conjunction with ergodic properties of the Gauss
map \cite{CG} to prove \eqref{1.1}. We will use it in Section 3 in
a suitable setting for our computations.
\par
One can also consider the phase space $\Sigma_\eps^+=\{ (x,\omega)
\in\partial Y_\eps \times \T \, ;\, \omega \cdot n_x>0\}$ with
$n_x$ the inward unit normal at $x\in \partial Y_\eps$ and the
probability measure $\nu_\eps$ on $\Sigma_\eps^+$ obtained by
normalizing the Liouville measure $\omega \cdot n_x \, dx\,
d\omega$ to mass one. Consider also the distribution
\begin{equation*}
\G_\eps (t) =\nu_\eps (\{ (x,\omega) \in \Sigma_\eps^+ \, ;\,
2\eps \tau_\eps (x,\omega)  >t \}) =\int_{\Sigma_\eps^+} e_t
(2\eps \tau_\eps)\, d\nu_\eps
\end{equation*}
of the \emph{geometric free path length} $\tau_\eps(x,\omega)$.
The first moment (\emph{geometric mean free path length}) of
$\tau_\eps$ with respect to $\nu_\eps$ can be expressed as
\begin{equation}\label{1.2}
\int_{\Sigma_\eps^+} \tau_\eps \, d\nu_\eps =\frac{\pi \vert
Y_\eps\vert}{\vert \partial Y_\eps \vert}=\frac{1-\pi
\eps^2}{2\eps}  .
\end{equation}
Equality \eqref{1.2} is a consequence of a more general formula of
Santal\' o \cite{San} who extended earlier work of P\' olya on the
mean visible distance in a forrest \cite{Po}. The formulation from
\eqref{1.2} appears in \cite{Ch1,Ch2,DDG}. Knowledge of the mean
free path does not give however any information on other moments
or on the limiting distribution of the free path in the
small-scatterer limit. Our number theoretical analysis leads to
the following solution of this limiting distribution problem,
proved in Sections 8-11 below.

\begin{theorem}\label{T1.2}
For every $t>0$ and $\delta >0$
\begin{equation*}
\G_\eps (t)=\G(t)+O_\delta (\eps^{1/8-\delta}) \qquad
(\varepsilon\rightarrow 0^+),
\end{equation*}
with
\begin{equation*}
\begin{split}
\G(t) & =\frac{6}{\pi^2}\begin{cases} \vspace{0.1cm}
\mbox{$\displaystyle \frac{\pi^2}{6}-t$}
 & \mbox{if $0<t\leq 1;$} \\
\vspace{0.1cm} \mbox{$\displaystyle \left(
-2+t+(t-1)\ln\frac{1}{t-1}\right) +\int_0^{t-1}
\widetilde{\psi}(x,t)\, dx+\int_{t-1}^1 \widetilde{\phi} (x,t)\,
dx$}  &
\mbox{if $1<t\leq 2;$} \\
\mbox{$\displaystyle \int_0^1 \widetilde{\psi}(x,t)\, dx$} &
\mbox{if $t>2,$} \end{cases}
\\
\widetilde{\psi}(x,t) & =\frac{(1-x)^2}{x^2}\ln
\frac{(t-x)^2}{t(t-2x)} ,\quad
\widetilde{\phi}(x,t)=\frac{1}{x}\ln \frac{1}{t-x}
+\frac{(1-x)^2}{x^2}\ln \frac{t-x}{t(1-x)} .
\end{split}
\end{equation*}
\end{theorem}

\begin{figure}[ht]
\includegraphics*[scale=0.57, bb=90 0 280 150]{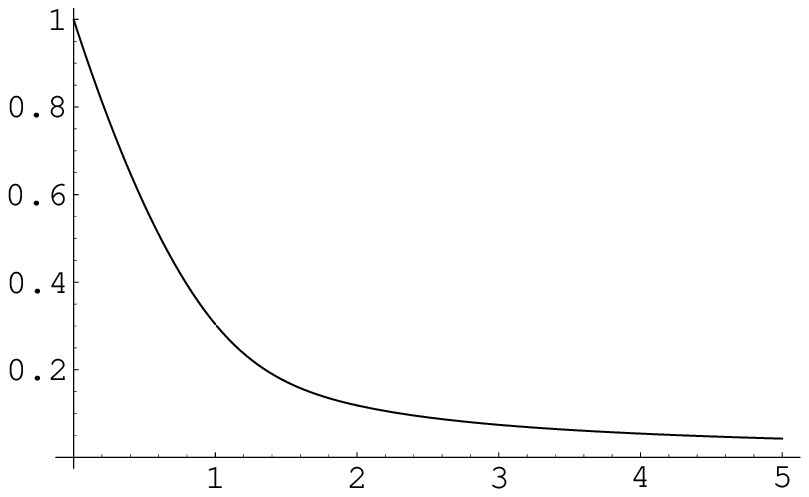}
\includegraphics*[scale=0.57, bb=90 0 280 150]{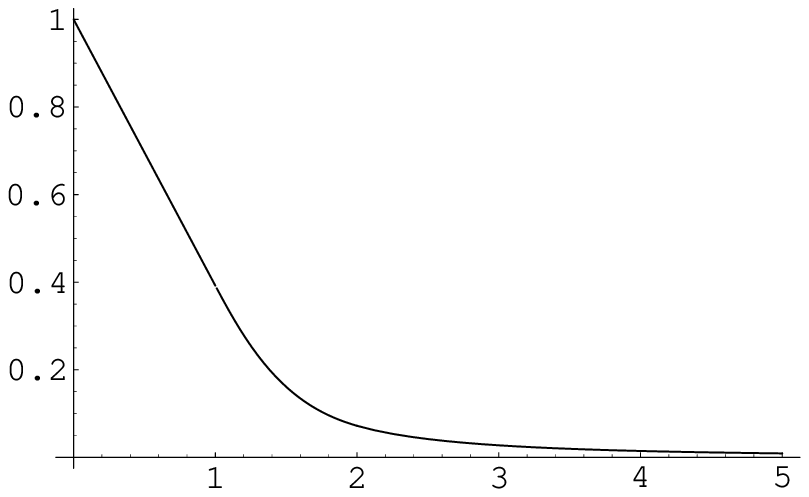}
\includegraphics*[scale=0.57, bb=90 0 280 150]{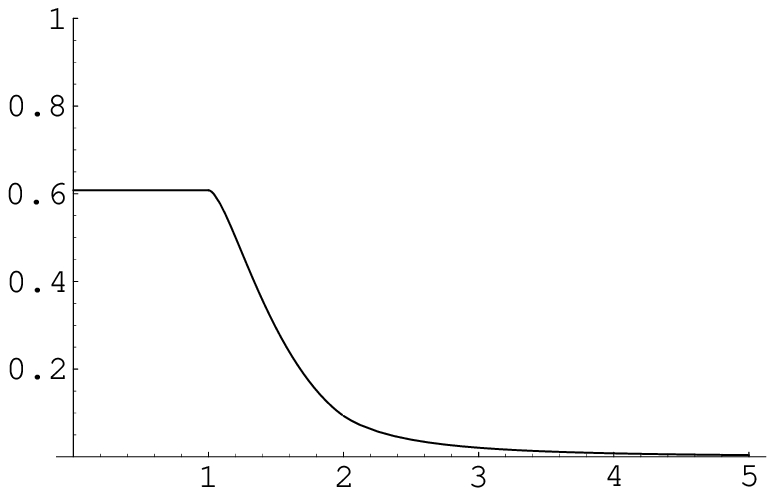}
\caption{The graphs of $\P(t)$, $\G(t)$, and respectively $g(t)$}
\end{figure}

We note the equalities
\begin{equation}\label{1.3}
\G(t)=-\P^\prime (t),\qquad t>0,
\end{equation}
and
\begin{equation}\label{1.4}
\begin{split}
& g(t):=-\G^\prime (t) =\P^{\prime \prime}(t) \\
& =\frac{6}{\pi^2}
\begin{cases}
1 & \mbox{\rm if $0<t\leq 1$;} \\
\mbox{\small $\displaystyle \frac{1}{t}+2\left(
1-\frac{1}{t}\right)^2 \ln \left( 1-\frac{1}{t} \right)
-\frac{1}{2}\left( 1-\frac{2}{t}\right)^2 \ln \left| 1-\frac{2}{t}
\right|$}  & \mbox{\rm if $t>1$.}
\end{cases}
\end{split}
\end{equation}
The latter also yields
\begin{equation}\label{1.5}
g(t)=\frac{24}{\pi^2 t^2} \sum\limits_{n=1}^\infty
\frac{2^n-1}{n(n+1)(n+2) t^n} ,\qquad t\geq 2.
\end{equation}
Remarkably, formulas \eqref{1.4} and \eqref{1.5} were found by
Dahlqvist \cite{Dah}. That approach however does not  provide a
rigorous proof for the existence of the limit distribution,
because it fails to control in a quantitative way the uniform
distribution of his variable $\Delta$ (see the comments after
formulas (75) and (86) in \cite{Dah}).
\par
In the final section we use some standard analysis arguments and
properties of the dilogarithm and trilogarithm to estimate
\begin{equation*}
C_\eps :=\ln \int_{\Sigma_\eps^+} \tau_\eps \, d\nu_\eps
-\int_{\Sigma_\eps^+} \ln \tau_\eps \, d\nu_\eps .
\end{equation*}
It was conjectured by Friedman, Kubo and Oono \cite{FOK} that
$C_\eps$ is convergent as $\eps \rightarrow 0^+$. Its hypothetical
limit $C$ was estimated to be $0.44\pm 0.001$ in \cite{FOK} and
$\approx 0.43$ in \cite{BoDo}. This conjecture was known to imply
\cite{Ch2,FOK} the asymptotic formula
\begin{equation*}
h(T_\eps)=-2\ln \eps+2-C+o(1) \qquad \mbox{\rm as $\eps
\rightarrow 0^+$}
\end{equation*}
for the KS entropy of the associated billiard map. In \cite{Ch1}
Chernov proved that $C_\eps$ remains bounded when $\eps
\rightarrow 0^+$, without giving however any estimate for the
bounds. The constant $C$ was identified by Dahlqvist \cite[formula
(73)]{Dah} as being $3\ln
2-\frac{9\zeta(3)}{4\zeta(2)}=0.43522513609\ldots$. The conjecture
of Friedman, Kubo and Oono, in the more precise form provided by
Dahlqvist, follows now from Theorem \ref{T1.2}.

\begin{theorem}\label{T1.3}
In the small scatterer limit $\eps \rightarrow 0^+$ the following
holds:
\begin{itemize}
\item[{\em (i)}] $\displaystyle C_\eps =3\ln
2-\frac{9\zeta(3)}{4\zeta(2)}+o(1).$ \item[{\em (ii)}]
$\displaystyle h(T_\eps)=-2\ln \eps+2 -3\ln
2+\frac{9\zeta(3)}{4\zeta(2)}+o(1).$
\end{itemize}
\end{theorem}
\par
These methods work for any convex scatterer due to the good error
control they give when integrating over the velocity in very short
intervals. To keep the presentation of the paper neat we have
chosen to only consider circular scatterers.
\par
In dimension $\geq 3$ the problem of the existence of the limiting
distribution of the free-path length in the small-scatterer limit
remains open and is manifestly difficult. Partial results in this
direction have appeared in \cite{BGW,Gol,GW}.

\section{Farey fractions and summation over primitive lattice
points} In this section we collect some basic properties of Farey
fractions and outline the summation method that will allow us to
estimate the limit distribution of the free path length when the
size of scatterers tends to zero.
\par
For each positive integer $Q$, let $\FQ$ denote the set of Farey
fractions of order $Q$. These are the rational numbers
$\gamma=\frac{a}{q}$ with coprime integers $a,q$ such that $1\leq
a\leq q\leq Q$. For each interval $I\subseteq [0,1]$ the number of
elements in the set
\begin{equation*}
\FI=I\cap \FQ
\end{equation*}
can be expressed, using elementary arguments on M\" obius and
Euler-Maclaurin summation, as
\begin{equation*}
\# \FI =\frac{Q^2\vert I\vert}{2\zeta(2)} +O(Q\ln Q).
\end{equation*}
\par
If $\gamma=\frac{a}{q}<\gamma^\prime =\frac{a^\prime}{q^\prime}$
are two consecutive elements in $\FQ$, then
\begin{equation}\label{2.1}
a^\prime q-aq^\prime =1 \qquad \mbox{\rm and} \qquad q+q^\prime
>Q.
\end{equation}
This shows on the one hand that the denominators of consecutive
Farey fractions of order $Q$ are exactly the primitive integer
points in the set
\begin{equation*}
Q\TT=\{(Qx,Qy)\, ;\ 0<x,y\leq 1,\ x+y>1\},
\end{equation*}
and on the other hand that denominators uniquely determine
consecutive Farey fractions. For instance, $a$ is the unique
integer in $[0,q]$ for which $(q-a)q^\prime =1\pmod{q}$.
\par
In many instances in this paper we will seek to estimate sums of
type
\begin{equation*}
S_{f,\Omega,I}(Q)=\sum\limits_{\substack{\gamma \in \FI \\
(q,q^\prime)\in Q\Omega}} f(q,q^\prime,a) ,
\end{equation*}
where $I\subseteq [0,1]$ is an interval, $\Omega \subseteq \TT$ a
region, and $f$ a $C^1$ function. These kinds of sums can be
roughly approximated by some integrals, with control on error
terms given by the following two results which will be
systematically used in this process. The first one is a standard
fact and is a plain consequence of the M\" obius summation (for a
proof see \cite[Lemma 2.3]{BCZ}).

\begin{lemma}\label{LEMMA1}
Let $0<a<b$ and $f$ be a $C^1$ function on $[a,b]$. Then
\begin{equation*}
\sum\limits_{a<k\leq b} \frac{\varphi(k)}{k}\, f(k)=
\frac{1}{\zeta(2)} \int_a^b f(x)\, dx+ O\left( \ln b \Big( \|
f\|_\infty +\int_a^b \vert f^\prime \vert \Big) \right) ,
\end{equation*}
where $\varphi$ denotes Euler's totient function.
\end{lemma}
\par
The second one is a consequence of Weil's type bounds for
Kloosterman sums (cf. \cite[Lemma 2.2]{BGZ}).

\begin{lemma}\label{LEMMA2}
Let $q\geq 1$ be an integer, $\II$ and $\JJ$ intervals with $\vert
\II \vert,\vert \JJ \vert <q$, $f$ a $C^1$ function on $\II \times
\JJ$, and $T\geq 1$ an integer. Then for all $\delta >0$
\begin{equation*}
\sum\limits_{\substack{a\in \II,\ b\in \JJ;\\ ab=1\hspace{-6pt}
\pmod{q}}} f(a,b)=\frac{\varphi(q)}{q^2} \iint\limits_{\II \times
\JJ} f(x,y)\, dx\, dy +{\mathcal E}, \end{equation*} with
\begin{equation*}
{\mathcal E}={\mathcal E}(q,T,f,\vert \II\vert,\vert
\JJ\vert,\delta) \ll_\delta T^2 q^{\frac{1}{2}+\delta} \|
f\|_\infty +Tq^{\frac{3}{2}+\delta} \| Df\|_\infty + \frac{\vert
\II\vert \, \vert \JJ \vert \, \| Df\|_\infty}{T},
\end{equation*}
where we denote $\| \cdot\|_\infty =\| \cdot\|_{\infty,\II \times
\JJ}$ and $Df=\vert \frac{\partial f}{\partial x}\vert+\vert
\frac{\partial f}{\partial y} \vert$.
\end{lemma}
\par
When $\Omega =\{ (x,y)\, ;\, \alpha <x\leq \beta, \ \xi(x)\leq
y\leq \eta(x)\}$ is a subset of $\TT =\{ (x,y) \, ;\, 0<x,y\leq 1,
x+y>1 \}$, the above mentioned properties of Farey fractions lead
to
\begin{equation*}
S_{f,\Omega,I}(Q)=\sum\limits_{\alpha Q< q\leq \beta Q}
\hspace{-6pt}
\sum\limits_{\substack{q^\prime \in QJ_\Omega (q/Q)\\
a\in qI \\ (q-a)q^\prime =1\hspace{-6pt} \pmod{q}}}\hspace{-10pt}
f(q,q^\prime ,a) = \sum\limits_{\alpha Q< q\leq \beta Q}
\sum\limits_{\substack{q^\prime \in QJ_\Omega (q/Q) \\
a\in q(1-I) \\ aq^\prime =1\hspace{-6pt} \pmod{q}}}\hspace{-10pt}
f(q,q^\prime ,q-a),
\end{equation*}
where we denote
\begin{equation*}
J_\Omega (x)=\{ y\, ;\, (x,y)\in \Omega\}=[\xi (x),\eta (x)]
\subseteq (1-x,1],\qquad x\in (\alpha,\beta].
\end{equation*}
\par
The inner sum above is mastered by Lemma \ref{LEMMA2}, being
approximated by
\begin{equation*}
\frac{\varphi(q)}{q^2} \iint_{QJ_\Omega (q/Q)\times q(1-I)}
\hspace{-18pt} f(q,q^\prime,q-a) \ dq^\prime \, da=
\frac{\varphi(q)}{q^2} \iint_{QJ_\Omega(q/Q)\times qI}
\hspace{-12pt} f(q,q^\prime,a)\ dq^\prime\, da.
\end{equation*}
Thus
\begin{equation*}
S_{f,\Omega,I}(Q)=\sum\limits_{\alpha Q< q\leq \beta Q}
\frac{\varphi (q)}{q}\, V(q)+\mbox{\rm error term},
\end{equation*}
where we take
\begin{equation*}
V(q)=\frac{1}{q} \iint_{QJ_\Omega (q/Q)\times qI}
f(q,q^\prime,a)\, dq^\prime \, da= Q\iint_{J_\Omega (q/Q)\times I}
 f(q,Qy,q\gamma )\, dy\, d\gamma .
\end{equation*}
When the error term is small enough, this sum is mastered by Lemma
\ref{LEMMA1}, giving
\begin{equation*}
\begin{split}
S_{f,\Omega,I}(Q) & =\frac{1}{\zeta(2)} \int_{\alpha Q}^{\beta Q}
V(q)\, dq +\mbox{\rm error}\\ & =\frac{Q}{\zeta(2)} \int_{\alpha
Q}^{\beta Q} dq \iint_{J_\Omega (q/Q)\times I}
 dy\, d\gamma \  f(q,Qy,q\gamma) +\mbox{\rm error} \\
& =\frac{Q^2}{\zeta(2)} \int_\alpha^\beta dx \iint_{J_\Omega
(x)\times I} dy\, d\gamma \ f(Qx,Qy,Qx\gamma )+\mbox{\rm error}
\\ & =\frac{1}{\zeta(2)} \iiint_{Q\Omega \times I}
f(v,w,v\gamma)\ dv\, dw\, d\gamma +\mbox{\rm error} .
\end{split}
\end{equation*}

\section{A partition of the unit square}
In this section we give an account on the three-strip partition
mentioned in the introduction. This approach, slightly different
from that in \cite{CG}, is suitable for computations involving
Farey fraction partitions of the unit interval.
\par
In the first part of this section we shall consider a fixed
(small) $\eps >0$ and let $Q=[\frac{1}{2\eps}]$ be the integer
part of $\frac{1}{2\eps}$. For each $\gamma=\frac{a}{q}\in \FQ$,
consider the points $N_\gamma(q,a+\eps)$ and $S_\gamma(q,a-\eps)$.
Consider also the points $N_0(0,\eps)$ and $S_0(0,-\eps)$, and
denote by $\SS_\gamma$ the strip determined by the lines
$N_0N_\gamma$ and $S_0S_\gamma$.
\par
A segment does not interfere with an open strip when their
intersection is empty. Throughout this section
$\gamma=\frac{a}{q}<\gamma^\prime=\frac{a^\prime}{q^\prime}$ will
be two consecutive fractions in $\FQ$, so that \eqref{2.1} is
fulfilled. In particular this gives
\begin{equation}\label{3.1}
2\eps \max \{ q,q^\prime \} \leq 2\eps Q\leq 1 \quad \mbox{\rm
and} \quad 2\eps (q+q^\prime)>2\eps (Q+1) >1.
\end{equation}
\par
The slope of a segment $AB$ is denoted by $\t_{AB}$. Set
$t_P=t_{OP}$.

\begin{figure}[ht]
\begin{center}
\unitlength 0.4mm
\begin{picture}(200,70)(15,0)
\path(0,0)(0,8)(165,65)(165,57)(0,0)(220,27)(220,35)(0,8)

\put(0,-5){\makebox(0,0){{\small $S_0(0,-\eps)$}}}
\put(-5,15){\makebox(0,0){{\small $N_0(0,\eps)$}}}

\put(185,70){\makebox(0,0){{\small $N_{\gamma'} (q',a'+\eps)$}}}
\put(185,52){\makebox(0,0){{\small $S_{\gamma'} (q',a'-\eps)$}}}

\put(235,40){\makebox(0,0){{\small $N_{\gamma} (q,a+\eps)$}}}
\put(235,22){\makebox(0,0){{\small $S_{\gamma} (q,a-\eps)$}}}
\end{picture}
\end{center}
\caption{The strips $\SS_\gamma$ and $\SS_{\gamma^\prime}$}
\label{Figure2}
\end{figure}

\begin{lemma}\label{LEMMA3}
The segment $N_\gamma S_\gamma$ does not interfere with the strip
$\SS_{\gamma^\prime}$, and the segment
$N_{\gamma^\prime}S_{\gamma^\prime}$ does not interfere with the
strip $\SS_\gamma$.
\end{lemma}
\par
\begin{proof} First, we show that $S_{\gamma^\prime}$ lies above the line
$N_0 N_\gamma$ of equation $ y-\eps-\frac{ax}{q}=0, $ which
amounts to $ a^\prime -2\eps -\frac{aq^\prime}{q}\geq 0. $ The
latter is equivalent to $1-2\eps q\geq 0$, which is true by
\eqref{3.1}.
\par
Furthermore, $N_\gamma$ lies below the line $S_0S_{\gamma^\prime}$
of equation $y+\eps-\frac{a^\prime x}{q^\prime}=0, $ as a result
of $ a+2\eps<\frac{a^\prime q}{q^\prime}$ being equivalent to
$2\eps q^\prime <1.$
\end{proof}
\par
For each $k\in \N_0 =\{ 0,1,2,\dots \}$ set
\begin{equation*}
\begin{split}
& q_k=q^\prime+kq,\qquad a_k=a^\prime+ka,\quad
q_k^\prime =q+kq^\prime ,\quad a_k^\prime =a+ka^\prime ,\\
& \gamma_k =\frac{a_k}{q_k} ,\quad t_k=\frac{a_k-2\eps}{q_k} ,
\quad u_k=\frac{a_k^\prime+2\eps}{q_k^\prime} ,\\
& \alpha_k =\arctan t_k ,\quad \beta_k=\arctan u_k .
\end{split}
\end{equation*}
\par
The following three relations hold for every $k\in \N = \{
1,2,\dots \}$
\begin{equation}\label{3.2}
a_{k-1}q_k-a_kq_{k-1}=1=a_k^\prime q_{k-1}^\prime - a_{k-1}^\prime
q_k^\prime ,
\end{equation}
\begin{equation}\label{3.3}
a_{k-1}q-aq_{k-1}=1 =a^\prime q_{k-1}^\prime -a^\prime_{k-1}
q^\prime ,
\end{equation}
\begin{equation}\label{3.4}
\min\limits_{k\geq 1} \{ 2\eps q_k^\prime ,2\eps q_k \} \geq 2\eps
(q^\prime+q)>1.
\end{equation}
As a result of \eqref{2.1} and \eqref{3.1}--\eqref{3.4}, it is
seen that
\begin{equation*}
\gamma^\prime =\gamma_0=\frac{a^\prime}{q^\prime}
>\gamma_1>\gamma_2 >\ldots >\gamma_k \
\stackrel{k\to \infty}{\longrightarrow} \  \frac{a}{q} =\gamma ,
\end{equation*}
and that
\begin{equation*}
\begin{split}
\gamma =\frac{a}{q} \ \stackrel{\infty \leftarrow
k}{\longleftarrow} & \ t_k \leq t_{k-1} \leq \ldots \leq t_1 \leq
t_0=\frac{a^\prime-2\eps}{q^\prime} \\ & <\frac{a+2\eps}{q}= u_0
\leq u_1\leq \ldots \leq u_{k-1} \leq u_k \ \stackrel{k\rightarrow
\infty}{\longrightarrow} \ \frac{a^\prime}{q^\prime}
=\gamma^\prime .
\end{split}
\end{equation*}
So putting
\begin{equation*}
I_{\gamma,0}=(t_0,u_0],\quad I_{\gamma,k}=(t_k,t_{k-1}], \quad
I_{\gamma,-k}=( u_{k-1},u_k],\qquad k\in \N ,
\end{equation*}
we end up with a partition $(I_{\gamma,k})_{k\in \Z}$ of the
interval $(\gamma,\gamma^\prime )$.
\par
Next we consider the points $N_{\gamma_k} (q_k,a_k+\eps)$ and
$S_{\gamma_k} (q_k,a_k-\eps)$, proving

\begin{lemma}\label{LEMMA4}
The following inequalities hold for every $k\geq 1${\em :}
\par
{\em (i)} $\quad \t_{N_0 N_{\gamma_{k-1}}}
>\t_{N_0 N_{\gamma_k}} > \t_{N_0 S_{\gamma_{k-1}}} \geq
\t_{N_0 S_{\gamma_k}} \geq \t_{N_0 N_\gamma} $.
\par
{\em (ii)} $\quad \t_{S_0 S_{\gamma^\prime}} \geq \t_{S_0
N_\gamma} >\t_{S_0 S_{\gamma_1}} \geq \t_{S_0 S_{\gamma_k}} $.
\end{lemma}
\par
\begin{proof} The inequalities in (i) are equivalent to
\begin{equation*}
\frac{a_{k-1}}{q_{k-1}} > \frac{a_k}{q_k}
> \frac{a_{k-1}-2\eps}{q_{k-1}} \geq
\frac{a_k-2\eps}{q_k} \geq \frac{a}{q}  ,
\end{equation*}
which follow from \eqref{3.2}, \eqref{3.3}, \eqref{3.1}, and from
\eqref{3.4}.
\par
The inequalities in (ii) are equivalent to
\begin{equation*}
\frac{a^\prime}{q^\prime} \geq \frac{a+2\eps}{q}>\frac{a_1}{q_1}=
\frac{a^\prime+a}{q^\prime +q}\geq \frac{a_k}{q_k}  ,
\end{equation*}
which follow from \eqref{2.1}, \eqref{3.2}, \eqref{3.4}, and $a_1
q_k-a_k q_1 =k-1$.
\end{proof}
\par
Consider the half-infinite strip
\begin{equation*}
\SS=\SS_\omega=\{ ( x,y+x\tan \omega ) \, ;\, x>0,\ -\eps \leq
y\leq \eps \}
\end{equation*}
of direction $\omega$, top line passing through $N_0$, and bottom
line passing through $S_0$. Assume that $\gamma <\tan \omega
<\gamma^\prime$. For each $y_0\in [-\eps,\eps]$ we wish to find
the first vertical segment of form $\{ m\} \times
[n-\eps,n+\eps]$, $m,n\in \N$, that intersects the line of slope
$\tan \omega$ passing through $(0,y_0)$. In other words, we wish
to calculate
\begin{equation*}
q(\omega,y_0)=\inf \{ n\in \N \, ;\, \| y_0+n\tan \omega\|\leq
\eps \},
\end{equation*}
where we denote $\| x\|=\dist (x,\Z)$, $x\in \R$. We shall assume
that $\tan \omega$ is irrational and split the discussion
according to the three cases where the slope of $\omega$ belongs
to one of the intervals $\big(\frac{a}{q},\frac{a^\prime
-2\eps}{q^\prime}\big]$, $\big( \frac{a^\prime -2\eps}{q^\prime},
\frac{a+2\eps}{q}\big]$ or $\big(
\frac{a+2\eps}{q},\frac{a^\prime}{q^\prime}\big]$.

\begin{proposition}\label{PROP1}
Let $\gamma=\frac{a}{q}<\gamma^\prime=\frac{a^\prime}{q^\prime}$
be consecutive fractions in $\FQ$. Suppose $\tan \omega \in
(\frac{a}{q},\frac{a^\prime-2\eps}{q^\prime}]$ is irrational and
$\tan \omega \in I_{\gamma,k}=(t_k,t_{k-1}]$ for some $k\in \N$.
Set
\begin{equation*}
\begin{split}
& w_{B_k}=w_{B_k}(\omega)=q_k \tan \omega -a_k+2\eps ,\\
& w_{C_k}=w_{C_k}(\omega)=-q_{k-1}\tan \omega +a_{k-1}-2\eps ,\\
& w_{A_+} =w_{A_+} (\omega)=-q\tan \omega +a+2\eps  ,\\
& I_{A_+}: = [-\eps,-\eps+w_{A_+}),\\ & I_{B_k}
:=(-\eps+w_{A_+}+w_{C_k},\eps]=(\eps -w_{B_k},\eps ],\\
& I_{C_k} :=[-\eps+w_{A_+},-\eps+w_{A_+}+w_{C_k}],
\end{split}
\end{equation*}
and
\begin{equation*} L(\omega,y_0)=\begin{cases}
L_{A_+}(\omega):=q & \mbox{if $y_0\in I_{A_+};$} \\
L_{C_k} (\omega):=q_{k+1} & \mbox{if $y_0\in I_{C_k}
;$} \\
L_{B_k} (\omega):=q_k & \mbox{if $y_0\in I_{B_k} .$}
\end{cases}
\end{equation*}
\par
Then for any $\star \in \{A_+,B_k,C_k\}$ we have $0\leq w_\star
\leq 2\eps$ and
\begin{equation*}
q(\omega,y_0)=L(\omega,y_0) =L_\star (\omega), \qquad y_0 \in
I_\star .
\end{equation*}
\par
Furthermore, if $S_\star$ denotes the parallelogram of height $\{
0\}\times I_\star$, angle $\omega$ between its side and the
horizontal direction, and side length $\frac{L_\ast(\omega)}{\cos
\omega}$, then
\begin{equation*}
\ar(S_{A_+})+\ar(S_{B_k})+\ar(S_{C_k})=w_{A_+}L_{A_+}+
w_{B_k}L_{B_k}+w_{C_k}L_{C_k}=1.
\end{equation*}
Moreover, $\{ S_{A_+},S_{B_k},S_{C_k}\}$ mod $\Z^2$ provides a
partition of the unit square $[0,1)^2$ {\em (}we allow the
boundaries of these three sets to intersect{\em )}.
\end{proposition}

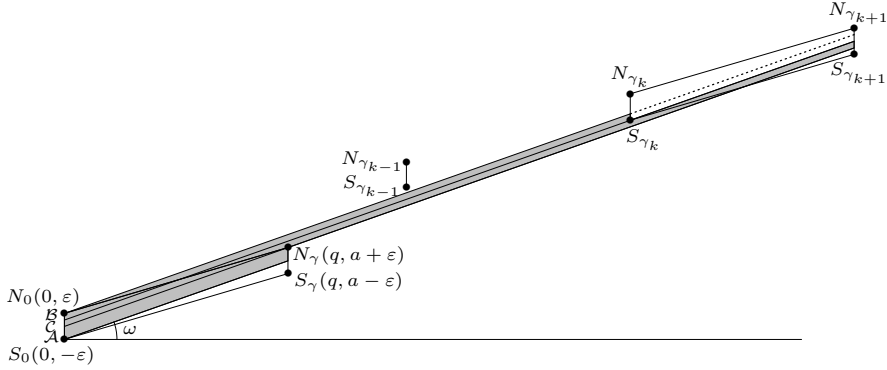
\begin{figure}[ht]
\begin{center}
\unitlength 0.35mm
\begin{picture}(200,115)(50,0)
\texture{c 0}
\shade\path(0,0)(85,30)(85,35)(300,110.882)(300,113.382)(215,83.382)(215,85.8823)(0,10)(0,0)

\path(0,0)(280,0) \path(0,0)(85,25)(85,35)(0,10)(0,0)
\path(0,0)(85,30) \path(85,35)(0,10)
\path(215,83.3823)(215,93.3823)(300,118.382)(300,108.382)(215,83.3823)
\dottedline{2}(215,85.882)(300,115.8823)
\path(130,57.8823)(130,67.8823) \path(0,7.5)(300,113.3823)
\path(0,5)(300,110.88236)

\put(-5,-6){\makebox(0,0){{\scriptsize $S_0(0,-\eps)$}}}
\put(-8,16){\makebox(0,0){{\scriptsize $N_0(0,\eps)$}}}
\put(-5,1){\makebox(0,0){{\tiny ${\mathcal A}$}}}
\put(-5,5){\makebox(0,0){{\tiny ${\mathcal C}$}}}
\put(-5,9.5){\makebox(0,0){{\tiny ${\mathcal B}$}}}
\put(0,0){\makebox(0,0){{\tiny $\bullet$}}}
\put(0,10){\makebox(0,0){{\tiny $\bullet$}}}
\put(85,25){\makebox(0,0){{\tiny $\bullet$}}}
\put(85,35){\makebox(0,0){{\tiny $\bullet$}}}

\put(300,108.382){\makebox(0,0){{\tiny $\bullet$}}}
\put(300,118.382){\makebox(0,0){{\tiny $\bullet$}}}
\put(215,83.3823){\makebox(0,0){{\tiny $\bullet$}}}
\put(215,93.3823){\makebox(0,0){{\tiny $\bullet$}}}

\put(130,58.05){\makebox(0,0){{\tiny $\bullet$}}}
\put(130,67.3823){\makebox(0,0){{\tiny $\bullet$}}}

\put(117,58.05){\makebox(0,0){{\scriptsize $S_{\gamma_{k-1}}$}}}
\put(117,67.3823){\makebox(0,0){{\scriptsize $N_{\gamma_{k-1}}$}}}

\put(108,32){\makebox(0,0){{\scriptsize $N_\gamma (q,a+\eps)$}}}
\put(108,22){\makebox(0,0){{\scriptsize $S_\gamma (q,a-\eps)$}}}

\put(302,101){\makebox(0,0){{\scriptsize $S_{\gamma_{k+1}}$}}}
\put(302,124){\makebox(0,0){{\scriptsize $N_{\gamma_{k+1}}$}}}

\put(215,100){\makebox(0,0){{\scriptsize $N_{\gamma_k}$}}}
\put(220,76){\makebox(0,0){{\scriptsize $S_{\gamma_k}$}}}
\put(25,3.5){\makebox(0,0){{\scriptsize $\omega$}}}
\put(0,0){\arc{40}{-0.34}{0}}
\end{picture}
\end{center}
\caption{The case $\tan \omega \in I_{\gamma,k}$, $k\in \N$}
\label{Figure3}
\end{figure}

\begin{proof} Taking stock on Lemma \ref{LEMMA4}, we notice that the line of
slope $\tan \omega$ through $S_0$ intersects the vertical line
$N_\gamma S_\gamma$ at a point between $N_\gamma$ and $S_\gamma$
(see Figure \ref{Figure3}). Also, because $\t_{N_0
S_{\gamma_k}}=t_k <\tan \omega \leq t_{k-1}= \t_{N_0
S_{\gamma_{k-1}}} <t_{N_0 N_{\gamma_k}}$, the line of slope $\tan
\omega$ through $N_0$ (respectively through $N_\gamma$) intersects
the line $N_{\gamma_k} S_{\gamma_k}$ (respectively
$N_{\gamma_{k+1}} S_{\gamma_{k+1}}$) between $N_{\gamma_k}$ and
$S_{\gamma_k}$ (respectively between $N_{\gamma_{k+1}}$ and
$S_{\gamma_{k+1}}$). The segment
$N_{\gamma_{k-1}}S_{\gamma_{k-1}}$ is placed above these two
parallel lines because $\tan \omega \leq t_{k-1} =t_{N_0
S_{\gamma_{k-1}}}$.
\par
Next, we find that the intersections with the vertical axis of the
lines $y-a-\eps=(x-q)\tan \omega$ and $y-a_k+\eps=(x-q_k)\tan
\omega$ which have slope $\tan \omega$ and pass through $N_\gamma$
and respectively $S_{\gamma_k}$, are $(0,\eps+a-q\tan \omega)$ and
respectively $(0,-\eps+a_k-q_k \tan \omega)$, whence the required
values of $w_{A_+}$, $w_{B_k}$ and $w_{C_k}$ follow. Notice that
\begin{equation*}
\begin{split}
& 2\eps > w_{A_+}=2\eps+a-q\tan \omega\geq
2\eps+a-q\frac{a^\prime-2\eps}{q^\prime}=
\frac{2\eps(q+q^\prime)-1}{q^\prime}>0 ,\\
& 2\eps >\frac{1-2\eps q}{q_{k-1}}=q_k
\frac{a_{k-1}-2\eps}{q_{k-1}} -a_k+2\eps \geq w_{B_k}
=q_k \tan \omega -a_k+2\eps >0 , \\
& 2\eps >\frac{1-2\eps q}{q_k}=a_{k-1}-2\eps -q_{k-1}
\frac{a_k-2\eps}{q_k} >a_{k-1} -2\eps -q_{k-1}\tan \omega =w_{C_k}
\geq 0.
\end{split}
\end{equation*}
Besides one clearly has
\begin{equation*}
w_{A_+}+w_{B_k}+w_{C_k}=2\eps ,
\end{equation*}
and it is easy to check by a direct calculation that
\begin{equation*}
\sum\limits_{\star \in \{ A_+,B_k,C_k \}} \hspace{-10pt} \ar
(S_\star) =\sum\limits_{\star \in \{ A_+,B_k,C_k \}}
\hspace{-10pt} w_\star L_\star =1.
\end{equation*}
It remains to check that the interiors of the subsets $S_\star
\hspace{-4pt} \mod \Z^2 \subseteq [0,1)^2$, $\star \in \{
A_+,B_k,C_k\}$, are disjoint. If not, there exist two points
$P,P^\prime$ inside $\cup_\star S_\star$ such that $P-P^\prime \in
\Z^2$. The latter is preserved by translating the segment
$PP^\prime$ to a parallel segment. Owing to the shape of
$\cup_\star S_\star$ we may thus assume that, say, $P$ lies on the
$y$-axis; hence $P=P(0,y_0)$ and $P^\prime =P^\prime (n,m+y_0)$
for some $y_0\in [-\eps,\eps]$, $m,n\in \N^*$. The line of slope
$\tan \omega$ which passes through $P^\prime$ intersects the
$y$-axis at $(0,m+y_0-n\tan \omega)$. Hence $-\eps \leq
m+y_0-n\tan \omega \leq \eps$, which shows that $\| y_0-n\tan
\omega \|=\| -y_0+n\tan \omega \| \leq \eps$. By the first part of
the proposition this gives $n\geq L(\omega,-y_0)$, thus $P^\prime$
must belong to the boundary, which is a contradiction.
\end{proof}

\begin{proposition}\label{PROP2}
Let $\gamma=\frac{a}{q}<\gamma^\prime=\frac{a^\prime}{q^\prime}$
be consecutive fractions in $\FQ$. Suppose $\tan \omega \in (
\frac{a^\prime -2\eps}{q^\prime},\frac{a^\prime}{q^\prime})$ is
irrational.

{\em (i)} If $\tan \omega \in I_{\gamma,0} =(t_0,u_0]$, then the
analog of Proposition {\rm \ref{PROP1}} holds true,
with\footnote{Note that in both cases $q<q^\prime$ or $q^\prime
<q$ we get $0\leq w_{C_0}<2\eps$.}
\begin{equation*}
\begin{split}
& w_{B_0}=w_{B_0} (\omega)=q^\prime \tan \omega -a^\prime+2\eps
\in (0,2\eps) ,\\
& w_{C_0}=w_{C_0} (\omega)=-(q^\prime-q)\tan \omega +a^\prime
-a-2\eps \in [0,2\eps) ,\\
& w_{A_0}=w_{A_0} (\omega)=w_{A_+}(\omega)= -q\tan \omega +a+2\eps
\in [0,2\eps) ,\\
& I_{A_0} :=[-\eps,\eps+w_{A_0}),\\ &
I_{B_0}=(-\eps+w_{A_0}+w_{C_0},\eps]=(\eps -w_{B_0},\eps],\\
& I_{C_0}:= [-\eps+w_{A_0},-\eps+w_{A_0}+w_{C_0}],\\
& L(\omega,y_0)=\begin{cases} L_{A_0} (\omega):=q & \mbox{if
$y_0 \in I_{A_0};$} \\
 L_{C_0}(\omega):=q^\prime +q & \mbox{if $y_0 \in I_{C_0};$} \\
L_{B_0}(\omega):=q^\prime & \mbox{if $y_0 \in I_{B_0}$.}
\end{cases}
\end{split}
\end{equation*}

{\em (ii)} If $k\in \N$ and $\tan \omega \in
I_{\gamma,-k}=(u_{k-1},u_k]$, then the analog of Proposition {\rm
\ref{PROP1}} holds true, with
\begin{equation*}
\begin{split}
& w_{B_-}=w_{B_-} (\omega)=w_{B_0}(\omega)=q^\prime\tan \omega
-a^\prime+2\eps \in (0,2\eps), \\
& w_{C_{-k}}= w_{C_{-k}} (\omega)=q^\prime_{k-1}\tan
\omega-a^\prime_{k-1}-2\eps\in (0,2\eps), \\
& w_{A_{-k}}=w_{A_{-k}} (\omega)=-q^\prime_k\tan \omega
+a^\prime_k +2\eps \in [0,2\eps) ,\\
& I_{A_{-k}} = [-\eps,-\eps+w_{A_{-k}}), \\ & I_{C_{-k}}
=[-\eps+w_{A_{-k}},-\eps+w_{A_{-k}}+w_{C_{-k}}],\\
& I_{B_-}
=(-\eps+w_{A_{-k}}+w_{C_{-k}},\eps]=(\eps-w_{B_-},\eps ],\\
 & L(\omega,y_0)=
\begin{cases} L_{A_{-k}}
(\omega):=q_k^\prime & \mbox{\sl if $ y_0 \in I_{A_{-k}};$} \\
L_{C_{-k}} (\omega):=q^\prime_{k+1} & \mbox{\sl if $ y_0 \in
I_{C_{-k}};$}
 \\ L_{B_-} (\omega):=q^\prime & \mbox{\sl if $ y_0 \in I_{B_-}
.$}
\end{cases}
\end{split}
\end{equation*}
\end{proposition}
\par
\begin{proof} (i) follows as in the proof of Proposition \ref{PROP1}, using
\begin{equation*}
\eps>a^\prime-\eps-q^\prime \tan \omega \geq a+\eps-q\tan \omega
\geq -\eps,\quad \tan \omega \in I_{\gamma,0}.
\end{equation*}

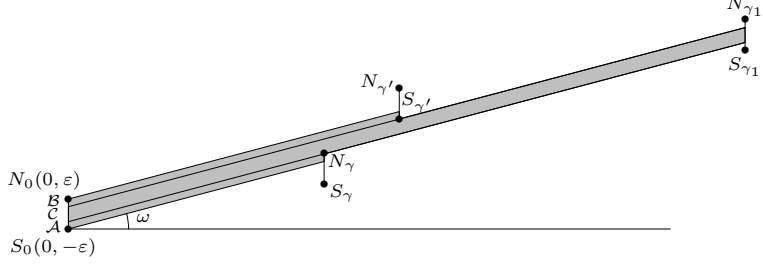
\begin{figure}[ht]
\begin{center}
\unitlength 0.4mm
\begin{picture}(200,60)(10,0)
\texture{c 0}
\shade\path(0,0)(85,22.5)(85,25)(225,62.05882)(225,67.058823)(110,36.61764)(110,39.11764)(0,10)(0,0)

\path(0,10)(0,0)(200,0) \path(0,2.5)(225,62.058823)
\path(0,7.5)(225,67.058823) \path(85,15)(85,25)
\path(225,59.5588)(225,69.5588) \path(110,36.6176)(110,46.6176)

\put(-5,-6){\makebox(0,0){{\scriptsize $S_0(0,-\eps)$}}}
\put(-8,16){\makebox(0,0){{\scriptsize $N_0(0,\eps)$}}}
\put(-5,1){\makebox(0,0){{\tiny ${\mathcal A}$}}}
\put(-5,5){\makebox(0,0){{\tiny ${\mathcal C}$}}}
\put(-5,9.5){\makebox(0,0){{\tiny ${\mathcal B}$}}}

\put(0,0){\makebox(0,0){{\tiny $\bullet$}}}
\put(0,10){\makebox(0,0){{\tiny $\bullet$}}}
\put(85,25){\makebox(0,0){{\tiny $\bullet$}}}
\put(85,15){\makebox(0,0){{\tiny $\bullet$}}}

\put(225,59.5588){\makebox(0,0){{\tiny $\bullet$}}}
\put(225,69.5588){\makebox(0,0){{\tiny $\bullet$}}}
\put(110,36.6176){\makebox(0,0){{\tiny $\bullet$}}}
\put(110,46.6176){\makebox(0,0){{\tiny $\bullet$}}}

\put(91,22){\makebox(0,0){{\scriptsize $N_\gamma$}}}
\put(91,12){\makebox(0,0){{\scriptsize $S_\gamma$}}}

\put(103,48){\makebox(0,0){{\scriptsize $N_{\gamma'}$}}}
\put(116,42){\makebox(0,0){{\scriptsize $S_{\gamma'}$}}}

\put(225,74){\makebox(0,0){{\scriptsize $N_{\gamma_1}$}}}
\put(225,54){\makebox(0,0){{\scriptsize $S_{\gamma_1}$}}}

\put(0,0){\arc{40}{-0.26}{0}}
\put(25,3.5){\makebox(0,0){{\scriptsize $\omega$}}}
\end{picture}
\end{center}
\caption{The case $\tan \omega \in I_{\gamma,0}$} \label{Figure4}
\end{figure}

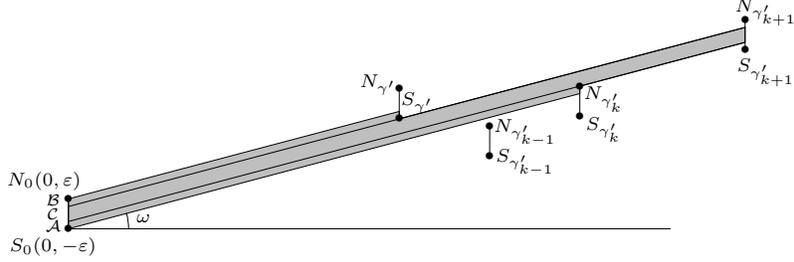
\begin{figure}[ht]
\begin{center}
\unitlength 0.4mm
\begin{picture}(200,60)(10,0)
\texture{c 0}
\shade\path(0,0)(170,45)(170,47.5)(225,62.05882)(225,67.058823)(110,36.61764)(110,39.11764)(0,10)(0,0)

\path(0,10)(0,0)(200,0) \path(0,2.5)(225,62.058823)
\path(0,7.5)(225,67.058823) \path(170,47.5)(170,37.5)
\path(225,59.5588)(225,69.5588) \path(110,36.6176)(110,46.6176)
\path(140,34.058823)(140,24.058823)

\put(-5,-6){\makebox(0,0){{\scriptsize $S_0(0,-\eps)$}}}
\put(-8,16){\makebox(0,0){{\scriptsize $N_0(0,\eps)$}}}
\put(-5,1){\makebox(0,0){{\tiny ${\mathcal A}$}}}
\put(-5,5){\makebox(0,0){{\tiny ${\mathcal C}$}}}
\put(-5,9.5){\makebox(0,0){{\tiny ${\mathcal B}$}}}

\put(0,0){\makebox(0,0){{\tiny $\bullet$}}}
\put(0,10){\makebox(0,0){{\tiny $\bullet$}}}
\put(170,47.5){\makebox(0,0){{\tiny $\bullet$}}}
\put(170,37.5){\makebox(0,0){{\tiny $\bullet$}}}

\put(225,59.5588){\makebox(0,0){{\tiny $\bullet$}}}
\put(225,69.5588){\makebox(0,0){{\tiny $\bullet$}}}
\put(110,36.6176){\makebox(0,0){{\tiny $\bullet$}}}
\put(110,46.6176){\makebox(0,0){{\tiny $\bullet$}}}
\put(140,34.058823){\makebox(0,0){{\tiny $\bullet$}}}
\put(140,24.058823){\makebox(0,0){{\tiny $\bullet$}}}

\put(152,32){\makebox(0,0){{\scriptsize
$N_{\gamma^\prime_{k-1}}$}}}
\put(152,22){\makebox(0,0){{\scriptsize
$S_{\gamma^\prime_{k-1}}$}}}
\put(178,43){\makebox(0,0){{\scriptsize $N_{\gamma^\prime_k}$}}}
\put(178,33){\makebox(0,0){{\scriptsize $S_{\gamma^\prime_k}$}}}

\put(103,48){\makebox(0,0){{\scriptsize $N_{\gamma'}$}}}
\put(116,42){\makebox(0,0){{\scriptsize $S_{\gamma'}$}}}

\put(232,72){\makebox(0,0){{\scriptsize
$N_{\gamma_{k+1}^\prime}$}}}

\put(232,52){\makebox(0,0){{\scriptsize
$S_{\gamma_{k+1}^\prime}$}}}

\put(0,0){\arc{40}{-0.26}{0}}
\put(25,3.5){\makebox(0,0){{\scriptsize $\omega$}}}
\end{picture}
\end{center}
\caption{The case $\tan \omega \in I_{\gamma,-k}$, $k\in \N$}
\label{Figure5}
\end{figure}

(ii) follows as in the proof of Proposition \ref{PROP1} using
\begin{equation*}
\eps>a^\prime -\eps-q^\prime \tan \omega > a_k^\prime +\eps
-q_k^\prime \tan \omega \geq -\eps,\qquad \tan \omega \in
I_{\gamma,-k}.
\end{equation*}
\end{proof}
\par
We now start investigating the case where the scatterers are
vertical slits. Propositions \ref{PROP1} and \ref{PROP2} will only
be applied for $\eps =\frac{1}{2Q}$, corresponding to the case of
vertical slits of height $\frac{1}{Q}$. The Lebesgue measure of a
Borel set $A$ in $\R^d$, $d=1,2,3$, will be denoted by $\vert
A\vert$.
\par
Throughout the paper $\widetilde{\tau}_\delta (x,\omega)$ will
denote the free path length in the periodic two-dimensional
Lorentz gas with vertical slits of height $2\delta$ as scatterers
centered at all integer lattice points. Given $\lambda >0$,
$I=[\tan \omega_0,\tan \omega_1]\subseteq [0,1]$ with $0\leq
\omega_0\leq \omega_1 \leq \frac{\pi}{4}$, and $Q\geq 1$ integer,
we denote
\begin{equation*}
\widetilde{\P}_{I,Q}(\lambda)=\big| \{ (x,\omega) \, ; \, x\in
[0,1)^2,\ \omega_0\leq \omega \leq \omega_1,\
\widetilde{\tau}_{1/(2Q)} (x,\omega)>\lambda \} \big| .
\end{equation*}
\par
Although the cases $0<t<1$, $1<t<2$, $t>2$, will be considered
separately, applying Propositions \ref{PROP1} and \ref{PROP2} to
$2\eps =\frac{1}{Q}$, can write for all $t,\varepsilon_\ast >0$
\begin{equation}\label{3.5}
\begin{split}
\widetilde{\P}_{I,Q} \bigg( \frac{t}{2\eps_\ast} \bigg) &
=\sum\limits_{\gamma \in \FI} \sum\limits_{k=1}^\infty
\int_{\alpha_k}^{\alpha_{k-1}} w_{C_k}(\omega)\max \left\{
q_{k+1}-\frac{t\cos \omega}{2\eps_\ast},0\right\} \, d\omega \\ &
\qquad +\mbox{\rm eight other similar terms where $\displaystyle
\frac{t\cos \omega}{2\eps_\ast}$ appears.}
\end{split}
\end{equation}

\begin{lemma}\label{LEMMA5}
For any interval $I=[\tan \omega_0,\tan \omega_1] \subseteq [0,1]$
such that $\vert I\vert \asymp \eps^c$ with fixed $0<c<1$ and
small $\eps >0$, and any {\em (}large{\em )} integer $Q=\frac{\cos
\omega_0}{2\eps}+O(\eps^{c-1})$, the estimate
\begin{equation}\label{3.6}
\widetilde{\P}_{I,Q} \bigg( \frac{t}{2\eps}\bigg)
=P_{I,Q}(t)+O(\eps^{2c}),
\end{equation}
holds uniformly in $t$ on compact subsets of $(0,\infty)$. Here
$P_{I,Q}(t)$ is obtained by substituting $tQ$ in place of
$\frac{t\cos \omega}{2\eps_\ast}$ in \eqref{3.5}, that is

\begin{equation}\label{3.7}
\begin{split}
P_{I,Q}(t) & :=\sum\limits_{\gamma \in \FI} \
\sum\limits_{k=1}^\infty \int_{\alpha_k}^{\alpha_{k-1}} \Big(
w_{C_k}(\omega) \max \{ q_{k+1}-tQ,0\}  \\ & \qquad
+w_{B_k}(\omega) \max
\{ q_k-tQ,0\} +w_{A_+} (\omega)\max \{ q-tQ,0\} \Big) \, d\omega \\
& +\sum\limits_{\gamma \in \FI} \ \int_{\alpha_0}^{\beta_0} \Big(
w_{C_0}(\omega)\max \{ q+q^\prime -tQ,0\} \\ & \qquad +w_{B_0}
(\omega)\max \{ q^\prime -tQ,0\} +w_{A_0}(\omega)\max \{ q-tQ,0\}
\Big) \, d\omega \\ & +\sum\limits_{\gamma \in \FI}
\sum\limits_{k=1}^\infty \ \int_{\beta_{k-1}}^{\beta_k} \Big(
w_{C_{-k}}(\omega)\max \{ q^\prime_{k+1}-tQ,0\} \\ & \qquad
+w_{A_{-k}} (\omega) \max \{ q^\prime_k-tQ,0\} +w_{B_-}(\omega)
\max \{ q^\prime -tQ,0\}\Big)\, d\omega ,
\end{split}
\end{equation}
with
\par
\begin{equation}\label{3.8}
\begin{split}
& w_{A_+}(\omega)=w_{A_0}(\omega)=Q^{-1}+a-q\tan \omega ,\\ &
w_{B-}(\omega)=w_{B_0}(\omega)=
q^\prime \tan \omega -a^\prime+Q^{-1} ,\\
& w_{C_0}(\omega)=Q^{-1} -w_{A_0}(\omega)-w_{B_0}(\omega) ,\\
& w_{B_k} (\omega)=q_k \tan \omega -a_k+Q^{-1},\quad w_{C_k}
(\omega)=a_{k-1}-Q^{-1}-q_{k-1} \tan \omega ,\\ & w_{C_{-k}}
(\omega)=q_{k-1}^\prime\tan \omega -a_{k-1}^\prime -Q^{-1}, \quad
w_{A_{-k}}(\omega) =Q^{-1}+a_k^\prime -q_k^\prime \tan
\omega ,\\
& \alpha_k=\arctan \frac{a_k-1/Q}{q_k},\quad \beta_k=\arctan
\frac{a_k^\prime+1/Q}{q_k^\prime} ,\quad k\in \N .
\end{split}
\end{equation}
\end{lemma}
\par
\begin{proof} Using the inequality $\max \{ w_A,w_B,w_C\} \leq
Q^{-1} \ll \eps$, which is a consequence of
$w_{A_+}+w_{B_k}+w_{C_k}=Q^{-1}$ and of the similar relations for
$k=0$ and $k\leq -1$, the estimate (see also \eqref{7.2})
\begin{equation*}
\sup\limits_{\omega \in I} \left| Q-\frac{\cos \omega}{2\eps}
\right| \ll \eps^{c-1}+\frac{\vert \cos \omega_0-\cos
\omega_1\vert}{2\eps} \leq \eps^{c-1}+\frac{\vert \tan
\omega_1-\tan \omega_0\vert}{2\eps} \ll \eps^{c-1},
\end{equation*}
and the inequalities
\begin{equation}\label{3.9}
\vert \max (x,0)-\max (y,0)\vert \leq \vert x-y\vert
\end{equation}
and
\begin{equation*}
\sum\limits_{\gamma \in \FI} \frac{1}{qq^\prime} \ll \vert I\vert
+\frac{1}{Q} \ll \eps^c,
\end{equation*}
it follows that we can replace $\frac{t\cos \omega}{2\eps}$ by
$tQ$ in \eqref{3.5} at a cost which is
\begin{equation*}
\ll \sum\limits_{\gamma \in \FI} \frac{1}{Q} \eps^{c-1}
\frac{1}{qq^\prime} \ll \eps^c \sum\limits_{\gamma \in \FI}
\frac{1}{qq^\prime} \ll \eps^{2c}.
\end{equation*}
\end{proof}
\par
Equality \eqref{3.7} will be at the center of most of the
forthcoming computations because it shows how the estimation of
distribution of the free path length reduces to estimates on sums
involving Farey fractions.
\par
There is an alternative approach to estimating
$\widetilde{\P}_{I,Q}$, by using a monotonicity argument instead
of the continuity argument which based on \eqref{3.9}. Such an
argument will be used in the proof of Theorem \ref{T1.2}.
\par
In the remainder of the paper given $I=[\tan \omega_0,\tan
\omega_1]\subseteq [0,1]$ we denote
\begin{equation}\label{3.10}
c_I=\int_I \frac{du}{1+u^2} =\omega_1-\omega_0 .
\end{equation}

\section{The case $0<t\leq 1$}
The aim of this section is to prove the following result

\begin{proposition}\label{PROP3}
Suppose $I$ is a subinterval of $[0,1]$ of size $\vert I\vert
\asymp Q^{-c}$ for some $0<c<1$. Then for every $c_1>0$ with
$c+c_1<1$ and $\delta>0$
\begin{equation*}
P_{I,Q}(t)=\left( 1-t+\frac{t^2}{2\zeta(2)} \right) c_I +O_\delta
( E_{c,c_1,\delta}(Q)) \qquad (Q\rightarrow \infty),
\end{equation*}
with
\begin{equation*}
E_{c,c_1,\delta}(Q)=Q^{\max \{ 2c_1-1/2+\delta,-c-c_1\}} .
\end{equation*}
The estimate is uniform in $t\in (0,1]$.
\end{proposition}
\par
Before starting to estimate $P_{I,Q}$, the following remark is in
order.

\newtheorem*{Note}{Remark 1}
\begin{Note}
{\em If $I\subseteq [0,1]$ is an interval with $\vert I\vert \geq
\frac{1}{Q}$, then as a consequence of $\gamma^\prime-\gamma
=\frac{1}{qq^\prime}\leq \frac{1}{Q}\leq \vert I\vert$ we have
\begin{equation*}
\sum\limits_{\substack{\gamma \in \FQ \\ \gamma \in I}} f(\gamma)=
\sum\limits_{\substack{\gamma \in \FQ \\ \gamma^\prime \in I}}
f(\gamma^\prime)+O( Q^{-1}\| f\|_\infty) .
\end{equation*}
As a result, replacing the condition $\gamma \in I$ by
$\gamma^\prime \in I$ only produces an error of order $Q^{-1}$,
which has no impact in any of the forthcoming estimates. Thus in
Propositions \ref{PROP3}, \ref{PROP4}, \ref{PROP5}, \ref{P8.2},
\ref{P9.1}, \ref{P10.1} the assumption $\vert I\vert \asymp
Q^{-c}$ can be replaced by the weaker assumption $\vert I\vert \ll
Q^{-c}$ and $\vert I\vert \geq Q^{-1}$.}
\end{Note}
\par
Then we notice that since $\min \{q_k,q^\prime_k\} \geq q+q^\prime
>tQ$ for all $k\geq 1$, we can write, according to \eqref{3.7} and
\eqref{3.8},
\begin{equation*}
\begin{split}
P_{I,Q}(t) = & \sum\limits_{\gamma \in \FI}
\sum\limits_{k=1}^\infty \int_{\alpha_k}^{\alpha_{k-1}}
S^{-}_{Q,\gamma,k}(\omega)\,d\omega + \sum\limits_{\gamma \in \FI}
\ \int_{\alpha_0}^{\beta_0} S^{(0)}_{Q,\gamma,k}(\omega)\,
d\omega\\ &  +\sum\limits_{\gamma \in \FI}
\sum\limits_{k=1}^\infty \ \int_{\beta_{k-1}}^{\beta_k}
S^{+}_{Q,\gamma,k}(\omega)\, d\omega ,
\end{split}
\end{equation*}
with
\begin{equation*}
\begin{split}
S^{-}_{Q,\gamma,k}(\omega) = & \ w_{C_k}(\omega)
(q_{k+1}-tQ)+w_{B_k}(\omega)(q_k-tQ)
+w_{A_+}(\omega)\max \{ q-tQ,0\} ,\\
S^{(0)}_{Q,\gamma,k}(\omega) = & \ w_{C_0} (\omega)(q+q^\prime
-tQ) +w_{A_+}(\omega) \max \{ q-tQ,0\}
+w_{B_-}(\omega) \max \{ q^\prime -tQ,0\}, \\
S^{+}_{Q,\gamma,k}(\omega) = & \
w_{C_{-k}}(\omega)(q^\prime_{k+1}-tQ)+w_{A_{-k}} (q^\prime_k -tQ)
+w_{B_-}(\omega)\max \{ q^\prime -tQ,0\}).
\end{split}
\end{equation*}
Here the formulas for the width of the strips are as in
\eqref{3.8}, and we take
\begin{equation}\label{4.1}
\begin{split}
& \alpha_k =\arctan \frac{a_k-1/Q}{q_k}, \quad \beta_k=\arctan
\frac{a_k^\prime+1/Q}{q_k^\prime} ,\\ & \alpha_\infty =\arctan
\frac{a}{q} ,\quad \beta_\infty =\arctan \frac{a^\prime}{q^\prime}
.
\end{split}
\end{equation}
\par
Taking into account the equalities
\begin{equation}\label{4.2}
q_{k+1} w_{C_k} +q_k w_{B_k}+qw_{A_+}=1=q_{k+1}^\prime
w_{C_{-k}}+q_k^\prime w_{A_{-k}}+qw_{B_-} ,
\end{equation}
\begin{equation}\label{4.3}
w_{C_k}+w_{B_k}+w_{A_+}=Q^{-1} =w_{C_{-k}}+w_{A_{-k}}+w_{B_-},
\end{equation}
and
\begin{equation*}
w_{C_0}(\omega)=Q^{-1}-w_{A_+}(\omega)-w_{B_-}(\omega),
\end{equation*}
we can write
\begin{equation}\label{4.4}
P_{I,Q}(t)=\sum\limits_{\gamma\in\FI}(
T^{(1)}_{Q,\gamma}+\cdots+T^{(5)}_{Q,\gamma})\, d\omega,
\end{equation}
with
\begin{equation*}
\begin{split}
T^{(1)}_{Q,\gamma} & =  \max \{ q-tQ,0\}
\int_{\alpha_\infty}^{\beta_0} w_{A_+}(\omega) \, d\omega ,\quad
T^{(2)}_{Q,\gamma}  = \max \{ q^\prime -tQ,0\}
\int_{\alpha_0}^{\beta_\infty} w_{B_-}(\omega) \, d\omega,\\
T^{(3)}_{Q,\gamma} & = (q+q^\prime -tQ) \int_{\alpha_0}^{\beta_0}
\left( \frac{1}{Q}
-w_{A_+}(\omega) -w_{B_-}(\omega) \right) d\omega ,\\
T^{(4)}_{Q,\gamma} & = \int_{\alpha_\infty}^{\alpha_0} \big(
(tQ-q)w_{A_+}(\omega)+1-t\big)\, d\omega ,\quad T^{(5)}_{Q,\gamma}
 = \int_{\beta_0}^{\beta_\infty} \big( (tQ-q)w_{B_-}(\omega)
+1-t\big) \, d\omega .
\end{split}
\end{equation*}
\par
Rewriting the terms in a convenient way we arrive at
\begin{equation*}
P_{I,Q}(t)=A_0+A_1+A_2+A_3,
\end{equation*}
where
\begin{equation*}
\begin{split}
A_0 & =(1-t)\sum\limits_{\gamma \in \FI}
(\beta_\infty-\alpha_\infty), \\
A_1 & =\sum\limits_{\gamma \in \FI} \Big( \max \{ q-tQ,0\}
+tQ-q\Big) \int_{\alpha_\infty}^{\beta_0} w_{A_+}(\omega) \,
d\omega \\ & =-\sum\limits_{\gamma \in \FI} \min\{ q-tQ,0\}
\int_{\alpha_\infty}^{\beta_0} w_{A_+}(\omega) \, d\omega ,
\\ A_2 & =\sum\limits_{\gamma \in \FI}\Big( \max \{ q^\prime -tQ,0\}
+tQ-q^\prime \Big) \int_{\alpha_0}^{\beta_\infty} w_{B_-}(\omega)
\, d\omega  \\ & =-\sum\limits_{\gamma \in \FI} \min \{ q^\prime
-tQ,0\} \int_{\alpha_0}^{\beta_\infty} w_{B_-}(\omega) \, d\omega
, \\
A_3 & =\sum\limits_{\gamma \in \FI}\ \int_{\alpha_0}^{\beta_0}
\bigg( \frac{q+q^\prime}{Q}-1 -q^\prime w_{A_+}(\omega)
-qw_{B_-}(\omega)\bigg) d\omega .
\end{split}
\end{equation*}

Remark first that $A_3=0$, as a result of
\begin{equation}\label{4.5}
\begin{split}
\frac{q+q^\prime}{Q}& -1 -q^\prime w_{A_+}(\omega)
-qw_{B_-}(\omega)
\\ &  =\frac{q+q^\prime}{Q}-1-q^\prime \left( \frac{1}{Q}+a-q\tan
\omega \right)-q\left( \frac{1}{Q}+q^\prime
\tan \omega -a^\prime \right) \\
& =a^\prime q-aq^\prime -1=0.
\end{split}
\end{equation}
\par
The next elementary statement will be repeatedly used.

\begin{lemma}\label{LEMMA6}
For any $\lambda,\mu \in \R$ we have, uniformly in $c\in [0,1]$ as
$h\rightarrow 0^+$,
\begin{equation}\label{4.6}
\begin{split}
\int_{\arctan c}^{\arctan(c+h)}(\lambda\tan \omega+\mu)\, d\omega
& =\left( \frac{h}{1+c^2}-\frac{h^2 c}{(1+c^2)^2}\right) (\lambda
c+\mu) \\ & \qquad + \frac{h^2 \lambda }{2(1+c^2)} +O(h^3 (\vert
\lambda\vert+\vert \mu \vert )),
\end{split}
\end{equation}
\begin{equation}\label{4.7}
\begin{split}
\int_{\arctan(c-h)}^{\arctan c} (\lambda\tan \omega +\mu) \,
d\omega & =\left( \frac{h}{1+c^2}+\frac{h^2 c}{(1+c^2)^2}\right)
(\lambda c+\mu) \\ & \qquad \qquad -\frac{h^2 \lambda}{2(1+c^2)}
+O(h^3 (\vert \lambda\vert+\vert \mu \vert )).
\end{split}
\end{equation}
\end{lemma}
\par
\begin{proof} Applying to our situation Taylor's formula
\begin{equation*}
\int_a^{a+\xi} f(x)\, dx=\xi f(a)+\frac{\xi^2}{2}\, f^\prime (a)
+O(\| f^{\prime \prime} \|_\infty \vert \xi \vert^3)
\end{equation*}
together with
\begin{equation}\label{4.8}
\xi=\arctan (c+h)-\arctan c=\frac{h}{1+c^2}-\frac{h^2
c}{(1+c^2)^2}+O(h^3),
\end{equation}
we get
\begin{equation*}
\begin{split}
\int_{\arctan c}^{\arctan (c+h)} \left( \tan
\omega+\frac{\mu}{\lambda}\right) d\omega & =\left(
\frac{h}{1+c^2}-\frac{h^2 c}{(1+c^2)^2}+O(h^3)\right) \left(
c+\frac{\mu}{\lambda}\right) \\
& \qquad \qquad +\frac{1}{2} \left( \frac{h}{1+c^2} -\frac{h^2
c}{(1+c^2)^2} +O(h^3)\right)^2 (1+c^2) +O(h^3) \\
& =\left( \frac{h}{1+c^2}-\frac{h^2 c}{(1+c^2)^2} \right) \left(
c+\frac{\mu}{\lambda}\right)+\frac{h^2}{2(1+c^2)} +O\left( h^3
\Big( \frac{\vert \mu\vert}{\vert \lambda\vert} +1\Big) \right) ;
\end{split}
\end{equation*}
whence \eqref{4.6} follows for $\lambda\neq 0$. The case $\lambda
=0$ is a direct consequence of \eqref{4.8}, while \eqref{4.7} is
derived from \eqref{4.6} by changing $h$ into $-h$.
\end{proof}
\par
This result will only be applied in cases where $\lambda c+\mu=0$.
We shall also use the following weaker form of \eqref{4.8}:
\begin{equation}\label{4.9}
\arctan (c+h)-\arctan c=\frac{h}{1+c^2}+O(h^2)
\end{equation}
\par
It remains to estimate $A_0$, $A_1$ and $A_2$. By \eqref{4.9} it
is immediate that
\begin{equation}\label{4.10}
A_0=(1-t)\sum\limits_{\gamma \in \FI}
 \left( \frac{1}{qq^\prime (1+\gamma^2)}
+O\Big( \frac{1}{q^2 q^{\prime 2}}\Big) \right) .
\end{equation}
This shows in conjunction with
\begin{equation}\label{4.11}
\sum\limits_{\gamma \in \FQ} \frac{1}{q^2 q^{\prime 2}} \ll
\sum\limits_{q=1}^Q \frac{1}{q^2} \sum\limits_{Q/2\leq q^\prime
\leq Q} \frac{1}{q^{\prime 2}} \ll \frac{1}{Q} \sum\limits_{q=1}^Q
\frac{1}{q^2} \ll \frac{1}{Q}
\end{equation}
and with the subsequent Lemma \ref{L4.4} that
\begin{equation}\label{4.12}
A_0 =c_I (1-t) +O_\delta ( E_{c,c_1,\delta}(Q)).
\end{equation}

\begin{lemma}\label{L4.4}
Let $c,c_1>0$ such that $c+c_1<1$. Then for any interval
$I\subseteq [0,1]$ with $\vert I\vert \asymp Q^{-c}$ and $\delta
>0$
\begin{equation*}
\sum\limits_{\gamma \in \FI} \frac{1}{qq^\prime (1+\gamma^2)} =c_I
+O_\delta ( E_{c,c_1,\delta}(Q)) .
\end{equation*}
\end{lemma}
\par
\begin{proof}
We decompose the sum above as $S_1+S_2$, according to whether
$q^\prime >q$ or $q>q^\prime$. Thus we can write
\begin{equation}\label{4.13}
S_1 =\sum\limits_{q=1}^Q \ \sum\limits_{\substack{q^\prime \in
\II:=( \max \{ Q-q,q\} ,Q] \\ a\in \JJ:=qI
\\ -aq^\prime =1 \hspace{-6pt} \pmod{q}}} \hspace{-20 pt}
f_q (q^\prime,a),
\end{equation}
where we put
\begin{equation*}
f_q (q^\prime ,a)=\frac{1}{qq^\prime ( 1+a^2/q^2)}, \qquad a\in
\II,\ q^\prime \in \JJ,\ q\in [1,Q-1].
\end{equation*}
The inclusion $\II \subseteq ( \frac{Q}{2},Q]$ gives
\begin{equation*}
\begin{split}
& 0\leq f_q (q^\prime,a) =\frac{q}{q^\prime(q^2+a^2)} \leq
\frac{1}{qq^\prime} \leq \frac{2}{qQ}  ;\\
& 0\leq Df_q (q^\prime,a) =\left| \frac{\partial f_q}{\partial
q^\prime}\, (q^\prime,a)\right| +\left| \frac{\partial
f_q}{\partial a}\, (q^\prime,a)\right| =
\frac{q}{q^\prime(q^2+a^2)} \left(
\frac{1}{q^\prime}+\frac{2a}{q^2+a^2}\right) \\
& \qquad \qquad \qquad \leq \frac{1}{qq^\prime} \left(
\frac{2}{Q}+\frac{2}{q}\right) \leq \frac{4}{q^2 q^\prime} \leq
\frac{8}{q^2 Q}  .
\end{split}
\end{equation*}
Applying Lemma \ref{LEMMA2} with $T=[Q^{c_1}]$, the inner sum in
\eqref{4.13} can be expressed as
\begin{equation*}
\begin{split}
& \frac{\varphi(q)}{q}\hspace{-2pt} \int_\II
\frac{dq^\prime}{qq^\prime} \int_{qI} \frac{da}{1+a^2/q^2} +
O_\delta \left( Q^{2c_1} q^{1/2+\delta} \frac{1}{qQ} + Q^{c_1}
q^{3/2+\delta}
\frac{1}{q^2 Q} +\frac{q^2 Q^{-c}}{Q^{c_1}q^2 Q} \right) \\
& =c_I \, \frac{\varphi(q)}{q}\, V(q) +O_\delta (Q^{2c_1-1}
q^{1/2+\delta} +Q^{c_1-1} q^{-1/2+\delta} +Q^{-1-c-c_1}),
\end{split}
\end{equation*}
where
\begin{equation*}
V(q) =\int_\II \frac{dq^\prime}{qq^\prime} =\frac{1}{q}\ln
\frac{Q}{\max \{ q,Q-q\}}  , \qquad q\in (0,Q].
\end{equation*}
The function
\begin{equation*}
W(x):=\begin{cases} \vspace{3pt} \displaystyle \frac{1}{x}
\ln \frac{1}{\max \{ x,1-x\}} & \mbox{\rm if $x\in (0,1];$}\\
1 & \mbox{\rm if $x=0$,}
\end{cases}
\end{equation*}
is bounded and has finite total variation on $[0,1]$, hence
\begin{equation*}
M:=\sup\limits_{x\in [0,1]} \vert W(x)\vert+\int_0^1 \vert
W^\prime (x)\vert \, dx\ =O(1).
\end{equation*}
Since $V(Qx)=\frac{W(x)}{Q}$, Lemma \ref{LEMMA1} yields
\begin{equation*}
\begin{split}
\sum\limits_{q=1}^Q \frac{\varphi (q)}{q}\ V(q) &
=\frac{1}{\zeta(2)} \int_0^Q V(q)\, dq+O\left( \ln Q \bigg(
\sup\limits_{q\in (0,Q]} \vert V(q)\vert +\int_0^Q \vert
V^\prime (q)\vert\, dq\bigg) \right) \\
& = \frac{1}{\zeta(2)}\int_0^1 W(x)\, dx+O(Q^{-1}\ln Q).
\end{split}
\end{equation*}
Hence
\begin{equation}\label{4.14}
\begin{split} S_1  & =c_I \sum\limits_{q=1}^{Q} \frac{\varphi(q)}{q}
\, V(q) +O_\delta ( E_{c,c_1,\delta}(Q))
 =\frac{c_I}{\zeta(2)} \int_0^1 W(x)\, dx+O_\delta (
E_{c,c_1,\delta}(Q)).
\end{split}
\end{equation}
Using a familiar identity of Euler (cf. formula (1.8) in
\cite{Lew}) we find that
\begin{equation*}
\begin{split}
\int_0^1 W(x)\, dx & =-\int_0^{1/2} \frac{\ln (1-x)}{x}\
dx-\int_{1/2}^1 \frac{\ln x}{x}\ dx \\ & = -\int_0^{1/2} \frac{\ln
(1-x)}{x}\ dx -\frac{\ln^2 2}{2} =\frac{\ln^2
2}{2}+\frac{\zeta(2)}{2}-\frac{\ln^2 2}{2} =\frac{\zeta(2)}{2} ,
\end{split}
\end{equation*}
which we combine with \eqref{4.14} to get
\begin{equation}\label{4.15}
S_1=\frac{c_I}{2}+O_\delta ( E_{c,c_1,\delta}(Q)).
\end{equation}
\par
Finally we employ
\begin{equation}\label{4.16}
\frac{1}{1+\gamma^{\prime 2}}=\frac{1}{1+\gamma^2}+O(\gamma^\prime
-\gamma) =\frac{1}{1+\gamma^2}+O\left( \frac{1}{qq^\prime}\right)
\end{equation}
and \eqref{4.11} to write
\begin{equation*}
S_2 =\sum\limits_{\substack{\gamma \in \FI \\ q>q^\prime}}
\frac{1}{qq^\prime (1+\gamma^2)} .
\end{equation*}
Using \eqref{4.16} and \eqref{4.11} we see that
\begin{equation*}
S_2 =\sum\limits_{\substack{\gamma \in \FI \\
q>q^\prime}} \frac{1}{qq^\prime
(1+\gamma^2)}=\sum\limits_{q^\prime=1}^Q\
\sum\limits_{\substack{q\in (\max \{ Q-q^\prime ,q^\prime \} ,Q]
\\ a^\prime \in q^\prime I \\ a^\prime q=1\hspace{-6pt}\pmod{q^\prime}}}
\frac{1}{qq^\prime ( 1+a^{\prime 2}/q^{\prime 2})} .
\end{equation*}
Changing $a^\prime$ to $q^\prime -a^\prime$, reversing the roles
of $q$ and $q^\prime$, and using
\begin{equation*}
\int_{q^\prime (1-I)} \frac{dx}{1+(1-x/q^\prime)^2} =c_I q^\prime,
\end{equation*}
it follows that $S_2$ is given by the same expression as in
\eqref{4.15}.
\end{proof}
\par
Next we estimate $A_1$ and find, taking $c=\frac{a+1/Q}{q}$,
$h=\frac{1}{qQ}$, $\lambda=-q$, $\mu=a+\frac{1}{Q}$ in
\eqref{4.7}, that
\begin{equation}\label{4.17}
\begin{split}
\int_{\alpha_\infty}^{\beta_0} & w_{A_+}(\omega)\, d\omega =
\int_{\arctan \frac{a}{q}}^{\arctan \frac{a+1/Q}{q}}
\left( \frac{1}{Q}+a-q\tan \omega \right) d\omega \\
& =\frac{1}{2qQ^2 \big( 1+ (a+1/Q)^2)/q^2\big)}+O\bigg(
\frac{1}{q^2 Q^3} \bigg) \\ & = \frac{1}{2qQ^2(1+\gamma^2)}
+O\left( \frac{1}{q^2Q^3}\right).
\end{split}
\end{equation}
Since
\begin{equation*}
\sum\limits_{\gamma \in \FQ} \frac{Q}{q^2Q^3} \leq \frac{1}{Q^2}
\sum\limits_{q=1}^Q \frac{\varphi(q)}{q^2} \ll \frac{\ln Q}{Q^2}
=O(Q^{-1}),
\end{equation*}
we infer from \eqref{4.17} and the definition of $A_1$ that
\begin{equation}\label{4.18}
\begin{split}
A_1+O(Q^{-1}) & =\sum\limits_{\substack{\gamma \in \FI \\ q\leq
tQ}} \frac{tQ-q}{2qQ^2(1+a^2/q^2)} \\ & =\sum\limits_{1\leq q\leq
tQ} \frac{tQ-q}{2qQ^2}
\sum\limits_{\substack{Q-q<q^\prime \leq Q \\ a\in qI \\
-a q^\prime =1\hspace{-6pt} \pmod{q}}} \frac{1}{1+a^2/q^2}  .
\end{split}
\end{equation}
\par
Applying Lemma \ref{LEMMA2} to $\II=(Q-q,Q]$, $\JJ=qI$, $f_q
(q^\prime ,a)=\frac{1}{1+a^2/q^2}$ for which $\| f_q\|_\infty \leq
1$ and $\| Df_q \|_\infty \leq \frac{2}{q}$, and taking
$T=[Q^{c_1}]$, the inner sum above becomes
\begin{equation*}
\begin{split}
& \frac{\varphi(q)}{q^2}\, q\int_{qI} \frac{da}{1+a^2/q^2}
+O_\delta \left( Q^{2c_1} q^{1/2+\delta} +Q^{c_1} q^{3/2+\delta}
\frac{1}{q}+\frac{q^2\vert I\vert }{Q^{c_1}q} \right) \\
& \qquad =c_I \varphi(q)+O_\delta (Q^{2c_1}q^{1/2+\delta}
+Q^{-c-c_1} q),
\end{split}
\end{equation*}
which inserted back into \eqref{4.18} gives that $A_1+O(Q^{-1})$
may be written as
\begin{equation}\label{4.19}
\begin{split}
\frac{c_I}{2Q^2} \sum\limits_{1\leq q\leq tQ} &
\frac{\varphi(q)}{q} (tQ-q) +O_\delta \left(\ \sum\limits_{q=1}^Q
\frac{Q}{qQ^2} \big(
Q^{2c_1}q^{1/2+\delta} +Q^{-c-c_1} q \big) \right) \\
& =\frac{c_I}{2Q^2} \sum\limits_{1\leq q\leq tQ}
\frac{\varphi(q)}{q} (tQ-q)+O_\delta (E_{c,c_1,\delta}(Q)) .
\end{split}
\end{equation}
Applying now Lemma \ref{LEMMA1} to the main term above with
$V(q)=tQ-q$, $q\in [1,tQ]$, we find that
\begin{equation}\label{4.20}
\begin{split}
A_1 & =\frac{c_I}{2Q^2 \zeta(2)} \int_0^{tQ}
(tQ-q)\, dq+O_\delta ( Q^{-1+\delta}+E_{c,c_1,\delta}(Q)) \\
& =\frac{c_I}{2\zeta(2)} \int_0^t (t-x)\, dx +O_\delta (
E_{c,c_1,\delta}(Q)) =\frac{c_I t^2}{4\zeta(2)}  + O_\delta (
E_{c,c_1,\delta}(Q)) .
\end{split}
\end{equation}
\par
In a similar way we find
\begin{equation*}
A_2 =\frac{c_I t^2}{4\zeta(2)} +O_\delta ( E_{c,c_1,\delta}(Q)) ,
\end{equation*}
and therefore
\begin{equation*}
P_{I,Q}(t)=\left( 1-t+\frac{t^2}{2\zeta(2)} \right) c_I +O_\delta
(E_{c,c_1,\delta}(Q)) ,
\end{equation*}
which proves Proposition \ref{PROP3}.

\section{The case $t>2$}
In this section we shall evaluate the contribution of the
integrals on $[\alpha_k,\alpha_{k-1}]$ in \eqref{3.7} when $k\geq
1$ and $t>2$. In this situation there is a unique nonnegative
integer, given by
\begin{equation*}
K=K(\gamma,t)=\left[ \frac{tQ-q^\prime}{q} \right] \geq 0,
\end{equation*}
for which
\begin{equation}\label{5.1}
q_K \leq tQ< q_{K+1} .
\end{equation}
When $t\geq 2$ it follows that $K\geq 1$, and we prove

\begin{proposition}\label{PROP4}
Suppose $I$ is a subinterval of $[0,1]$ of size $\vert I\vert
\asymp Q^{-c}$ for some $0<c <1$. Then for every $c_1>0$ with
$c+c_1<1$ and $\delta >0$
\begin{equation*}
P_{I,Q}(t) =\frac{c_I}{\zeta(2)}
 \int_0^1 \psi(x,t)\, dx
+O_\delta ( E_{c,c_1,\delta}(Q)) \qquad (Q\rightarrow \infty),
\end{equation*}
with $\psi$ as in Theorem {\em \ref{T1.1}} and $E_{c,c_1,\delta}$
as in Proposition {\em \ref{PROP3}}. The estimate is uniform in
$t$ on compacts of $(2,\infty)$.
\end{proposition}
\par
Next, $\alpha_k$ and $\beta_k$ will be as in \eqref{4.1} and the
widths $w$ as in \eqref{3.8}. Since $t>2$, then $q+q^\prime <tQ$
and the second sum in \eqref{3.7} is zero.
\par
In the beginning we fix $\gamma =\frac{a}{q}\in \FI$ and estimate
\begin{equation*}
S_2(\gamma,t) :=\sum\limits_{k=K+1}^\infty
\int_{\alpha_k}^{\alpha_{k-1}} \Big( q_{k+1} w_{C_k}(\omega)+q_k
w_{B_k} (\omega) -tQ \big( w_{C_k}(\omega)+w_{B_k} (\omega)\big)
\Big)\, d\omega .
\end{equation*}
Using \eqref{4.2} and \eqref{4.3}, taking $c=\frac{a}{q}$,
$h=\frac{a_K-1/Q}{q_K}-\frac{a}{q}=\frac{1-q/Q}{qq_K}$, $\lambda
=q$, $\mu=-a$ in \eqref{4.6}, and also owing to
\begin{equation*}
\alpha_K-\alpha_\infty =\arctan \frac{a_K-1/Q}{q_K} -\arctan
\frac{a}{q}= \frac{1-q/Q}{qq_K( 1+a^2/q^2)} +O\left( \frac{1}{q^2
q_1^2}\right),
\end{equation*}
we infer that
\begin{equation*}
\begin{split}
S_2(\gamma,t) & =\int_{\alpha_\infty}^{\alpha_K} \left( 1
-qw_{A_+}(\omega)-tQ \Big( \frac{1}{Q}-w_{A_+}(\omega)\Big)
\right)\, d\omega \\
& =\int_{\alpha_\infty}^{\alpha_K} \left( q^2 \tan \omega+1-\Big(
\frac{1}{Q}+a\Big) q\right)\, d\omega
-tQ\int_{\alpha_\infty}^{\alpha_K}
(q\tan \omega-a)\, d\omega \\
& =\left( 1-\frac{q}{Q}\right) (\alpha_K-\alpha_\infty)+
(q-tQ)\int\limits_{\alpha_\infty}^{\alpha_K}
(q\tan \omega-a)\, d\omega \\
 & =\left( 1-\frac{q}{Q}\right)
(\alpha_K-\alpha_\infty) +(q-tQ) \left( \frac{(1-q/Q)^2 q}{2q^2
q_K^2 (1+\gamma^2)} +O\Big( \frac{Q}{q^2 q_K^3} \Big) \right)
\\ & =\frac{( 1-q/Q)^2}{qq_K (1+\gamma^2)} +(q-tQ)
\frac{( 1-q/Q)^2}{2qq_K^2 (1+\gamma^2)}
+O\left( \frac{1}{q^2 q_1^2}\right) \\
& =\frac{( 1-q/Q)^2 (q+2q_K-tQ)}{2qq_K^2(1+\gamma^2)} +O\left(
\frac{1}{q^2 q_1^2} \right) .
\end{split}
\end{equation*}
\par
On the other hand, taking $c=\frac{a_{K-1}-1/Q}{q_{K-1}}$,
$h=\frac{a_{K-1}-1/Q}{q_{K-1}}-\frac{a_K-1/Q}{q_K}
=\frac{1-q/Q}{q_{K-1}q_K}$, $\lambda=-q_{K-1}$,
$\mu=a_{K-1}-\frac{1}{Q}$ in \eqref{4.7}, and also using
\begin{equation*}
0\leq \frac{1}{1+c^2}-\frac{1}{1+\gamma^2} \ll
\frac{a_{K-1}-\frac{1}{Q}}{q_{K-1}} -\frac{a}{q} \leq
\frac{1}{qq_{K-1}} ,
\end{equation*}
we estimate
\begin{equation*}
\begin{split}
\int_{\alpha_K}^{\alpha_{K-1}} w_{C_k}(\omega)\, d\omega &
=\int_{\alpha_K}^{\alpha_{K-1}}  \left(
a_{K-1}-\frac{1}{Q}-q_{K-1}\tan \omega \right)\, d\omega \\ &
=\frac{q_{K-1}( 1-q/Q)^2}{2q_{K-1}^2 q_K^2 (1+c^2)}+O\left(
\frac{q_{K-1}}{q_{K-1}^3 q_K^3}\right)
\\ & =\frac{( 1-q/Q)^2}{2q_{K-1} q_K^2
(1+\gamma^2)}+O\left( \frac{1}{q_{K-1}^2
q_K^3}+\frac{1}{qq_{K-1}}\cdot \frac{1}{q_{K-1}q_K^2} \right)
\\ & =\frac{( 1-q/Q)^2}{2q_{K-1} q_K^2
(1+\gamma^2)}+O\left( \frac{1}{qq_{K-1}^2 q_K^2}\right) .
\end{split}
\end{equation*}
\par
Using also $0<q_{K+1}-tQ\leq q$, this gives whenever $t>2$ (so
$K\geq 1$)
\begin{equation*}
\begin{split}
S_1(\gamma,t) & :=(q_{K+1}-tQ) \int_{\alpha_K}^{\alpha_{K-1}}
 w_{C_k}(\omega)\, d\omega  \\  &= \frac{(q_{K+1}-tQ)(
1-q/Q)^2}{2q_{K-1}^2 q_K^2 (1+\gamma^2)}+O\left(
\frac{1}{q^{\prime 2}(q+q^\prime)^2}\right) .
\end{split}
\end{equation*}
\par
Since $t>2$, the sum of integrals on $[\alpha_k,\alpha_{k-1}]$ in
\eqref{3.7} becomes
\begin{equation*}
P^+_{I,Q}(t):=\sum\limits_{\gamma \in \FI} \big(
S_1(\gamma,t)+S_2(\gamma,t)\big) .
\end{equation*}
Making use of
\begin{equation*}
\sum\limits_{\gamma \in \FQ} \frac{1}{q^{\prime 2}(q+q^\prime)^2}
\leq \sum\limits_{q^\prime =1}^Q \frac{1}{q^{\prime 2}}
\sum\limits_{q=Q-q^\prime}^Q \frac{1}{(q^\prime +q)^2} \leq
\sum\limits_{q^\prime =1}^Q \frac{1}{q^{\prime 2}}
\sum\limits_{k=Q+1}^\infty \frac{1}{k^2}\ll \frac{1}{Q}
\end{equation*}
and of
\begin{equation*}
\begin{split}
& \frac{(1-q/Q)^2 (q_{K+1}-tQ)}{2q_{K-1} q_K^2 (1+\gamma^2)}
+\frac{(1-q/Q)^2
(q+2q_K-tQ)}{2qq_K^2(1+\gamma^2)} \\
& \qquad =\frac{( 1-q/Q)^2\big( q(q_{K+1}+q_{K-1})+2q_{K-1}q_K
-tQ(q+q_{K-1})\big)}{2qq_{K-1}q_K^2
(1+\gamma^2)}   \\
 & \qquad =\frac{(1-q/Q)^2 (2qq_K+2q_{K-1}q_K-tQq_K )
 }{2qq_{K-1} q_K^2 (1+\gamma^2)} =
\frac{(1-q/Q)^2 (2q_K-tQ)}{2qq_{K-1} q_K(1+\gamma^2)} ,
\end{split}
\end{equation*}
we find
\begin{equation*}
P^+_{I,Q}(t)=\sum\limits_{\gamma \in \FI} \frac{(1-q/Q)^2
(2q_K-tQ)}{2qq_{K-1} q_K (1+\gamma^2)} +O(Q^{-1}).
\end{equation*}
\par
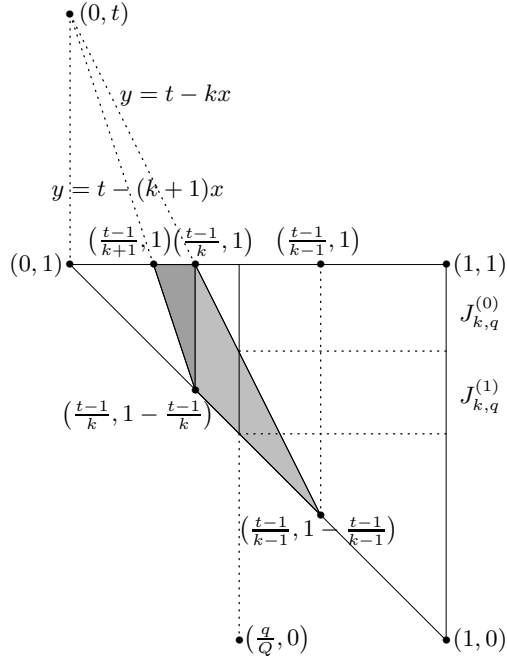
\begin{figure}[ht]
\begin{center}
\unitlength 0.5mm
\begin{picture}(100,170)(0,0)
\texture{ccc 0}
\shade\path(22.2222,100)(33.3333,100)(33.3333,66.6666)(22.2222,100)

\texture{c 0}
\shade\path(33.3333,100)(33.3333,66.6666)(66.6666,33.3333)(33.3333,100)

\path(0,100)(100,100)(100,0)(0,100)
\path(22.2222,100)(33.3333,66.6666)
\path(33.3333,100)(66.6666,33.3333) \path(45,55)(45,100)
\dottedline{2}(66.6666,100)(66.6666,33.3333)
\dottedline{2}(45,55)(45,0)
\dottedline{2}(0,100)(0,166.6666)(22.2222,100)
\dottedline{2}(0,166.6666)(33.3333,100)
\dottedline{2}(45,77)(100,77) \dottedline{2}(45,55)(100,55)

\put(0,166.6666){\makebox(0,0){{\tiny $\bullet$}}}
\put(0,100){\makebox(0,0){{\tiny $\bullet$}}}
\put(22.2222,100){\makebox(0,0){{\tiny $\bullet$}}}
\put(33.3333,100){\makebox(0,0){{\tiny $\bullet$}}}
\put(66.6666,100){\makebox(0,0){{\tiny $\bullet$}}}
\put(100,100){\makebox(0,0){{\tiny $\bullet$}}}
\put(33.3333,66.6666){\makebox(0,0){{\tiny $\bullet$}}}
\put(66.6666,33.3333){\makebox(0,0){{\tiny $\bullet$}}}
\put(100,0){\makebox(0,0){{\tiny $\bullet$}}}
\put(45,0){\makebox(0,0){{\tiny $\bullet$}}}

\put(-9,100){\makebox(0,0){{\small $(0,1)$}}}
\put(9,166.66){\makebox(0,0){{\small $(0,t)$}}}
\put(109,100){\makebox(0,0){{\small $(1,1)$}}}
\put(109,0){\makebox(0,0){{\small $(1,0)$}}}
\put(55,0){\makebox(0,0){{\small $\big( \frac{q}{Q},0\big)$}}}
\put(18,60){\makebox(0,0){{\small $\big(
\frac{t-1}{k},1-\frac{t-1}{k}\big)$}}}
\put(66,29){\makebox(0,0){{\small $\big(
\frac{t-1}{k-1},1-\frac{t-1}{k-1}\big)$}}}
\put(28,145){\makebox(0,0){{\small $y=t-kx$}}}
\put(18,120){\makebox(0,0){{\small $y=t-(k+1)x$}}}

\put(66,106){\makebox(0,0){{\small $\big(
\frac{t-1}{k-1},1\big)$}}} \put(38,106){\makebox(0,0){{\small
$\big( \frac{t-1}{k},1\big)$}}} \put(16,106){\makebox(0,0){{\small
$\big( \frac{t-1}{k+1},1\big)$}}}
\put(109,87){\makebox(0,0){{\small $J_{k,q}^{(0)}$}}}
\put(109,65){\makebox(0,0){{\small $J_{k,q}^{(1)}$}}}
\end{picture}
\end{center}
\caption{The set $\Omega_k \cap \TT$} \label{Figure6}
\end{figure}

Next for each integer $k\geq 1$ consider the sets
\begin{equation*}
\Omega_k=\left\{ (x,y)\in \R^2 \, ;\, \bigg[ \frac{t-y}{x} \bigg]
=k\right\} \quad \mbox{\rm and} \quad
 I_k=\bigg[ \frac{t-1}{k},\frac{t-1}{k-1}\bigg) \cap [0,1),
\end{equation*}
and for $q\in QI_k$ and $k\geq 1$, respectively $k\geq 2$, the
intervals (see Figure \ref{Figure6})
\begin{equation*}
\begin{split}
& J_{k,q}^{(0)}=\bigg( t-\frac{kq}{Q},1\bigg]=\left\{
\frac{q^\prime}{Q} \, ;\, \bigg( \frac{q}{Q},\frac{q^\prime}{Q}
\bigg) \in \Omega_{k-1} \cap \TT \right\} \subseteq
\bigg( 1-\frac{q}{Q},1\bigg] , \\
& J_{k,q}^{(1)} =\bigg( 1-\frac{q}{Q},t-\frac{kq}{Q} \bigg]
=\left\{ \frac{q^\prime}{Q}\, ;\, \bigg(
\frac{q}{Q},\frac{q^\prime}{Q} \bigg) \in \Omega_k \cap \TT
\right\} \subseteq \bigg( 1-\frac{q}{Q},1\bigg] .
\end{split}
\end{equation*}
Note that $\vert QJ^{(0)}_{k,q}\vert,\vert QJ^{(1)}_{k,q}\vert
<q$, that $\min \{ k\, ;\, \vert \Omega_k \cap \TT \vert
>0\}=[t]-1\geq 1$, and that $\vert I_k\vert =0$ unless $k\geq
[t]\geq 2$.

We also consider the function
\begin{equation*}
\begin{split}
Q\Omega_k \times [0,q] \ni (q,q^\prime,a) \mapsto
f_k(q,q^\prime,a) & =\frac{(1-q/Q)^2 (2q_k-tQ)}{2qq_k
q_{k-1}(1+\gamma^2)} \\ & = \frac{(1-q/Q)^2
(2q^\prime+2kq-tQ)}{2q(q^\prime +kq)( q^\prime+(k-1)q)
(1+a^2/q^2)} .
\end{split}
\end{equation*}
Using the one-to-one correspondence between the primitive integer
points in $Q(\Omega_k \cap \TT)$ and the set of consecutive Farey
fractions $\gamma=\frac{a}{q}$ and
$\gamma^\prime=\frac{a^\prime}{q^\prime}$ in $\FQ$ with
$[\frac{Q-q^\prime}{q}]=k$, we derive using the summation method
described in Section 2 that
\begin{equation*}
P^+_{I,Q}(t) = \sum\limits_{k=1}^\infty \sum_{\substack{\gamma
=\frac{a}{q} \in \FI \\ (q,q^\prime)\in Q(\Omega_k \cap \TT)}}
\hspace{-12pt} f_k(q,q^\prime,a) =\sum\limits_{k=2}^\infty
\sum\limits_{q\in QI_k} \big( S_k(q)+T_k(q)\big) +O(Q^{-1}),
\end{equation*}
with
\begin{equation*}
S_k(q)=\hspace{-6pt} \sum\limits_{\substack{q^\prime \in
QJ_{k,q}^{(1)} ,\, a\in qI \\ -aq^\prime =1\hspace{-6pt}
\pmod{q}}} \hspace{-10pt} f_k(q,q^\prime,a),\quad T_k(q)=
\sum\limits_{\substack{q^\prime \in QJ_{k,q}^{(0)} ,\, a\in qI \\
-aq^\prime =1\hspace{-6pt} \pmod{q}}} \hspace{-10pt} f_{k-1}
(q,q^\prime,a) .
\end{equation*}
We aim to estimate $S_k(q)$ and $T_k(q)$ applying Lemma
\ref{LEMMA2} to the intervals $\II=QJ_{k,q}^{(1)}$, $\JJ =qI$ and
the function $f=f_k(q,\cdot,\cdot)$, and respectively to
$\II=QJ_{k,q}^{(0)}$, $\JJ=qI$ and $f=f_{k-1}(q,\cdot,\cdot)$. For
$(q,q^\prime)\in Q(\Omega_k \cap \TT)$ we have $q_k\leq
tQ<q_{k+1}$, or equivalently
\begin{equation*}
q_{k-1} <2q_k-tQ\leq q_k .
\end{equation*}
As a result, we see that (here $k\geq 2$)
\begin{equation*}
\begin{split}
\| f_k(q,\cdot,\cdot)\|_{\infty} & \leq \sup\limits_{q^\prime \in
QJ_{k,q}^{(1)}} \frac{q_k}{qq_k q_{k-1}} \leq
\sup\limits_{q^\prime \in (Q-q,Q]} \frac{1}{q(q+q^\prime)}
<\frac{1}{qQ} ,\\
\| f_{k-1}(q,\cdot,\cdot)\|_{\infty} & \leq \sup\limits_{q^\prime
\in QJ_{k,q}^{(0)}} \frac{q_{k-1}}{qq_{k-1}q_{k-2}} \leq
\sup\limits_{q^\prime
>(t-2)Q} \frac{1}{qq^\prime} \leq \frac{1}{(t-2)qQ}
\ll_t \frac{1}{qQ} .
\end{split}
\end{equation*}
The last estimate holds without the factor $\frac{1}{t-2}$
whenever $k>[t]$. In the remainder of this section we will simply
write $\frac{1}{t-2} \ll 1$ with the understanding that this holds
uniformly in $t$ on compacts of $(2,\infty)$.
\par
We also need to estimate the $L_\infty$-norm of $Df_k$. It is
easily seen that
\begin{equation*}
\begin{split}
\left\| \frac{\partial f_k}{\partial a} (q,\cdot,\cdot)
\right\|_{\infty} & \leq \frac{2q}{q^2}  \|
f_k(q,\cdot,\cdot)\|_{\infty} \ll \frac{1}{q^2 Q} , \\
\left\| \frac{\partial f_{k-1}}{\partial a}(q,\cdot,\cdot)
\right\|_{\infty} & \leq \frac{2q}{q^2} \| f_{k-1}
(q,\cdot,\cdot)\|_{\infty} \ll \frac{1}{(t-2)q^2 Q} \ll
\frac{1}{q^2 Q} ,\\
\left\| \frac{\partial f_k}{\partial q^\prime} (q,\cdot,\cdot)
\right\|_{\infty} & \leq \frac{1}{2q} \sup\limits_{q^\prime \in
QJ_{k,q}^{(1)}} \frac{\vert
2q_kq_{k-1}-(2q_k-tQ)(q_k+q_{k-1})\vert}{q_k^2 q_{k-1}^2} \\
& \ll \sup\limits_{q^\prime \in QJ_{k,q}^{(1)}} \left(
\frac{1}{qq_k q_{k-1}} +\frac{q_k(q_k+q_{k-1})}{qq_k^2q_{k-1}^2}
\right) \ll \sup\limits_{q^\prime \in QJ_{k,q}^{(1)}}
\frac{q_k+q_{k-1}}{qq_kq_{k-1}^2} \\
& \ll \sup\limits_{q^\prime \in QJ_{k,q}^{(1)}}
\frac{1}{qq_{k-1}^2} \leq \sup\limits_{q^\prime \in (Q-q,Q]}
\frac{1}{q(q+q^\prime)^2} <\frac{1}{qQ^2} \leq \frac{1}{q^2 Q} ,
\end{split}
\end{equation*}
and similarly
\begin{equation*}
\left\| \frac{\partial f_{k-1}}{\partial q^\prime}
(q,\cdot,\cdot)\right\|_{\infty} \ll \frac{1}{q}
\sup\limits_{q^\prime \in [(t-2)Q,Q]} \frac{1}{q^{\prime 2}} \leq
\frac{1}{(t-2)^2 qQ^2} \ll_t \frac{1}{q Q^2} \leq \frac{1}{q^2 Q}
.
\end{equation*}
Applying now Lemma \ref{LEMMA2} to this situation with $T=[
Q^{c_1}]$, where $0<c_1<\frac{1}{2}$ is to be determined later, we
approximate $S_k(q)+T_k(q)$ within error
\begin{equation*}
\begin{split}
E_k(q) & \ll_\delta Q^{2c_1}q^{1/2+\delta}\frac{1}{Qq}+Q^{c_1}
q^{3/2+\delta} \frac{1}{Qq^2}+Q^{-c-c_1}q^2\frac{1}{Qq^2}
\\ & \ll_\delta Q^{2c_1-1} q^{-1/2+\delta} +Q^{-1-c-c_1}
\end{split}
\end{equation*}
by
\begin{equation*}
\begin{split}
& \frac{\varphi(q)}{q^2} \iint_{QJ_{k,q}^{(1)}\times qI} f_k
(q,q^\prime,a)\, dq^\prime \, da+\frac{\varphi(q)}{q^2}
\iint_{QJ_{k,q}^{(0)}\times qI} f_{k-1}(q,q^\prime,a) \,
dq^\prime\, da \\ & \qquad =c_I\, \frac{\varphi(q)}{q}\cdot
\frac{(1-q/Q)^2}{2q}\, W_k (q),
\end{split}
\end{equation*}
where $c_I$ is as in \eqref{3.10} and
\begin{equation*}
\begin{split}
W_k (q) & =\int_{QJ_{k,q}^{(1)}} g_k (q,q^\prime)\, dq^\prime
+\int_{QJ_{k,q}^{(0)}} g_{k-1} (q,q^\prime)\, dq^\prime ,\\
g_k(q,q^\prime) & =\frac{2q_k-tQ}{q_k q_{k-1}}  , \quad
(q,q^\prime)\in Q(\Omega_k \cap \TT).
\end{split}
\end{equation*}
\par
By a direct computation we find that
\begin{equation*}
\begin{split}
W_k (q) =W(q) & =\int_{Q-q}^{tQ-kq} \frac{2q_k-tQ}{q_k q_{k-1}} \,
dq^\prime +\int_{tQ-kq}^{Q}
\frac{2q_{k-1}-tQ}{q_{k-1} q_{k-2}} \,  dq^\prime \\
& =\int_{Q+(k-1)q}^{tQ} \frac{2y-tQ}{y(y-q)}\, dy
+\int_{tQ-q}^{Q+(k-1)q} \frac{2y-tQ}{y(y-q)}\, dy \\
& =\int_{tQ-q}^{tQ} \frac{2y-tQ}{y(y-q)}\, dy= \left( 2\ln
(y-q)-\frac{tQ}{q} \ln \frac{y-q}{y} \right) \Bigg|_{y=tQ-q}^{tQ} \\
& =2\ln \frac{tQ-q}{tQ-2q} -\frac{tQ}{q}  \ln
\frac{(tQ-q)^2}{tQ(tQ-2q)}
\end{split}
\end{equation*}
is independent of $k$. Since the error terms sum up to
\begin{equation*}
\begin{split}
\sum\limits_{k=2}^\infty \sum\limits_{q\in QI_k} (Q^{2c_1-1}
q^{-1/2+\delta} +Q^{-1-c-c_1} ) & =Q^{2c_1-1} \sum\limits_{q=1}^Q
q^{-1/2+\delta} +Q^{1-1-c-c_1} \ll E_{c,c_1,\delta}(Q) ,
\end{split}
\end{equation*}
we arrive at
\begin{equation}\label{5.2}
P^+_{I,Q}(t)= c_I \sum\limits_{q=1}^Q \frac{\varphi(q)}{q}\, V(q)
+O_\delta ( E_{c,c_1,\delta}(Q)) ,
\end{equation}
with
\begin{equation*}
V(q)=\frac{( 1-q/Q)^2}{2q}  W(q) ,\quad q\in (0,Q].
\end{equation*}
\par
For $t>2$ consider the function
\begin{equation*}
f_t(x)=\frac{\psi (x,t)}{2}=\frac{(1-x)^2}{2x} \left( 2\ln \Big(
1+\frac{x}{t-2x}\Big)-\frac{t}{x}\ln \Big(
1+\frac{(t-x)^2}{t(t-2x)}\Big)\right), \quad x\in (0,1].
\end{equation*}
Using the Taylor series of the logarithm we obtain for small $x$
\begin{equation*}
\begin{split}
f_t(x) & =\frac{(1-x)^2}{2x} \left(
\frac{2x}{t-2x}-\frac{x^2}{(t-2x)^2}-\frac{t}{x}\cdot
\frac{x^2}{t(t-2x)}+O(x^3)\right) \\
& =(1-x)^2 \left(
\frac{1}{2(t-2x)}-\frac{x}{2(t-2x)^2}+O(x^2)\right),
\end{split}
\end{equation*}
which shows that $f$ extends to a $C^1$ function on $[0,1]$, and
so
\begin{equation*}
\int_0^1 \vert f_t^\prime (x)\vert \, dx \ll 1,
\end{equation*}
uniformly for $t$ in compacts of $(2,\infty)$.
\par
The equality $V(Qx)=Q^{-1} f(x)$, $x\in (0,1]$, implies now that
both $\| V\|_\infty$ and the total variation of $V$ on $(0,Q]$ are
$\ll Q^{-1}$. Thus we may apply Lemma \ref{LEMMA1} to \eqref{5.2}
and conclude, also using $c+c_1<1$, that
\begin{equation}\label{5.3}
\begin{split}
P^+_{I,Q} & =\frac{c_I}{\zeta(2)} \int_0^Q V(q)\, dq+O_\delta (
E_{c,c_1,\delta}(Q)) =\frac{c_I}{\zeta(2)}
\int_0^1 f_t(x)\, dx +O_\delta ( E_{c,c_1,\delta}(Q)) \\
& =\frac{c_I}{2\zeta(2)} \int_0^1 \psi (x,t)\, dx+O_\delta (
E_{c,c_1,\delta}(Q)) .
\end{split}
\end{equation}
\par
One can see in a similar way that the contribution of integrals on
the intervals $[\beta_{k-1},\beta_k]$ in \eqref{3.7} for $k\geq 1$
and $t>2$ is
\begin{equation}\label{5.4}
P^-_{I,Q}(t)=\frac{c_I}{2\zeta(2)} \int_0^1 \psi(x,t)\, dx+O (
E_{c,c_1,\delta}(Q)) .
\end{equation}
Proposition \ref{PROP4} now follows from \eqref{5.3} and
\eqref{5.4}.

\section{The case $1<t< 2$}
In this section we prove

\begin{proposition}\label{PROP5}
Suppose $I$ is a subinterval of $[0,1]$ of size $\vert I\vert
\asymp Q^{-c}$ for some $0<c<1$. Then for any $c_1$ with $c+c_1<1$
and $\delta >0$
\begin{equation*}
P_{I,Q}(t) =\frac{c_I}{\zeta(2)}\left(
 \int_0^{t-1} \hspace{-6pt}\psi(x,t)\, dx +
 \int_{t-1}^1 \phi(x,t)\, dx\right)
+O_\delta( E_{c,c_1,\delta}(Q)) \quad (Q \rightarrow \infty),
\end{equation*}
with $\phi$ and $\psi$ as in Theorem {\em \ref{T1.1}}. The
estimate holds uniformly in $t$ on compacts of $(1,2)$.
\end{proposition}
\par
In this case \eqref{3.7} gives
\begin{equation*}
\begin{split}
P_{I,Q}(t) = & \sum\limits_{\gamma \in \FI}
\int_{\alpha_0}^{\beta_0} w_{C_0}(\omega) \max \{ q+q^\prime
-tQ,0\} \, d\omega \\ & +\hspace{-8pt}\sum\limits_{\gamma \in \FI}
\sum\limits_{k=1}^\infty \int_{\alpha_k}^{\alpha_{k-1}} \Big(
w_{C_k}(\omega)\max \{ q_{k+1}-tQ,0\} +w_{B_k}(\omega) \max
\{ q_k-tQ,0\} \Big) d\omega \\
 & +\hspace{-8pt}\sum\limits_{\gamma \in \FI}
\sum\limits_{k=1}^\infty \int_{\beta_{k-1}}^{\beta_k}
 \Big( w_{C_{-k}}(\omega) \max \{
q^\prime_{k+1}-tQ,0\} +w_{A_{-k}}(\omega)\max \{
q^\prime_k-tQ,0\}\Big)\, d\omega .
\end{split}
\end{equation*}
\par
We break the main term above according as to whether $q+q^\prime
>tQ$ or $q+q^\prime \leq tQ$. Thus we first estimate
\begin{equation*}
\begin{split}
P^{>}_{I,Q}(t) := & \sum\limits_{\substack{\gamma \in \FI \\
q+q^\prime >tQ}}\int_{\alpha_0}^{\beta_0} \left(
\frac{1}{Q}-w_{A_+}(\omega) -w_{B_-}(\omega) \right) (q+q^\prime
-tQ)\, d\omega \\ & +\hspace{-5pt}
\sum\limits_{\substack{\gamma \in \FI \\
q+q^\prime >tQ}}\ \sum\limits_{k=1}^\infty
\int_{\alpha_k}^{\alpha_{k-1}} \Big(
w_{C_k}(\omega)q_{k+1}+w_{B_k}(\omega)q_k -tQ\big(
w_{C_k}(\omega)+w_{B_k} (\omega)\big) \Big) d\omega \\
& +\hspace{-5pt}\sum\limits_{\substack{\gamma \in \FI \\
q+q^\prime >tQ}}\ \sum\limits_{k=1}^\infty \,
\int_{\beta_{k-1}}^{\beta_k} \Big(
w_{C_{-k}}(\omega)q^\prime_{k+1} +w_{A_{-k}}(\omega) q^\prime_k-tQ
\big( w_{C_{-k}}(\omega)+w_{A_{-k}}(\omega)\big) \Big) d\omega .
\end{split}
\end{equation*}
\par
Using \eqref{4.2} and \eqref{4.3} we may also write
\begin{equation*}
\begin{split}
P^{>}_{I,Q}(t) & =\sum\limits_{\substack{\gamma \in \FI \\
q+q^\prime >tQ}} \int_{\alpha_0}^{\beta_0} \left(
\frac{1}{Q}-w_{A_+}(\omega)-w_{B_-}(\omega)\right) (q+q^\prime
-tQ)\, d\omega \\ & +\sum\limits_{\substack{\gamma \in \FI \\
q+q^\prime >tQ}} \int_{\alpha_\infty}^{\alpha_0} \left(
1-w_{A_+}(\omega)
q-tQ\Big( \frac{1}{Q} -w_{A_+}(\omega)\Big) \right)d\omega \\
 & + \sum\limits_{\substack{\gamma \in \FI \\
q+q^\prime >tQ}} \int_{\beta_0}^{\beta_\infty} \left(
1-w_{B_-}(\omega)q^\prime -tQ \Big(
\frac{1}{Q}-w_{B_-}(\omega)\Big) \right) d\omega \\
& =\widetilde{A}_0+\widetilde{A}_1+\widetilde{A}_2+\widetilde{A}_3
,
\end{split}
\end{equation*}
with
\begin{equation*}
\begin{split}
\widetilde{A}_0 & =(1-t) \sum\limits_{\substack{\gamma \in \FI \\
q+q^\prime >tQ}} (\beta_\infty-\alpha_\infty) ,
\quad \widetilde{A}_1 =\sum\limits_{\substack{\gamma \in \FI \\
q+q^\prime >tQ}} (tQ-q)\int_{\alpha_\infty}^{\beta_0}
w_{A_+}(\omega) \,
d\omega ,\\ \widetilde{A}_2 & =\sum\limits_{\substack{\gamma \in \FI \\
q+q^\prime >tQ}} (tQ-q)\int_{\alpha_0}^{\beta_\infty}
w_{B_-}(\omega) \, d\omega ,\\
\widetilde{A}_3 & =\sum\limits_{\substack{\gamma \in \FI \\
q+q^\prime
>tQ}}\ \int_{\alpha_0}^{\beta_0} \left(
\frac{q+q^\prime}{Q}-1-q^\prime
w_{A_+}(\omega)-qw_{B_-}(\omega)\right) d\omega .
\end{split}
\end{equation*}

\begin{figure}[ht]
\begin{center}
\unitlength 0.5mm
\begin{picture}(100,100)(0,0)
\texture{ccc 0}
\shade\path(66.6666,100)(100,66.6666)(100,0)(66.6666,33.3333)(66.6666,100)

\texture{c 0}
\shade\path(0,100)(66.6666,33.3333)(66.6666,100)(0,100)
\path(0,100)(100,100)(100,0)(0,100)
\path(22.2222,100)(33.333,66.6666)
\path(33.333,100)(66.666,33.3333) \path(66.666,100)(100,66.6666)
\dottedline{2}(33.3333,100)(33.3333,66.6666)

\put(0,100){\makebox(0,0){{\tiny $\bullet$}}}
\put(22.2222,100){\makebox(0,0){{\tiny $\bullet$}}}
\put(33.3333,100){\makebox(0,0){{\tiny $\bullet$}}}
\put(66.6666,100){\makebox(0,0){{\tiny $\bullet$}}}
\put(100,100){\makebox(0,0){{\tiny $\bullet$}}}
\put(33.3333,66.6666){\makebox(0,0){{\tiny $\bullet$}}}
\put(66.6666,33.3333){\makebox(0,0){{\tiny $\bullet$}}}
\put(66.6666,100){\makebox(0,0){{\tiny $\bullet$}}}
\put(100,0){\makebox(0,0){{\tiny $\bullet$}}}
\put(100,66.6666){\makebox(0,0){{\tiny $\bullet$}}}

\put(-9,100){\makebox(0,0){{\small $(0,1)$}}}
\put(109,100){\makebox(0,0){{\small $(1,1)$}}}
\put(109,0){\makebox(0,0){{\small $(1,0)$}}}
\put(113,66.6666){\makebox(0,0){{\small $(1,t-1)$}}}
\put(66,106){\makebox(0,0){{\small $(t-1,1)$}}}
\put(53,29){\makebox(0,0){{\small $(t-1,2-t)$}}}
\put(22,61){\makebox(0,0){{\small $\big(
\frac{t-1}{2},\frac{3-t}{2}\big)$}}}
\put(40,106){\makebox(0,0){{\small $\big( \frac{t-1}{2},1\big)$}}}
\put(16,106){\makebox(0,0){{\small $\big( \frac{t-1}{3},1\big)$}}}
\end{picture}
\end{center}
\caption{The set $\cup_{k=1}^\infty \Omega_k \cap \TT$ when
$1<t<2$} \label{Figure7}
\end{figure}
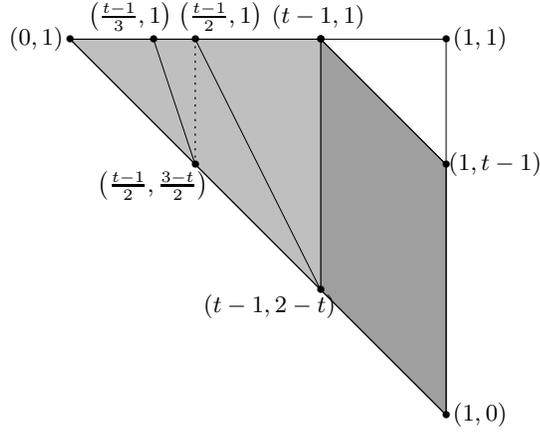

We proceed to estimate $\widetilde{A}_0$, $\widetilde{A}_1$,
$\widetilde{A}_2$ and $\widetilde{A}_3$ by noticing that
\eqref{4.5} yields
\begin{equation}\label{6.1}
\widetilde{A}_3=0.
\end{equation}
\par
Next $\widetilde{A}_1$ is estimated in a similar way as $A_1$ was
in Section 4, only with the difference that the summation over
$\gamma \in \FI$ is being done under the additional requirement
$q+q^\prime >tQ$. This is not going to produce any change in the
error, and will only affect the main terms. As in \eqref{4.10} and
\eqref{4.11} we obtain
\begin{equation*}
\widetilde{A}_0 =(1-t)\sum\limits_{\substack{\gamma \in \FI \\
q+q^\prime >tQ}} \frac{1}{qq^\prime (1+\gamma^2)} +O( Q^{-1}) .
\end{equation*}
Then, as in the proof of Lemma \ref{L4.4}, we find that
\begin{equation*}
\widetilde{A}_0 =(1-t)\hspace{-1pt} \sum\limits_{q=1}^\infty \
\sum\limits_{\substack{q^\prime \in \II:=(tQ-q,Q]\\ a\in \JJ:=qI
\\ -aq^\prime =1\hspace{-6pt}\pmod{q}}} \hspace{-15pt}
f_q(q^\prime,a) =c_I(1-t)\hspace{-13pt} \sum\limits_{(t-1)Q<q\leq
Q}\hspace{-6pt} \frac{\varphi(q)}{q} V(q)+O_\delta (
E_{c,c_1,\delta}(Q)),
\end{equation*}
where this time we take
\begin{equation*}
V(q)=\frac{1}{q} \ln \frac{Q}{tQ-q} ,\quad q\in ( (t-1)Q,Q].
\end{equation*}
But $V(Qx)=Q^{-1}\widetilde{V}(x)$ and the function
\begin{equation*}
\widetilde{V}(x)=\frac{1}{x}\ln \frac{1}{t-x} \ ,\qquad x\in
[t-1,1],
\end{equation*}
is $C^1$ on $[t-1,1]$. Hence both the $L^\infty$-norm and the
total variation of $V$ on the interval $[ (t-1)Q,Q]$ are $\ll
Q^{-1}$, uniformly in $t$ on compacts of $(1,2)$. Lemma
\ref{LEMMA1} applies now and yields
\begin{equation}\label{6.2}
\begin{split}
\widetilde{A}_0 & =\frac{c_I(1-t)}{\zeta(2)} \int_{(t-1)Q}^Q
V(q)\, dq+O_\delta ( E_{c,c_1,\delta}(Q)) \\ & =
\frac{c_I(1-t)}{\zeta(2)} \int_{t-1}^1 \frac{1}{x}  \ln
\frac{1}{t-x}\ dx+O_\delta ( E_{c,c_1,\delta}(Q)) .
\end{split}
\end{equation}
\par
Proceeding as in Section 4 (see \eqref{4.18}--\eqref{4.20}) we
find
\begin{equation*}
\begin{split}
\widetilde{A}_1 & =\sum\limits_{(t-1)Q<q\leq Q} \frac{tQ-q}{2qQ^2}
\hspace{-8pt}
\sum\limits_{\substack{tQ-q<q^\prime \leq Q \\
a\in qI \\ -aq^\prime =1\hspace{-6pt} \pmod{q}}} \hspace{-8pt}
\frac{1}{1+a^2/q^2} +O(Q^{-1})  \\
& =\frac{c_I}{2Q^2} \sum\limits_{(t-1)Q<q\leq Q}
\frac{\varphi(q)}{q^2} \, (tQ-q)\big( q-(t-1)Q\big) +
O_\delta ( E_{c,c_1,\delta}(Q)) \\
& =\frac{c_I}{2Q^2 \zeta(2)} \int_{(t-1)Q}^Q \frac{(tQ-q)\big(
q-(t-1)Q\big)}{q}\ dq+O_\delta ( E_{c,c_1,\delta}(Q)).
\end{split}
\end{equation*}
This immediately gives
\begin{equation}\label{6.3}
\widetilde{A}_1 =\frac{c_I}{2\zeta(2)} \int_{t-1}^1
\frac{(t-x)(x-t+1)}{x}\ dx+O_\delta ( E_{c,c_1,\delta}(Q)).
\end{equation}
In a similar way we find
\begin{equation}\label{6.4}
\widetilde{A}_2 =\frac{c_I}{2\zeta(2)} \int_{t-1}^1
\frac{(t-x)(x-t+1)}{x}\ dx+O_\delta ( E_{c,c_1,\delta}(Q)).
\end{equation}
From \eqref{6.1}--\eqref{6.4} we now collect
\begin{equation}\label{6.5}
\begin{split}
P^{>}_{I,Q}(t) & =\frac{c_I(1-t)}{\zeta(2)} \int_{t-1}^1
\frac{1}{x}\ln \frac{1}{t-x}\ dx +\frac{c_I}{\zeta(2)}
\int_{t-1}^1 \frac{(t-x)(x-t+1)}{x}\ dx \\ & \qquad \qquad
+O_\delta ( E_{c,c_1,\delta}(Q)).
\end{split}
\end{equation}
\par
It remains to estimate the contribution of Farey fractions of
order $Q$ with $q+q^\prime \leq tQ$ to $P_{I,Q}(t)$, which is
\begin{equation*}
P^{<}_{I,Q}(t):=B_1+B_2,
\end{equation*}
where $B_1$ denotes
\begin{equation*}
\sum\limits_{\substack{\gamma \in \FI \\ q+q^\prime \leq tQ}}\
\sum\limits_{k=1}^\infty \int_{\alpha_k}^{\alpha_{k-1}} \Big(
w_{C_k}(\omega)\max \{ q_{k+1}-tQ,0\} +w_{B_k} (\omega) \max \{
q_k-tQ,0\} \Big)\, d\omega ,
\end{equation*}
and $B_2$ denotes
\begin{equation*}
\sum\limits_{\substack{\gamma \in \FI \\ q+q^\prime \leq tQ}}\
\sum\limits_{k=1}^\infty \ \int_{\beta_{k-1}}^{\beta_k} \Big(
w_{C_{-k}}(\omega) \max \{ q^\prime_{k+1}-tQ,0\}
+w_{A_{-k}}(\omega) \max \{ q^\prime_k-tQ,0 \} \Big)\, d\omega .
\end{equation*}
In this case one also has
\begin{equation*}
K=\left[ \frac{tQ-q^\prime}{q}\right] \geq 1,
\end{equation*}
and as in Section 5 we find
\begin{equation*}
\begin{split}
B_1 + & O_\delta ( E_{c,c_1,\delta}(Q))
\\ & =\sum\limits_{k=1}^\infty\hspace{-5pt}
\sum\limits_{\substack{\gamma \in \FI\\ (q,q^\prime) \in
Q(\Omega_k \cap \TT)}}\hspace{-10pt}
\int_{\alpha_k}^{\alpha_{k-1}}  \Big( w_{C_k}(\omega)
\max \{ q_{k+1}-tQ,0\} +w_{B_k} \max \{ q_k-tQ,0\} \Big)\, d\omega \\
& =\sum\limits_{k=2}^\infty \sum\limits_{q\in QI_k} \big(
S_k(q)+T_k(q)\big) +\sum\limits_{(t-1)Q<q\leq Q}
\ \sum\limits_{\substack{Q-q<q^\prime <tQ-q \\
a\in qI \\ -aq^\prime =1\hspace{-6pt} \pmod{q}}}
f_1(q,q^\prime,a) \\
& =\frac{c_I}{2\zeta(2)} \int_0^{t-1} \hspace{-8pt} \psi(x,t)\,
dx+\frac{c_I}{2\zeta(2)} \int_{t-1}^1 \hspace{-4pt}
\frac{(1-x)^2}{x}
\int_{1-x}^{t-x} \frac{2(x+y)-t}{xy(x+y)} \ dy\, dx \\
& =\frac{c_I}{2\zeta(2)} \int_0^{t-1} \hspace{-8pt} \psi(x,t)\,
dx+ \frac{c_I}{2\zeta(2)} \int_{t-1}^1 \hspace{-4pt}
\frac{(1-x)^2}{x} \int_{1-x}^{t-x} \left(
\frac{2}{y}-\frac{t}{x}\Big( \frac{1}{y}-\frac{1}{y+x}\Big)
\right)\, dy\, dx,
\end{split}
\end{equation*}
and thus
\begin{equation}\label{6.6}
\begin{split}
B_1 & =\frac{c_I}{2\zeta(2)} \int_{t-1}^1 \frac{(1-x)^2}{x} \left(
2\ln \frac{t-x}{1-x}-\frac{t}{x} \ln
\frac{t-x}{t(1-x)} \right) dx \\
& \qquad +\frac{c_I}{2\zeta(2)} \int_0^{t-1} \psi(x,t)\, dx
+O_\delta ( E_{c,c_1,\delta}(Q)) .
\end{split}
\end{equation}
In a similar way we find that $B_2$ can too be expressed as in
\eqref{6.6}, and thus
\begin{equation}\label{6.7}
\begin{split}
P^{<}_{I,Q}(t) & =\frac{c_I}{\zeta(2)} \int_{t-1}^1
\frac{(1-x)^2}{x} \left( 2\ln \frac{t-x}{1-x}-\frac{t}{x} \ln
\frac{t-x}{t(1-x)} \right) dx \\
& \qquad +\frac{c_I}{\zeta(2)} \int_0^{t-1} \psi(x,t)\, dx
+O_\delta ( E_{c,c_1,\delta}(Q)) .
\end{split}
\end{equation}
Proposition \ref{PROP5} follows now from \eqref{6.5} and
\eqref{6.7}.

\section{Proof of Theorem \ref{T1.1}}
We may assume without loss of generality that $\omega \in
\big[0,\frac{\pi}{4}\big]$, thus estimate for small $\eps >0$ the
quantity
\begin{equation*}
\P_\eps (t)=\frac{4}{\pi} \left| \bigg\{ (x,\omega) \in Y_\eps
\times \Big[ 0,\frac{\pi}{4}\Big) \, ;\, \tau_\eps
(x,\omega)>\frac{t}{2\eps} \bigg\} \right| .
\end{equation*}
We partition the interval $[0,1]$ as a union of $N$ intervals
$I_j=[\tan \omega_j,\tan \omega_{j+1}]$ of equal size, with $N=[
\eps^{-c}]$, thus with $\vert I_j \vert =\frac{1}{N} \asymp
\eps^c$, where $0<c<1$ is to be chosen later. For each $j$ we set
\begin{equation*}
Q_j^-=\left[ \frac{\cos \omega_{j+1}}{2\eps+2\eps^{c+1}} \right] ,
\quad Q_j^+=\left[ \frac{\cos
\omega_j}{2\eps-2\eps^{c+1}}\right]+1.
\end{equation*}
Since $\omega_j\in [ 0,\frac{\pi}{4}]$, we have $Q_j^\pm \asymp
\eps^{-1}$, and thus $\vert I_j \vert \asymp \eps^c$. Moreover,
for $\omega \in [\omega_j,\omega_{j+1}]$ we have
\begin{equation}\label{7.1}
\frac{1}{2Q_j^+} < \frac{\eps-\eps^{c+1}}{\cos \omega_j} \leq
\frac{\eps}{\cos \omega_j} \leq \frac{\eps}{\cos \omega} \leq
\frac{\eps}{\cos \omega_{j+1}} \leq \frac{\eps+\eps^{c+1}}{\cos
\omega_{j+1}} \leq \frac{1}{2Q_{j}^-}  .
\end{equation}
From the definition of $Q_j^\pm$ and from
\begin{equation}\label{7.2}
\vert \cos y-\cos x\vert \leq \vert \sin (x-y)\vert \leq \vert
\tan x-\tan y\vert , \qquad x,y\in [ 0,\pi/4 ],
\end{equation}
we infer that
\begin{equation*}
Q_j^+-Q_j^- \ll \frac{\cos \omega_j}{2\eps -2\eps^{c+1}}
-\frac{\cos \omega_{j+1}}{2\eps +2\eps^{c+1}} \ll
\frac{\eps^{c+1}}{\eps}+\frac{\cos \omega_j-\cos
\omega_{j+1}}{2\eps} \ll \eps^{c-1}
\end{equation*}
and
\begin{equation}\label{7.3}
Q_j^\pm =\frac{\cos \omega_j}{2\eps}+O(\eps^{c-1}).
\end{equation}
\par
\newtheorem*{Note2}{Remark 2}
\begin{Note2}
{\em If $\omega \in \big[ 0,\frac{\pi}{4}\big]$ and $\lambda_\pm$
are such that $\lambda_-<\frac{\cos \omega}{2\eps}<\lambda_+$,
then for all $x\in Y_\eps$ we have
\begin{equation*}
\widetilde{\tau}_{1/(2\lambda_+)} (x,\omega)+\eps >\tau_\eps
(x,\omega) >\widetilde{\tau}_{1/(2\lambda_-)} (x,\omega)-\eps.
\end{equation*}
This shows in turn that if for each interval $I=[\tan
\omega_0,\tan \omega_1] \subseteq [0,1]$ we denote
\begin{equation*}
\P_{\eps,I}:=\left| \bigg\{ (x,\omega)\, ;\, x\in Y_\eps,\ \tan
\omega \in I,\ \tau_\eps (x,\omega)>\frac{t}{2\eps} \bigg\}
\right|,
\end{equation*}
then for any integers $Q_\pm$ such that $Q_-<\frac{\cos
\omega_1}{2\eps}<\frac{\cos \omega_0}{2\eps}<Q_+$ we have
\begin{equation*}
\widetilde{\P}_{I,Q_-} \left( \frac{t+\eps}{2\eps}\right)-\pi
\eps^2 \leq \P_{\eps,I}(t) \leq \widetilde{\P}_{I,Q_+}\left(
\frac{t-\eps}{2\eps}\right) .
\end{equation*}}
\end{Note2}
\par
By the previous remark we infer
\begin{equation}\label{7.4}
\widetilde{\P}_{I_j,Q_j^-} \left( \frac{t+\eps}{2\eps}\right) -\pi
\eps^2 \leq \P_{\eps,I_j}(t) \leq \widetilde{\P}_{I_j,Q_j^+}\left(
\frac{t-\eps}{2\eps}\right) ,\qquad j=1,\ldots,N.
\end{equation}
\par
For small $\eps >0$ we also have
\begin{equation*}
\frac{t}{2\eps+2\eps^{c+1}}<\frac{t-\eps}{2\eps}<\frac{t+\eps}{2\eps}
<\frac{t}{2\eps-2\eps^{c+1}}  ,
\end{equation*}
uniformly in $t$ on compacts of $(0,\infty)$. Thus \eqref{7.4},
\eqref{7.3}, and Lemma \ref{LEMMA5} yield
\begin{equation}\label{7.5}
\P_{\eps,I_j} (t)\leq \widetilde{\P}_{I_j,Q_j^+} \left(
\frac{t}{2\eps +2\eps^{c+1}} \right)=P_{I_j,Q_j^+} \left(
\frac{t}{1+\eps^c}\right) +O(\eps^{2c})
\end{equation}
and
\begin{equation}\label{7.6}
\P_{\eps,I_j} (t)\geq \widetilde{\P}_{I_j,Q_j^-} \left(
\frac{t}{2\eps-2\eps^{c+1}}\right)-\pi \eps^2 =P_{I_j,Q_j^-}
\left( \frac{t}{1-\eps^c} \right) +O(\eps^{2c}).
\end{equation}
\par
By the definition of $\P$ we see that for any compact interval
$K\subset (0,\infty)\setminus \{ 1,2\}$, there exists a constant
$C_K>0$ such that
\begin{equation}\label{7.7}
\vert \P (t_1)-\P(t_2)\vert \leq C_K \vert t_1-t_2 \vert, \quad
t_1,t_2\in K.
\end{equation}
Now by Propositions \ref{PROP3}, \ref{PROP4}, \ref{PROP5} we know
that for any $j\in \{1,2,\cdots,N\}$ we have for small $\eps >0$
\begin{equation*}
\begin{split}
P_{I_j,Q_j^\pm} \left( \frac{t}{1\pm \eps^c}\right) &
=P_{I_j,Q_j^\pm} \big( t(1+O(\eps^c))\big) \\ & =c_{I_j}(
\P(t)+O(\eps^c))+O_\delta (\eps^{2c}+\eps^{c+c_1}
+\eps^{1/2-2c_1-\delta}) ,
\end{split}
\end{equation*}
uniformly in $t$ on compacts of $(0,\infty)\setminus \{ 1,2\}$.
Here $\P(t)$ is defined as in Theorem 1.1. Summing over $j$ the
inequalities \eqref{7.5} and \eqref{7.6}, and using also
\begin{equation}\label{7.8}
\sum\limits_{j=1}^N c_{I_j} =\int_0^1 \frac{du}{1+u^2}
=\frac{\pi}{4},
\end{equation}
$N\leq \eps^{-c}$, and \eqref{7.7}, we gather
\begin{equation*}
\sum\limits_{j=1}^N \P_{\eps,I_j} (t)=\frac{\pi}{4}\, \P
(t)+O_\delta (\eps^c+\eps^{c_1}+\eps^{1/2-2c_1-c-\delta}).
\end{equation*}
\par
For obvious symmetry reasons we can only consider $\omega \in [
0,\frac{\pi}{4}]$. Thus, after normalizing the Lebesgue measure
$\mu_\eps$ on $Y_\eps$ by dividing by $\frac{\pi}{4}\ar
(Y_\eps)=\frac{\pi (1-\pi \eps^2)}{4}$, we get
\begin{equation*}
\P_\eps (t)=\P(t) +O_\delta
(\eps^c+\eps^{c_1}+\eps^{1/2-2c_1-c-\delta}).
\end{equation*}
The proof of Theorem \ref{T1.1} is completed by taking
$c=c_1=\frac{1}{8}$.

\section{The geometric free path length in the case $0<t\leq 1$}
In this and the next two sections we shall take $\omega \in [
0,\frac{\pi}{4}]$, and analyze the geometric free path length in
the case of vertical scatterers of height $2\delta$ centered at
integer lattice points. In this setup we will consider the phase
space $(\widetilde{\Sigma}_{\delta,I}, \frac{d\mu}{2\delta})$,
where $\delta >0$, $I=[\tan \omega_0,\tan \omega_1]\subseteq
[0,1]$ is an interval,
$\widetilde{\Sigma}_{\delta,I}=[-\delta,\delta] \times
[\omega_0,\omega_1]$ and $d\mu$ is the (non-normalized) Lebesgue
measure on $\widetilde{\Sigma}_{\delta,I}$. The trajectory will
therefore start at a point $(0,y)$, $y\in [-\delta,\delta]$, under
angle $\omega$, with $\tan \omega \in I$. Recall that the free
path length is denoted by $\widetilde{\tau}_\delta$ in this case.
Given $\lambda >0$, consider
\begin{equation}\label{8.1}
\begin{split}
\widetilde{G}_{\delta,I}(\lambda) & =\frac{1}{2\delta} \left| \{
(y,\omega)\in \widetilde{\Sigma}_{\delta,I} \, ;\,
\widetilde{\tau}_\delta (y,\omega)>\lambda \}\right| =
\frac{1}{2\delta} \int_{\omega_0}^{\omega_1} \int_{-\delta}^\delta
e_\lambda ( \widetilde{\tau}_\delta (y,\omega))\, dy\, d\omega .
\end{split}
\end{equation}
\par
Actually it will suffice to take $\delta=\frac{1}{2Q}$ for
properly chosen integers $Q$. The first goal will be to estimate
the distribution of the free path length
$\widetilde{\tau}_{1/(2Q)} (x,\omega)$ when we average over
$(x,\omega)\in \widetilde{\Sigma}_{1/(2Q),I}$, under the
assumptions that $I=[\tan \omega_0,\tan \omega_1] \subseteq [0,1]
$ is a short interval of length $\vert I\vert \asymp \eps^{1/8}$
for small $\eps$, and that $Q$ is a (large) integer such that
$Q=\frac{\cos \omega_0}{2\eps}+O(\eps^{1/8-1})$. Concretely, we
will be interested in the quantity
\begin{equation*}
\widetilde{\G}_{I,Q}(t):=\widetilde{G}_{1/(2Q),I}(t)=Q\left| \{
(x,\omega)\in \widetilde{\Sigma}_{1/(2Q),I} \, ;\,
\widetilde{\tau}_{1/(2Q)} (x,\omega)>t\} \right| .
\end{equation*}
\par
In the remainder of the paper we take $c=c_1=\frac{1}{8}$. We set
\begin{equation*}
\Delta_\lambda (x)=e_{(-\infty,\lambda)}(x) =
\begin{cases} 1 & \mbox{\rm if $\ x<\lambda$;} \\
0 & \mbox{\rm if $\ x\geq \lambda$.}
\end{cases}
\end{equation*}
\par
A direct application of Propositions \ref{PROP1} and \ref{PROP2},
with widths $w$ given by \eqref{3.8} and $\alpha_k$, $\beta_k$ by
\eqref{4.1}, yields the following formula, derived from
\eqref{3.5} by replacing $\max \{ q-x,0\}$ with $Q\Delta_q (x)$,
and valid for any $t,\eps_\ast >0$:
\begin{equation}\label{8.2}
\begin{split}
& \widetilde{\G}_{I,Q}\left( \frac{t}{2\eps_\ast}\right)
=Q\sum\limits_{\gamma \in \FI} \sum\limits_{k=1}^\infty
\int_{\alpha_k}^{\alpha_{k-1}} w_{C_k}(\omega) \Delta_{q_{k+1}}
\left( \frac{t\cos \omega}{2\eps_\ast} \right)\, d\omega \\ &
\qquad +\mbox{\rm eight other terms where $\frac{t\cos
\omega}{2\eps_\ast}$ appears in an analogous way.}
\end{split}
\end{equation}
\par
This quantity will be compared with the one obtained by
substituting $tQ$ in place of $\frac{t\cos \omega}{2\eps_\ast}$ in
\eqref{8.2}, as in Section 3. For this purpose we shall consider
\begin{equation}\label{8.3}
\begin{split}
G_{I,Q}(t):= Q \sum\limits_{\gamma \in \FI}
\sum\limits_{k=1}^\infty & A^-_{Q,\gamma,k}(t,\omega)\, d\omega
+Q\sum\limits_{\gamma \in \FI} \int_{\alpha_0}^{\beta_0}
A^{(0)}_{Q,\gamma,k} (t,\omega) \, d\omega \\ &
+Q\sum\limits_{\gamma \in \FI} \sum\limits_{k=1}^\infty \
\int_{\beta_{k-1}}^{\beta_k} A^+_{Q,\gamma,k}(t,\omega) \,
d\omega,
\end{split}
\end{equation}
with
\begin{equation*}
\begin{split}
A^-_{Q,\gamma,k} (t,\omega) & =w_{C_k}(\omega) \Delta_{q_{k+1}}
(tQ) +w_{B_k}(\omega) \Delta_{q_k}(tQ) +w_{A_+} (\omega)\Delta_q
(tQ),\\
A^{(0)}_{Q,\gamma,k}(t,\omega) &
=w_{C_0}(\omega)\Delta_{q+q^\prime} (tQ) +w_{B_0}
(\omega)\Delta_{q^\prime} (tQ) +w_{A_0}(\omega)\Delta_q (tQ),\\
A^+_{Q,\gamma,k}(t,\omega) & =w_{C_{-k}}(\omega)
\Delta_{q^\prime_{k+1}} (tQ) +w_{A_{-k}} (\omega)
\Delta_{q_k^\prime} (tQ) +w_{B_-}(\omega) \Delta_{q^\prime} (tQ).
\end{split}
\end{equation*}

\newtheorem*{Note3}{Remark 3}
\begin{Note3}
{\em If $I=[\tan \omega_0,\tan \omega_1] \subseteq [0,1]$ and
$0<\lambda_- \leq \frac{\cos \omega_1}{2\eps} <\frac{\cos
\omega_0}{2\eps}\leq \lambda_+$, then owing to \eqref{8.2},
\eqref{8.3} and to the fact that $x\mapsto \Delta_\lambda (x)$ is
monotonically decreasing we have
\begin{equation*}
G_{I,Q} \left( \frac{t\lambda_+}{Q}\right)\leq
\widetilde{\G}_{I,Q} \left( \frac{t}{2\eps}\right) \leq
G_{I,Q}\left( \frac{t\lambda_-}{Q}\right) .
\end{equation*}}
\end{Note3}
\par
The argument, based on inequality \eqref{3.9} used to compare
$\widetilde{\P}_{I,Q} ( \frac{t}{2\eps})$ with $P_{I,Q}(t)$ in
Lemma \ref{LEMMA5}, is not going to apply here because
$\Delta_\lambda$ is not a Lipschitz function. Nevertheless, we can
overcome this problem by appealing again to a soft monotonicity
argument, based on Remark 3 and on the fact (which can be seen
directly from the definition of the function $\G(t)$) that for any
compact $K\subset (0,\infty)\setminus \{ 1,2\}$, there exists a
constant $C_K >0$ such that
\begin{equation}\label{8.4}
\vert \G (t_1)-\G (t_2)\vert \leq C_K \vert t_1-t_2 \vert, \quad
t_1,t_2 \in K.
\end{equation}
\par
In this and the the next two sections we will analyze the
asymptotic of the quantity $G_{I,Q}(t)$ for large integers $Q$ and
short intervals $I$ such that $\vert I\vert \asymp Q^{-1/8}$. We
note at this point that the relation \eqref{1.3} is hinted by
formula \eqref{8.3} and by
\begin{equation*}
\frac{d}{dt} \max \{ q-tQ,0\}=-Q \Delta_q (tQ),\qquad t\neq
\frac{q}{Q} .
\end{equation*}
\par
For the sake of space, the error estimates which are similar to
the ones already derived in the first part of the paper are going
to be more sketchy.

\begin{proposition}\label{P8.2}
For every interval $I\subseteq [0,1]$ of size $\vert I\vert \asymp
Q^{-1/8}$ and every $\delta >0$
\begin{equation*}
G_{I,Q}(t)=\left(1 -\frac{t}{\zeta(2)}\right) c_I +O_\delta
(Q^{-1/4+\delta}) \qquad (Q\rightarrow \infty).
\end{equation*}
The estimate holds uniformly in $t\in (0,1]$.
\end{proposition}
\par
\begin{proof} Since $0<t\leq 1$, we have $\min \{ q_k,q_k^\prime\}\geq
q+q^\prime >tQ$ for all $k\geq 1$. Thus we infer from \eqref{8.3},
as in formula \eqref{4.4}, that
\begin{equation*}
G_{I,Q} (t)= G^{(1)}_{I,Q} (t) +G^{(2)}_{I,Q}(t)+G^{(3)}_{I,Q}(t),
\end{equation*}
with
\begin{equation*}
\begin{split}
& G^{(1)}_{I,Q}(t) :=Q\hspace{-6pt} \sum\limits_{\gamma \in \FI}
\sum\limits_{k=1}^\infty \int_{\alpha_k}^{\alpha_{k-1}}
\hspace{-3pt}\big( w_{C_k}(\omega)+w_{B_k}(\omega)\big) \, d\omega
+ Q\hspace{-6pt} \sum\limits_{\gamma \in \FI}
\int_{\alpha_0}^{\beta_0}
w_{C_0}(\omega)\, d\omega \\
& \qquad \qquad \qquad +Q\hspace{-6pt} \sum\limits_{\gamma \in
\FI} \sum\limits_{k=1}^\infty \int_{\beta_{k-1}}^{\beta_k}
\hspace{-2pt} \big( w_{C_{-k}}(\omega)+w_{A_{-k}}(\omega) \big)\, d\omega \\
& \qquad =Q\hspace{-6pt} \sum\limits_{\gamma \in \FI}
\int_{\alpha_\infty}^{\alpha_0} \left(
\frac{1}{Q}-w_{A_+}(\omega)\right) d\omega
+Q\hspace{-6pt}\sum\limits_{\gamma \in \FI}
\int_{\beta_0}^{\beta_\infty}
\left( \frac{1}{Q}-w_{B_-}(\omega) \right) d\omega \\
& \qquad \qquad \qquad +Q\hspace{-6pt}\sum\limits_{\gamma \in \FI}
\int_{\alpha_0}^{\beta_0} \left( \frac{1}{Q}
-w_{A_+}(\omega)-w_{B_-}(\omega) \right) d\omega \\
& \qquad =\sum\limits_{\gamma \in \FI} \hspace{-5pt}(\beta_\infty
-\alpha_\infty) -Q\hspace{-6pt}\sum\limits_{\gamma \in \FI}
\int_{\alpha_\infty}^{\beta_0} \hspace{-5pt} w_{A_+} (\omega) \,
d\omega -Q\hspace{-6pt} \sum\limits_{\gamma \in \FI}
\int_{\alpha_0}^{\beta_\infty} \hspace{-5pt} w_{B_-}(\omega)\,
d\omega ,
\end{split}
\end{equation*}

\begin{equation*}
\begin{split}
G^{(2)}_{I,Q}(t) & :=Q \hspace{-6pt}\sum\limits_{\substack{\gamma \in \FI \\
q>tQ}} \sum\limits_{k=1}^\infty \int_{\alpha_k}^{\alpha_{k-1}}
\hspace{-3pt} w_{A_+}(\omega)\, d\omega + Q\hspace{-6pt}
\sum\limits_{\substack{\gamma \in \FI
\\ q>tQ}} \int\limits_{\alpha_0}^{\beta_0} w_{A_+}(\omega)\,
d\omega =Q\hspace{-6pt} \sum\limits_{\substack{\gamma \in \FI \\
q>tQ}} \int_{\alpha_\infty}^{\beta_0} w_{A_+} (\omega) ,
\\
 G^{(3)}_{I,Q}(t) & :=Q \hspace{-6pt}
\sum\limits_{\substack{\gamma \in \FI \\ q^\prime >tQ}}
\sum\limits_{k=1}^\infty \int_{\beta_{k-1}}^{\beta_k}
w_{B_-}(\omega)\, d\omega +Q\hspace{-6pt}
\sum\limits_{\substack{\gamma \in \FI \\ q^\prime
>tQ}} \int_{\alpha_0}^{\beta_0} w_{B_-}(\omega)\, d\omega
\\ &  =Q\hspace{-6pt} \sum\limits_{\substack{\gamma \in \FI \\
q^\prime >tQ}} \int_{\alpha_0}^{\beta_\infty} w_{B_-}(\omega)\,
d\omega .
\end{split}
\end{equation*}
\par
From \eqref{4.12} we gather
\begin{equation}\label{8.5}
\sum\limits_{\gamma \in \FI} \hspace{-4pt} (\beta_\infty
-\alpha_\infty) = \hspace{-8pt} \sum\limits_{\gamma \in \FI}
 \left( \arctan \frac{a^\prime}{q^\prime} -\arctan
\frac{a}{q}\right) =c_I+O_\delta (Q^{-1/4+\delta}).
\end{equation}
On the other hand, \eqref{4.17} gives
\begin{equation}\label{8.6}
Q\int_{\alpha_\infty}^{\beta_0} w_{A_+}(\omega)\, d\omega
=\frac{1}{2qQ(1+\gamma^2)}+O\left( \frac{1}{q^2Q^2}\right).
\end{equation}
We can show in a similar way that
\begin{equation}\label{8.7}
Q\int_{\alpha_0}^{\beta_\infty} w_{B_-}(\omega)\, d\omega
=\frac{1}{2q^\prime Q(1+\gamma^{\prime 2})} +O\left(
\frac{1}{q^{\prime 2} Q^2}\right) .
\end{equation}
From the formulas for $G_{I,Q}^{(1)}$, $G_{I,Q}^{(2)}$,
$G_{I,Q}^{(3)}$ and from \eqref{8.5}--\eqref{8.7} we infer
\begin{equation}\label{8.8}
\begin{split}
G_{I,Q}(t)  & =c_I -\sum\limits_{\substack{ \gamma \in \FI \\
q\leq tQ}} \left( \frac{1}{2qQ(1+\gamma^2)} +O\Big( \frac{1}{q^2
Q^2}\Big) \right) \\ & \qquad \qquad
-\sum\limits_{\substack{ \gamma \in \FI \\
q^\prime \leq tQ}} \left( \frac{1}{2q^\prime Q(1+\gamma^{\prime
2})} +O\Big( \frac{1}{q^{\prime 2} Q^2}\Big) \right)+O_\delta
(Q^{-1/4+\delta})\\ & =c_I-\sum\limits_{\substack{\gamma \in \FI
\\ q\leq tQ}} \frac{1}{2qQ(1+\gamma^2)} -
\sum\limits_{\substack{\gamma \in \FI \\ q^\prime \leq tQ}}
\frac{1}{2q^\prime Q(1+\gamma^{\prime 2})} +O_\delta
(Q^{-1/4+\delta}).
\end{split}
\end{equation}
\par
Finally we show as at the end of Section 4 that
\begin{equation*}
\begin{split}
\sum\limits_{\substack{\gamma \in \FI \\ q\leq tQ}}
\frac{1}{2qQ(1+\gamma^2)} & =\frac{1}{2Q} \sum\limits_{1\leq q\leq
tQ}\frac{1}{q} \sum\limits_{\substack{Q-q<q^\prime \leq Q \\ a\in qI \\
-aq^\prime =1\hspace{-6pt}\pmod{q}}} \frac{1}{1+a^2/q^2} \\
& =\frac{1}{2Q} \sum\limits_{1\leq q\leq tQ} \frac{1}{q}\cdot
\frac{\varphi(q)}{q^2}  q^2c_I +O_\delta
(Q^{-1/4+\delta}) \\
& =\frac{c_I}{2Q} \sum\limits_{\substack{1\leq q\leq tQ}}
\frac{\varphi(q)}{q}+O_\delta (Q^{-1/4+\delta})\\ & =\frac{c_I
t}{2\zeta(2)} +O_\delta (Q^{-1/4+\delta}).
\end{split}
\end{equation*}
A similar formula holds for the second sum in \eqref{8.8}, and
therefore we get
\begin{equation*}
G_{I,Q}(t)=c_I-\frac{c_I t}{\zeta(2)}+O_\delta (Q^{-1/4+\delta}).
\end{equation*}
It is clear that these estimate hold uniformly in $t\in [0,1]$.
\end{proof}

\section{The geometric free path in the case $t>2$}
In this section we prove in the setting of Section 8 the following
result

\begin{proposition}\label{P9.1}
For every interval $I\subseteq [0,1]$ of size $\vert I\vert \asymp
Q^{-\frac{1}{8}}$ and every $\delta >0$
\begin{equation*}
G_{I,Q}(t)=\frac{c_I}{\zeta(2)} \int_0^1 \frac{(1-x)^2}{x^2} \ln
\frac{(t-x)^2}{t(t-2x)} \, dx +O_\delta (Q^{-1/4+\delta}) \qquad
(Q\rightarrow \infty).
\end{equation*}
The estimate holds uniformly in $t$ on compacts of $(2,\infty)$.
\end{proposition}
\par
\begin{proof} We proceed as in Section 5, estimating first
\begin{equation}\label{9.1}
\begin{split}
& \widetilde{S}_2 (\gamma,t): =Q\sum\limits_{k=K+1}^\infty
\int_{\alpha_k}^{\alpha_{k-1}} \big(
w_{C_k}(\omega)+w_{B_k}(\omega)\big) \, d\omega \\ & \qquad =
Q\int_{\alpha_\infty}^{\alpha_K} \left(
\frac{1}{Q}-w_{A_+}(\omega)\right) d\omega  =Q \int_{\arctan
\frac{a}{q}}^{\arctan \frac{a_K-1/Q}{q_K}} (q\tan \omega -a)\,
d\omega \\ & \qquad =\frac{Q( 1-q/Q)^2 q}{2q^2
q_K^2(1+\gamma^2)}+O\left( \frac{Q}{q^2 q_K^3}\right) =\frac{Q(
1-q/Q)^2}{2qq_K^2 (1+\gamma^2)}+O\left( \frac{Q}{q^2
q_K^3}\right),
\end{split}
\end{equation}
and then
\begin{equation}\label{9.2}
\begin{split}
\widetilde{S}_1 (\gamma,t) & =Q\int_{\alpha_K}^{\alpha_{K-1}}
w_{C_K}(\omega)\, d\omega =Q\int_{\arctan
\frac{a_K-1/Q}{q_K}}^{\arctan \frac{a_{K-1}-1/Q}{q_{K-1}}}
 \left( a_{K-1}-q_{K-1}\tan \omega
-\frac{1}{Q}\right) \, d\omega \\ & =\frac{Q( 1-q/Q)^2
q_{K-1}}{2(1+\gamma^2) q_{K-1}^2 q_K^2} +O\left(
\frac{Qq_{K-1}}{q_{K-1}^3 q_K^3}\right)  =\frac{Q(
1-q/Q)^2}{2(1+\gamma^2)q_{K-1} q_K^2} +O\left( \frac{1}{q^2
q^2_1}\right) .
\end{split}
\end{equation}
\par
In this case we find from \eqref{8.3}, \eqref{9.1} and
\eqref{9.2}, as in Section 5, that
\begin{equation*}
\begin{split}
G_{I,Q}(t) & =2\sum\limits_{\gamma \in \FI} \frac{Q( 1-q/Q)^2
(q_{K-1}+q)}{2(1+\gamma^2) qq_{K-1} q_K^2} +O( Q^{-1} \ln Q)  \\
& =\sum\limits_{\gamma \in \FI}
\frac{Q(1-q/Q)^2}{(1+\gamma^2)qq_{K-1} q_K} +O(Q^{-1}\ln Q)\\ &
=\sum\limits_{k=2}^\infty \sum\limits_{q\in QI_k} \big(
\widetilde{S}_k (q)+\widetilde{T}_k (q) \big) +O( Q^{-1}\ln Q) ,
\end{split}
\end{equation*}
with
\begin{equation*}
\begin{split}
\widetilde{S}_k(q) & = \sum\limits_{\substack{q^\prime \in
QJ_{k,q}^{(1)},\,  a\in qI \\ -aq^\prime =1 \hspace{-6pt}
\pmod{q}}} \hspace{-10pt} \widetilde{f}_k (q,q^\prime,a) ,\quad
\widetilde{T}_k(q) =\hspace{-10pt} \sum\limits_{\substack{q^\prime
\in QJ_{k,q}^{(0)} ,\, a\in qI \\ -aq^\prime =1 \hspace{-6pt}
\pmod{q}}}\hspace{-10pt} \widetilde{f}_{k-1} (q,q^\prime,q-a),\\
 \widetilde{f}_k (q,q^\prime,a) & = \frac{Q( 1-q/Q)^2}{(
1+a^2/q^2) qq_{k-1} q_k} ,\quad q^\prime \in \II=QJ_{k,q}^{(1)} ,
\ a \in \JJ=qI.
\end{split}
\end{equation*}
\par
Employing the same technique as in Section 5 and the fact that one
gets a similar result while integrating between $\beta_{k-1}$ and
$\beta_k$, we find that $G_{I,Q}(t)$ can be expressed, up to an
error term of order $O_\delta (Q^{-1/4+\delta})$, as
\begin{equation*}
\begin{split}
\sum\limits_{k=2}^\infty \sum\limits_{q\in QI_k} & \frac{\varphi
(q)}{q^2} \left( \hspace{5pt} \iint_{QJ_{k,q}^{(1)}\times qI}
\hspace{-8pt} \widetilde{f}_k(q,q^\prime,a) \, dq^\prime \, da+
\iint_{QJ_{k,q}^{(0)}\times qI} \widetilde{f}_{k-1}
(q,q^\prime,a)\, dq^\prime\, da   \right)
\\ &
=c_I \sum\limits_{k=2}^\infty \sum\limits_{q\in QI_k}
\frac{\varphi (q)}{q} \cdot \frac{Q( 1-q/Q)^2}{q} \left(
\int_{Q-q}^{tQ-kq} \frac{dq^\prime}{q_{k-1}q_k} +\int_{tQ-kq}^Q
\frac{dq^\prime}{q_{k-2}q_{k-1}}  \right)
\\ &
=c_I \sum\limits_{k=2}^\infty \sum\limits_{q\in QI_k}
\frac{\varphi (q)}{q} \cdot \frac{Q( 1-q/Q)^2}{q} \left(
\int_{Q+(k-1)q}^{tQ} \frac{dy}{y(y-q)}
+\int_{tQ-q}^{Q+(k-1)q}\hspace{-6pt} \frac{dy}{y(y-q)} \right).
\\
& =c_I \sum\limits_{q=1}^Q \frac{\varphi(q)}{q} \cdot \frac{Q(
1-q/Q)^2}{q^2} \ln \frac{(tQ-q)^2}{tQ(tQ-2q)}.
\end{split}
\end{equation*}
This is further equal to
\begin{equation*}
\begin{split}
\frac{c_I}{\zeta(2)} \int_0^Q  \frac{Q( 1-q/Q)^2}{q^2} \ln
\frac{(tQ-q)^2}{tQ(tQ-2q)} \, dq & = \frac{c_I}{\zeta(2)} \int_0^Q
\frac{Q( 1-\frac{q}{Q})^2}{q^2} \ln \frac{(tQ-q)^2}{tQ(tQ-2q)} \,
dq \\ &  =\frac{c_I}{\zeta(2)} \int_0^1 \frac{(1-x)^2}{x^2} \ln
\frac{(t-x)^2}{t(t-2x)} \, dx ,
\end{split}
\end{equation*}
which is the desired conclusion.
\end{proof}

\newpage

\section{The geometric free path in the case $1<t<2$}
In this section we prove in the setting of Section 8 the following
result

\begin{proposition}\label{P10.1}
For every interval $I\subseteq [0,1]$ of size $\vert I\vert \asymp
Q^{-\frac{1}{8}}$ and $\delta >0$
\begin{equation*}
\begin{split}
G_{I,Q}(t) & =\frac{c_I}{\zeta(2)} \left( \int_{t-1}^1
\frac{1}{x}\ln\frac{1}{t-x}\, dx  -2+t+
(t-1)\ln\frac{1}{t-1}+\int_{t-1}^1
\frac{(1-x)^2}{x^2}\ln\frac{t-x}{t(1-x)}\right.
 \\ &  \left. \qquad \qquad \qquad+\int_0^{t-1}
\frac{(1-x)^2}{x^2}\ln\frac{(t-x)^2}{t(t-2x)}\, dx\right)+O_\delta
(Q^{-1/4+\delta}) \qquad \qquad (Q\rightarrow \infty).
\end{split}
\end{equation*}
The estimate holds uniformly in $t$ on compacts of $(1,2)$.
\end{proposition}
\par
\begin{proof} Since $1<t<2$, we have $\max \{ q,q^\prime \}\leq tQ$ and we
infer from \eqref{8.2} that
\begin{equation*}
G_{I,Q}(t) = G^>_{I,Q}(t)+ G^<_{I,Q}(t),
\end{equation*}
where $G^>_{I,Q}(t)$, respectively $G^<_{I,Q}(t)$, contains the
contribution of Farey fractions in $\FI$ with $q+q^\prime >tQ$,
respectively with $q+q^\prime \leq tQ$.
\par
When $q+q^\prime >tQ$ we have $\min \{ q_k,q_k^\prime \}>tQ$,
$k\geq 1$, and therefore
\begin{equation*}
\begin{split}
G^>_{I,Q}(t) := & Q\hspace{-6pt}\sum\limits_{\substack{ \gamma \in \FI \\
q+q^\prime >tQ}} \sum\limits_{k=1}^\infty
\int_{\alpha_k}^{\alpha_{k-1}} \big( w_{C_k}(\omega) +
w_{B_k}(\omega) \big) \, d\omega \\ &
+ Q\hspace{-6pt}\sum\limits_{\substack{ \gamma \in \FI \\
q+q^\prime >tQ}} \int_{\alpha_0}^{\beta_0}  w_{C_0}(\omega)\,
d\omega +Q\hspace{-6pt}\sum\limits_{\substack{ \gamma \in \FI \\
q+q^\prime >tQ}} \sum\limits_{k=1}^\infty
\int_{\beta_{k-1}}^{\beta_k} \big( w_{C_{-k}}(\omega)
+w_{A_{-k}}(\omega) \big) \, d\omega \\ &
=Q\hspace{-6pt}\sum\limits_{\substack{\gamma \in \FI \\ q+q^\prime
>tQ}}\ \int_{\alpha_\infty}^{\alpha_0} \left( \frac{1}{Q}-w_{A_+}
(\omega) \right) d\omega \\ &
+Q\hspace{-6pt}\sum\limits_{\substack{ \gamma \in \FI \\
q+q^\prime >tQ}}  \int_{\alpha_0}^{\beta_0} \left(
\frac{1}{Q}-w_{A_+}(\omega)-w_{B_-}(\omega) \right)
d\omega +Q\hspace{-6pt}\sum\limits_{\substack{ \gamma \in \FI \\
q+q^\prime >tQ}} \left( \frac{1}{Q}-w_{B_-}
(\omega) \right) d\omega \\
& =\sum\limits_{\substack{\gamma \in \FI \\ q+q^\prime
>tQ}}\hspace{-5pt} (\beta_\infty -\alpha_\infty)
-Q\hspace{-6pt}\sum\limits_{\substack{\gamma \in \FI \\ q+q^\prime
>tQ}} \int_{\alpha_\infty}^{\beta_0} w_{A_+}(\omega) \, d\omega
-Q\hspace{-6pt}\sum\limits_{\substack{\gamma \in \FI
\\ q+q^\prime >tQ}} \int_{\alpha_0}^{\beta_\infty}
w_{B_-} (\omega)\, d\omega .
\end{split}
\end{equation*}
\par
Standard considerations as in Sections 6 and 8 show that,
uniformly in $t$ on compacts of $(1,2)$ and up to an error term of
order $O(Q^{-1}\ln Q)$, $G_{I,Q}^>(t)$ can be expressed as

\begin{equation*}
\begin{split}
\sum\limits_{\substack{\gamma \in \FI \\
q+q^\prime >tQ}} & \frac{1}{qq^\prime (1+\gamma^2)}
-\sum\limits_{\substack{\gamma \in \FI \\ q+q^\prime >tQ}} \left(
\frac{1}{2qQ(1+\gamma^2)} +\frac{1}{2q^\prime Q(1+\gamma^{\prime
2})}\right)  \\
& =\sum\limits_{\substack{\gamma \in \FI \\
q+q^\prime >tQ}} \left( \frac{1}{qq^\prime (1+\gamma^2)} -
\frac{1}{qQ(1+\gamma^2)}\right) =\frac{1}{Q}
\sum\limits_{\substack{\gamma \in \FI \\ q+q^\prime
>tQ}} \frac{Q-q^\prime}{qq^\prime (1+\gamma^2)}  \\
& =\frac{1}{Q} \sum\limits_{(t-1)Q<q\leq Q} \frac{1}{q}
\sum\limits_{\substack{tQ-q<q^\prime \leq Q\\ a\in qI \\
-aq^\prime =1\hspace{-6pt}\pmod{q}}} \frac{Q-q^\prime}{q^\prime
( 1+a^2/q^2)} \\
 & =\frac{1}{Q} \sum\limits_{(t-1)Q<q\leq Q} \frac{1}{q}\cdot
\frac{\varphi(q)}{q^2} \, qc_I \int_{tQ-q}^Q
\frac{Q-q^\prime}{q^\prime}\ dq^\prime
+O_\delta (Q^{-1/4+\delta}) \\
& =\frac{c_I}{Q} \sum\limits_{(t-1)Q<q\leq Q} \frac{\varphi(q)}{q}
\cdot \frac{1}{q} \int_{tQ-q}^Q \frac{Q-q^\prime}{q^\prime}\
dq^\prime +O_\delta (Q^{-1/4+\delta})
\\ & =\frac{c_I}{Q\zeta(2)} \int_{(t-1)Q}^Q \frac{1}{q}
\int_{tQ-q}^Q \frac{Q-q^\prime}{q^\prime}\ dq^\prime \ dq+O_\delta
(Q^{-1/4+\delta}) \\ & =\frac{c_I}{\zeta(2)} \int_{t-1}^1
\frac{1}{x} \ln \frac{1}{t-x}\ dx-\frac{c_I}{\zeta(2)} \left(
2-t+(1-t)\ln \frac{1}{t-1}\right)+O_\delta (Q^{-1/4+\delta}).
\end{split}
\end{equation*}
\par
Next we estimate
\begin{equation*}
\begin{split}
G^<_{I,Q}(t) & :=Q \sum\limits_{\substack{\gamma \in \FI
\\ q+q^\prime \leq tQ}}\sum\limits_{k=1}^\infty
\int_{\alpha_k}^{\alpha_{k-1}}  \big( w_{C_k}
(\omega)+w_{B_k}(\omega)\big)\, d\omega  \\ & \qquad +Q
\sum\limits_{\substack{\gamma \in \FI
\\ q+q^\prime \leq tQ}} \sum\limits_{k=1}^\infty\
\int_{\beta_{k-1}}^{\beta_k}  \big(
w_{C_{-k}}(\omega)+w_{A_{-k}}(\omega)\big) d\omega ,
\end{split}
\end{equation*}
and find as in Sections 5, 6 and 9 that $G^<_{I,Q}(t)$ can be
expressed, up to an error term of order $O_\delta
(Q^{-1/4+\delta})$, as
\begin{equation*}
\begin{split}
& 2\sum\limits_{k=2}^\infty \sum\limits_{q\in QI_k} \big(
\widetilde{S}_k (q)+ \widetilde{T}_k (q)\big) +2\sum\limits_{q\in
QI_1} \sum\limits_{\substack{q^\prime \in QJ_{1,q}^{(1)},\,a\in qI \\
-aq^\prime =1\hspace{-6pt} \pmod{q}}}
\widetilde{f}_1 (q,q^\prime,a)  \\
&  \qquad =\frac{c_I}{\zeta(2)} \int_0^{t-1} \frac{(1-x)^2}{x^2}
\ln \frac{(t-x)^2}{t(t-2x)}\, dx + \frac{c_I}{\zeta(2)}
\int_{(t-1)Q}^Q \int_{Q-q}^{tQ-q} \frac{Q( 1-q/Q)^2}{qq^\prime
(q+q^\prime)}\ dq^\prime \, dq  \\
& \qquad = \frac{c_I}{\zeta(2)} \left( \int_0^{t-1}
\frac{(1-x)^2}{x^2} \, \ln \frac{(t-x)^2}{t(t-2x)}\, dx
 +\int_{t-1}^1 \frac{(1-x)^2}{x^2} \,
 \ln \frac{t-x}{t(1-x)} \right) .
\end{split}
\end{equation*}
We conclude the proof by adding the formulas for $G^<_{I,Q}(t)$
and $G^>_{I,Q}(t)$.
\end{proof}

\section{Proof of Theorem \ref{T1.2}}
Identifying $\Sigma_\eps^+$ with
\begin{equation*}
\left\{ (\eps \mathrm{e}^{\mathrm{i} \alpha},\omega)\, ;\,
-\omega-\pi/2\leq \alpha\leq \omega+\pi/2 \right\} =\left\{ (\eps
\mathrm{e}^{\mathrm{i}(\omega+\beta)},\omega)\, ;\, \beta \in [
-\pi/2,\pi/2] \right\} ,
\end{equation*}
the (non-normalized) Liouville measure on the phase space
$\Sigma_\eps^+$ is expressed as
\begin{equation*}
d\lambda_\eps =\eps \big< (\cos \alpha,\sin \alpha),(\cos
\omega,\sin \omega)\big> \ d\alpha\, d\omega = \eps \cos (\omega
-\alpha)\, d\alpha\, d\omega = \eps \cos \beta \, d\beta\, d\omega
.
\end{equation*}

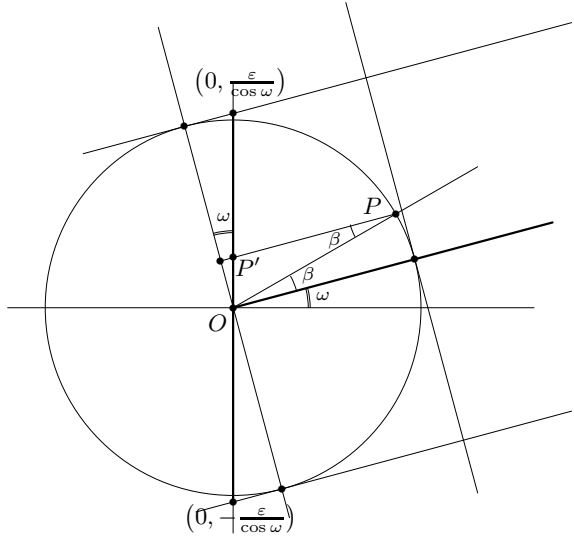
\begin{figure}[ht]
\begin{center}
\unitlength 0.5mm
\begin{picture}(0,130)(20,-60)
\put(0,0){\arc{100}{-6.28}{0}} \path(-60,0)(80,0)
\path(0,-60)(0,60) \path(0,0)(65,37.52777)
\path(-3.3495,12.5005)(43.30127,25)
\path(15,-55.9807)(-20,74.64102) \path(90,-27.64838)(-10,-54.4433)
\path(90,75.8792)(-40,41.04584) \path(65,-49,39814)(30,81.22364)

\put(0,51.7638){\makebox(0,0){{\tiny $\bullet$}}}
\put(0,-51.7638){\makebox(0,0){{\tiny $\bullet$}}}
\put(43.30127,25){\makebox(0,0){{\tiny $\bullet$}}}
\put(48.2963,12.94095){\makebox(0,0){{\tiny $\bullet$}}}
\put(12.94095226,-48.2963){\makebox(0,0){{\tiny $\bullet$}}}
\put(-12.94095226,48.2963){\makebox(0,0){{\tiny $\bullet$}}}
\put(0,13.3975){\makebox(0,0){{\tiny $\bullet$}}}
\put(-3.34936,12.5005){\makebox(0,0){{\tiny $\bullet$}}}
\put(0,0){\makebox(0,0){{\tiny $\bullet$}}}

\put(37,27){\makebox(0,0){{\small $P$}}}
\put(4,11){\makebox(0,0){{\small $P'$}}}
\put(-4,-4){\makebox(0,0){{\small $O$}}}
\put(-2.5,24){\makebox(0,0){{\scriptsize $\omega$}}}
\put(27,18){\makebox(0,0){{\scriptsize $\beta$}}}
\put(20,8.5){\makebox(0,0){{\scriptsize $\beta$}}}
\put(24,3){\makebox(0,0){{\scriptsize $\omega$}}}
\put(2,-56){\makebox(0,0){{\small $\big( 0,-\frac{\eps}{\cos
\omega}\big)$}}} \put(2,60){\makebox(0,0){{\small $\big(
0,\frac{\eps}{\cos \omega}\big)$}}}

\put(0,0){\arc{40}{-1.83}{-1.57}}
\put(0,0){\arc{41}{-1.83}{-1.57}}
\put(43.30127,25){\arc{25}{-3.65}{-3.38}}
\put(0,0){\arc{35}{-.5}{-.25}} \put(0,0){\arc{40}{-0.25}{0}}
\put(0,0){\arc{41}{-0.25}{0}}

\thicklines \path(0,-51.7638)(0,51.7638) \path(0,0)(85,22.7757)
\end{picture}
\end{center}
\caption{The parametrization of $\Sigma_\eps^+$} \label{Figure8}
\end{figure}

Next we shall consider a fixed interval $I=[\tan \omega_0,\tan
\omega_1]\subseteq [0,1]$, define
\begin{equation*}
\Sigma_{\eps,I}^+:=\left\{ (\eps
\mathrm{e}^{\mathrm{i}(\omega+\beta)},\omega)\, ;\, \vert
\beta\vert \leq \pi/2,\ \omega_0\leq \omega \leq \omega_1 \right\}
,
\end{equation*}
and estimate
\begin{equation*}
\G_{\eps,I}(t):=\frac{1}{\eps}\, \lambda_\eps \left( \bigg\{
(x,\omega)\in \Sigma_{\eps,I}^+ \, ;\, \tau_\eps
(x,\omega)>\frac{t}{2\eps} \bigg\}\right) .
\end{equation*}
To each point $P=\eps \mathrm{e}^{\mathrm{i}(\omega+\beta)}$ we
associate (see Figure \ref{Figure8}) the point $P^\prime (0,y)$,
where $y=\frac{\eps \sin \beta}{\cos \omega}\in \big[
-\frac{\eps}{\cos \omega},\frac{\eps}{\cos \omega}\big]$. Note
that
\begin{equation*}
\lambda_\eps (\Sigma_{\eps,I}^+)=\eps \int_{\omega_0}^{\omega_1}
\int_{-\pi/2}^{\pi/2} \cos \beta\, d\beta \, d\omega =2\eps c_I .
\end{equation*}
\par
Since $PP^\prime$ has slope $\tan \omega$, we have the obvious
inequality
\begin{equation*}
\left| \tau_\eps (\eps
\mathrm{e}^{\mathrm{i}(\omega+\beta)},\omega)-\widetilde{\tau}_{\eps/\cos
\omega} \bigg( \frac{\eps \sin \beta}{\cos \omega},\omega \bigg)
\right| <2\eps,
\end{equation*}
and as a consequence we can write
\begin{equation*}
\begin{split}
\G_{\eps,I}(t) & =\int_{\omega_0}^{\omega_1} \int_{-\pi/2}^{\pi/2}
e_{t/(2\eps)} \big( \tau_\eps (\eps
\mathrm{e}^{\mathrm{i}(\omega+\beta)},\omega)\big) \cos \beta \,
d\beta\, d\omega \\ & \leq \int_{\omega_0}^{\omega_1}
\int_{-\pi/2}^{\pi/2} e_{t/(2\eps)-2\eps} \left(
\widetilde{\tau}_{\eps/\cos \omega} \Big( \frac{\eps \sin \beta
}{\cos \omega} ,\omega \Big) \right) \cos \beta\, d\beta\,
d\omega \\
& =\int_{\omega_0}^{\omega_1} \int_{-\eps/\cos \omega}^{\eps/\cos
\omega} \frac{\cos \omega}{\eps} \, e_{(t-4\eps^2)/(2\eps)} \big(
\widetilde{\tau}_{\eps/\cos \omega} (y,\eps )\big)\, dy\, d\omega
.
\end{split}
\end{equation*}
When $0<\lambda_- \leq \frac{\cos \omega_1}{2\eps} <\frac{\cos
\omega_0}{2\eps}\leq \lambda_+$, obvious monotonicity properties
yield
\begin{equation}\label{11.1}
\begin{split}
\G_{\eps,I}(t) & \leq 2\lambda_+ \int_{\omega_0}^{\omega_1}
\int_{-1/(2\lambda_-)}^{1/(2\lambda_-)} e_{(t-4\eps^2)/(2\eps)}
\big( \widetilde{\tau}_{1/(2\lambda_+)}(y,\omega)\big)\, dy\,
d\omega
\\ & =2\lambda_+ \int_{\omega_0}^{\omega_1}
\int_{-1/(2\lambda_+)}^{1/(2\lambda_+)} e_{(t-4\eps^2)/(2\eps)}
\big( \widetilde{\tau}_{1/(2\lambda_+)} (y,\omega)\big) \, dy\,
d\omega +O\left( \lambda_+ \Big(
\frac{1}{\lambda_-}-\frac{1}{\lambda_+} \Big) \right) \\
& =2\widetilde{G}_{1/(2\lambda_+),I} \left(
\frac{t-4\eps^2}{2\eps}\right) +O\left(
\frac{\lambda_+}{\lambda_-}-1\right) ,
\end{split}
\end{equation}
with $\widetilde{G}_{\delta,I}$ as defined in \eqref{8.1}. Using
similar arguments we infer
\begin{equation}\label{11.2}
\G_{\eps,I}(t) \geq 2\widetilde{G}_{1/(2\lambda_-),I} \left(
\frac{t+4\eps^2}{2\eps}\right) +O\left(
1-\frac{\lambda_-}{\lambda_+} \right).
\end{equation}
\par
Take now $\eps >0$ small, and suppose that $\vert I\vert \asymp
\eps^{1/8}$ and $Q^\pm$ are two integers such that
\begin{equation*}
Q^- \leq\frac{\cos \omega_1}{2\eps} \leq\frac{\cos
\omega_0}{2\eps} \leq Q^+, \quad Q^\pm =\frac{\cos
\omega_0}{2\eps}+O(\eps^{1/8-1}), \quad
\frac{Q^\pm}{Q^\mp}=1+O(\eps^{1/8}) .
\end{equation*}
Such integers can be chosen for instance as at the beginning of
Section 7 with $c=\frac{1}{8}$. Fix also a compact $K\subset
(0,\infty)\setminus \{ 1,2\}$. Applying successively \eqref{11.1},
Remark 3, Propositions \ref{P8.2}, \ref{P9.1}, \ref{P10.1}, and
inequality \eqref{8.4}, we infer that
\begin{equation}\label{11.3}
\begin{split}
\G_{\eps,I}(t) & \leq 2\widetilde{G}_{1/(2Q_+),I} \left(
\frac{t-4\eps^2}{2\eps}\right)+O\left( \frac{Q^+}{Q^-}-1\right)
\\ & = 2\widetilde{\G}_{I,Q^+} \left(
\frac{t-4\eps^2}{2\eps}\right)+O\left( \frac{Q^+}{Q^-}-1\right) \\
& \leq 2G_{I,Q^+}\left( \frac{(t-4\eps^2)Q^-}{Q^+} \right)
+O\left( \frac{Q^+}{Q^-}-1\right) \\ & =2G_{I,Q^+} \Big(
(t-4\eps^2)\big(
1+O(\eps^{1/8})\big) \Big) +O(\eps^{1/8}) \\
& =2G_{I,Q^+} \big( t+O(\eps^{1/8})\big) +O(\eps^{1/8})\\
& =2c_I \G (t)+O_\delta (\eps^{1/8-\delta}) \qquad \mbox{\rm
uniformly in $t\in K$}.
\end{split}
\end{equation}
\par
In a similar way we infer from \eqref{11.2} and the previous
arguments that
\begin{equation}\label{11.4}
\G_{\eps,I}(t) \geq 2c_I \G (t)+O_\delta (\eps^{1/8-\delta})
\qquad \mbox{\rm uniformly in $t\in K$.}
\end{equation}
\par
Consider now a partition of $[0,1]$ with intervals
$\{I_j\}_{j=1}^N$, where $N =[\eps^{-1/8}]$ and $\vert I_j\vert
=\frac{1}{N}\asymp \eps^{1/8}$. Summing over $j$ we find as a
result of \eqref{11.3}, \eqref{11.4} and \eqref{7.8} that
\begin{equation*}
\G_{\eps,[0,1]} (t)=\sum\limits_{j=1}^N \G_{\eps,I_j}(t)
=\frac{\pi}{2}\, \G(t)+ O_\delta (\eps^{1/8-\delta}),
\end{equation*}
and thus
\begin{equation*}
\begin{split}
\frac{\G_{\eps,[0,1]}(t)}{\lambda_\eps ( \Sigma^+_{\eps,[0,1]})} &
= \frac{\lambda_\eps ( \{ (x,\omega)\in \Sigma^+_{\eps,[0,1]} \,
;\, 2\eps \tau_\eps
(x,\omega)>t \})}{\lambda_\eps ( \Sigma_{\eps,[0,1]}^+)} \\
& =\frac{\pi}{2}\cdot \frac{\eps \G(t)}{2\eps c_{[0,1]}}+O_\delta
(\eps^{1/8-\delta}) =\G (t)+O_\delta (\eps^{1/8-\delta}).
\end{split}
\end{equation*}
For obvious symmetry reasons we can only consider $\omega \in
\big[ 0,\frac{\pi}{4}\big]$, therefore
\begin{equation*}
\G_\eps (t)=\G(t) +O_\delta (\eps^{1/8-\delta}),
\end{equation*}
which ends the proof of Theorem \ref{T1.2}.

\section{Estimates of $C_\eps =\ln \langle \tau_\eps \rangle -
\langle \ln \tau_\eps \rangle$} In this section we prove Theorem
\ref{T1.3} (i). Part (ii) then follows from (i) and from relation
(2.8) in \cite{Ch2}.
\par
We consider the probability measures $\nu_0$ and
$\widetilde{\nu}_\eps$ on $[0,\infty)$ defined by
\begin{equation*}
\begin{split}
& \int_0^\infty f(u)\, d\nu_0(u)=\int_0^\infty f(u) g(u)\, du ,\\
& \int_0^\infty f(u)\, d\widetilde{\nu}_\eps (u) =
\int_{\Sigma_\eps^+} f(2\eps \tau_\eps)\, d\nu_\eps ,\quad f\in
C_c ([0,\infty)) .
\end{split}
\end{equation*}
\par
As a result of Theorem \ref{T1.2}
\begin{equation*}
\lim_{\eps \rightarrow 0^+} \int_t^\infty d\widetilde{\nu}_\eps
(u)=\int_t^\infty d\nu_0 (u),\quad t>0,
\end{equation*}
which implies
\begin{equation*}
\lim\limits_{\eps \rightarrow 0^+} \int_0^\infty f(u) \,
d\widetilde{\nu}_\eps (u)=\int_0^\infty f(u)\, d\nu_0(u) ,\quad
f\in C_c ( [0,\infty )),
\end{equation*}
meaning that $\widetilde{\nu}_\eps \rightarrow \nu_0$ vaguely as
$\varepsilon \rightarrow 0^+$. Since
\begin{equation*}
\lim\limits_{\eps \rightarrow 0^+} \frac{1}{x} \int_x^\infty
d\widetilde{\nu}_\eps (u)=\frac{1}{x} \int_x^\infty
d\nu_0(u),\quad x\geq 1,
\end{equation*}
and the map
\begin{equation*}
x\mapsto \frac{1}{x} \int_x^\infty d\nu_0 (u)=\frac{1}{x}
\int_x^\infty g(u)\, du
\end{equation*}
belongs to $L^1 ( [1,\infty),dx)$ because $g(u)=O(u^{-3})$, $u\geq
1$, the Lebesgue Dominated Convergence theorem yields
\begin{equation*}
\lim\limits_{\eps \rightarrow 0^+} \int_1^\infty \frac{1}{x}
\int_x^\infty d\widetilde{\nu}_\eps (u) \ dx =\int_1^\infty
\frac{1}{x} \int_x^\infty d\nu_0(u)\ dx <\infty .
\end{equation*}
\par
Using Fubini's theorem, these double integrals can also be
expressed as
\begin{equation*}
\begin{split}
\int_1^\infty \int_1^\infty \frac{1}{x}  e_{[1,u]} (x)\,
d\widetilde{\nu}_\eps (u)\ dx & = \int_1^\infty \int_1^\infty
\frac{1}{x} e_{[1,u]} (x) \, dx\ d\widetilde{\nu}_\eps (u) \\ &
=\int_1^\infty \int_1^u \frac{dx}{x} \ d\widetilde{\nu}_\eps (u)
=\int_1^\infty \ln u\, d\widetilde{\nu}_\eps (u),
\end{split}
\end{equation*}
and respectively as
\begin{equation*}
\int_1^\infty \int_1^\infty \frac{1}{x} e_{[1,u]}(x)\, d\nu_0(u)\
dx =\int_1^\infty \ln u\, d\nu_0 (u).
\end{equation*}
It follows that for any (small) $\eps >0$
\begin{equation}\label{12.1}
\int_1^\infty \ln u\, d\widetilde{\nu}_\eps (u)<\infty ,
\end{equation}
and also that
\begin{equation}\label{12.2}
\lim\limits_{\eps \rightarrow 0^+} \int_1^\infty \ln u\,
d\widetilde{\nu}_\eps (u)\, du=\int_1^\infty g(u)\ln u\, du.
\end{equation}
\par
We show in a similar way that
\begin{equation*}
\lim\limits_{\eps \rightarrow 0^+} \int_0^1 \ln u\,
d\widetilde{\nu}_\eps (u) =\int_0^1 g(u)\ln u\, du =
\frac{6}{\pi^2} \int_0^1 \ln u\, du=-\frac{6}{\pi^2}
\end{equation*}
by using Fubini's theorem which gives in turn
\begin{equation*}
\begin{split}
\int_0^1 \ln u\, d\widetilde{\nu}_\eps (u) & = -\int_0^1 \int_u^1
\frac{1}{x}\ dx\ d\widetilde{\nu}_\eps (u)=-\int_0^1 \int_0^1
\frac{1}{x} e_{[u,1]} (x)\, dx\ d\widetilde{\nu}_\eps (u)
\\ & =-\int_0^1 \int_0^1 \frac{1}{x}
e_{[u,1]} (x)\, d\widetilde{\nu}_\eps (u)\ dx =-\int_0^1
\frac{1}{x} \int_0^x d\widetilde{\nu}_\eps (u) \ dx.
\end{split}
\end{equation*}
\par
By \eqref{12.1} and \eqref{12.2} we get
\begin{equation*}
\begin{split}
-C & :=\int_0^\infty g(u)\ln u\, du = \lim\limits_{\eps
\rightarrow 0^+} \int_0^\infty \ln u\, d\widetilde{\nu}_\eps (u) =
\lim\limits_{\eps \rightarrow 0^+} \int_{\Sigma_\eps^+}
\ln (2\eps \tau_\eps)\, d\nu_\eps \\
& =\ln 2+\lim_{\eps \rightarrow 0^+} \left( \ln \eps
+\int_{\Sigma_\eps^+} \ln \tau_\eps \, d\nu_\eps \right).
\end{split}
\end{equation*}
\par
Since \eqref{1.2} yields
\begin{equation*}
\lim\limits_{\eps \rightarrow 0^+} \left( \ln \int_{\Sigma_\eps^+}
\tau_\eps \, d\nu_\eps +\ln \eps +\ln 2\right)= \lim\limits_{\eps
\rightarrow 0^+} \left( \ln \int_{\Sigma_\eps^+} \tau_\eps \,
d\nu_\eps -\ln \frac{1}{2\eps} \right) =0,
\end{equation*}
we collect
\begin{equation*}
\lim\limits_{\eps \rightarrow 0^+} \left( \ln \int_{\Sigma_\eps^+}
\tau_\eps \, d\nu_\eps -\int_{\Sigma_\eps^+} \ln \tau_\eps \,
d\nu_\eps \right)=C.
\end{equation*}
\par
Finally we outline the proof of the identity
\begin{equation}\label{12.3}
C=3\ln 2-\frac{9\zeta(3)}{4\zeta(2)} \, .
\end{equation}
First we note that
\begin{equation*}
\int_0^1 g(t)\ln t\, dt=\frac{6}{\pi^2} \int_0^1 \ln t\,
dt=-\frac{6}{\pi^2}  ,
\end{equation*}
so we may write
\begin{equation}\label{12.4}
-C=-\frac{6}{\pi^2} +C_1+C_2+C_3 ,
\end{equation}
where
\begin{equation*}
\begin{split}
& C_1=\frac{6}{\pi^2} \int_1^\infty \left( \frac{2}{t} +2\Big(
1-\frac{1}{t}\Big)^2 \ln \Big( 1-\frac{1}{t} \Big) \right)
\ln t\, dt ,\\
& C_2 =\frac{6}{\pi^2} \int_1^2 \left( -\frac{1}{t}-\frac{1}{2}
\Big( 1-\frac{2}{t}\Big)^2 \ln \Big(
\frac{2}{t}-1\Big) \right) \ln t\, dt, \\
& C_3 =\frac{6}{\pi^2} \int_2^\infty \left(
-\frac{1}{t}-\frac{1}{2}\Big( 1-\frac{2}{t} \Big)^2 \ln \Big(
1-\frac{2}{t}\Big)\right) \ln t\, dt.
\end{split}
\end{equation*}
The substitution $t=2u$ leads to
\begin{equation*}
C_3=-\frac{1}{2} C_1 -\frac{6}{\pi^2}\, \ln 2 \int_1^\infty \left(
\frac{1}{u}+\Big( 1-\frac{1}{u}\Big)^2 \ln \Big(
1-\frac{1}{u}\Big) \right) du.
\end{equation*}
By a direct computation, the integral above is equal to $2(
\frac{\pi^2}{6}-1)$, thus
\begin{equation}\label{12.5}
C_3=-\frac{C_1}{2}-\left( \frac{6}{\pi^2} \ln 2\right) 2\left(
\frac{\pi^2}{6}-1\right) =-\frac{C_1}{2} -(2\ln 2) \left(
1-\frac{6}{\pi^2} \right) .
\end{equation}
Next, a direct computation shows that
\begin{equation}\label{12.6}
C_1 =\frac{12}{\pi^2} \big( 2\zeta(3)-1\big).
\end{equation}
The relations \eqref{12.4}--\eqref{12.6} provide
\begin{equation}\label{12.7}
-C=C_2+\frac{12}{\pi^2} \zeta(3)- \frac{12}{\pi^2} -2\left(
1-\frac{6}{\pi^2} \right) \ln 2.
\end{equation}
But
\begin{equation*}
C_2 =-\frac{6}{\pi^2}\cdot \frac{\ln^2 2}{2} -\frac{3}{\pi^2}
\int_1^2 \left( 1-\frac{2}{t}\right)^2 \ln \left(
\frac{2}{t}-1\right) \ln t\, dt,
\end{equation*}
thus we get
\begin{equation}\label{12.8}
-C=\frac{12}{\pi^2} \zeta(3)-\frac{12}{\pi^2} -2\left(
1-\frac{6}{\pi^2} \right) \ln 2- \frac{6}{\pi^2} \cdot \frac{\ln^2
2}{2} -\frac{3}{\pi^2} C_4,
\end{equation}
with
\begin{equation*}
C_4 =\int_1^2 \left( 1-\frac{2}{t}\right)^2 \ln \left(
\frac{2}{t}-1\right) \ln t\, dt =C_5-C_6+C_7 ,
\end{equation*}
where
\begin{equation*}
\begin{split}
& C_5 =\ln 2\int_1^2 \left( 1-\frac{2}{t}\right)^2 \ln t\, dt
,\quad C_6=\int_1^2 \left( 1-\frac{2}{t}\right)^2 \ln^2 t\, dt ,\\
& C_7=\int_1^2 \left( 1-\frac{2}{t} \right)^2 \ln \left(
1-\frac{t}{2} \right) \ln t\, dt.
\end{split}
\end{equation*}
\par
By a direct computation we find
\begin{equation*}
C_5=\ln 2-2\ln^3 2,\quad C_6=6-8\ln 2-\frac{4}{3}\ln^3 2.
\end{equation*}
As a result we gather
\begin{equation*}
C_4=\ln 2-2\ln^3 2-6+8\ln 2+ \frac{4}{3}\ln^3 2 +C_7 ,
\end{equation*}
and so
\begin{equation}\label{12.9}
\begin{split}
-C & =\frac{12}{\pi^2}\zeta(3)-\frac{12}{\pi^2} -2\ln 2
+\frac{12\ln 2}{\pi^2} -\frac{3}{\pi^2}\ln^2 2
-\frac{3}{\pi^2}\ln 2+\frac{6}{\pi^2} \ln^3 2    \\
& \qquad \qquad +\frac{18}{\pi^2} -\frac{24\ln 2}{\pi^2}
-\frac{4}{\pi^2} \ln^3 2 -\frac{3}{\pi^2} C_7 \\
& =\frac{12}{\pi^2} \zeta(3) +\frac{6}{\pi^2} -2\ln 2
-\frac{15}{\pi^2} \ln 2 -\frac{3}{\pi^2}  \ln^2 2 +\frac{2}{\pi^2}
\ln^3 2  -\frac{3}{\pi^2} C_7 .
\end{split}
\end{equation}
By a careful computation we find
\begin{equation*}
C_7 =-\ln^2 2 -5\ln 2+\frac{2\pi^2 \ln 2}{3}+4\mbox{\rm Li}_3
\left( \frac{1}{2}\right) -4\zeta (3)+2,
\end{equation*}
where $\mbox{\rm Li}_3$ denotes the trilogarithm function
\begin{equation*}
\mbox{\rm Li}_3 (z)=\sum\limits_{m=1}^\infty \frac{z^m}{m^3}
,\quad \vert z\vert \leq 1.
\end{equation*}
Using the equality (cf. \cite[formula (6.12)]{Lew})
\begin{equation*}
\mbox{\rm Li}_3 \left( \frac{1}{2}\right) =\frac{7}{8}\zeta(3)-
\frac{\pi^2}{12}\ln 2+\frac{\ln^3 2}{6}
\end{equation*}
we infer
\begin{equation*}
C_7=-\ln^2 2 -5\ln 2+\frac{\pi^2}{3}\ln 2-\frac{\zeta(3)}{2}
+\frac{2}{3} \ln^3 2 +2 .
\end{equation*}
Inserting this back into \eqref{12.9} we finally find
\begin{equation*}
\begin{split}
-C & =\frac{12}{\pi^2}\zeta(3)+\frac{6}{\pi^2} -2\ln 2
-\frac{15}{\pi^2}\ln 2 -\frac{3}{\pi^2}\ln^2 2
+\frac{2}{\pi^2} \ln^3 2 +\frac{3}{\pi^2}\ln^2 2  \\
& \qquad \qquad +\frac{15}{\pi^2}\ln 2-\ln 2+ \frac{3}{2\pi^2}
\zeta(3)-\frac{2}{\pi^2}\ln^3 2 -\frac{6}{\pi^2} \\ & =-3\ln
2+\frac{27}{2\pi^2}\zeta(3) =-3\ln 2+\frac{9\zeta(3)}{4\zeta(2)} .
\end{split}
\end{equation*}

\section*{Acknowledgments} We are grateful to Professor Giovanni Gallavotti for bringing to
our attention reference \cite{BGW} in 2002.

\end{document}